\documentclass[aop]{imsart}

\RequirePackage{amsthm,amsmath,amsfonts,amssymb}
\RequirePackage[numbers]{natbib}

\startlocaldefs
\numberwithin{equation}{section}

\theoremstyle{plain} 
\newtheorem{theorem}{Theorem}[section] 
\newtheorem{lemma}[theorem]{Lemma} 
\newtheorem{prop}[theorem] {Proposition} 
\newtheorem{cor}[theorem]  {Corollary} 

\theoremstyle{remark}
\newtheorem{rem}[theorem] {Remark} 

\usepackage[utf8]{inputenc}
\usepackage[dvipsnames]{xcolor}
\usepackage{dsfont}
\usepackage[colorlinks=true,linkcolor=blue,citecolor=blue]{hyperref}
\usepackage[mathscr]{eucal}
\usepackage{mathrsfs} 
\usepackage{xspace}
\usepackage{upgreek}
\usepackage{url}      
\usepackage{mathtools}
\usepackage{color}    
\usepackage{subfig}   
\usepackage{enumerate} 
\usepackage{centernot} 

\newcommand{\G}{{\rm NG}}
\newcommand{\Aold}{{\beta}}
\def\1{{\mathchoice {1\mskip-4mu\mathrm l}      
{1\mskip-4mu\mathrm l} 
{1\mskip-4.5mu\mathrm l} {1\mskip-5mu\mathrm l}}} 
\renewcommand{\d}{{\rm d}}
\newcommand{\R}     {\mathbb{R}} 

\newcommand{\N}     {\mathbb{N}} 
\renewcommand{\P}   {\mathbb{P}} 
\newcommand{\D}     {\mathbb{D}} 
\newcommand{\E}     {\mathbb{E}} 
\newcommand{\Q}     {\mathbb{Q}}

\newcommand{\smfrac}[2]{{\textstyle{\frac {#1}{#2}}}}

\newcommand{\supp}{\mathrm{supp}}
\newcommand{\Poi}{\mathrm{Poi}}
 
\newcommand{\ssup}[1] {{{\scriptscriptstyle{({#1}})}}} 

\newcommand{\eps}{\varepsilon} 
\newcommand{\Acal}  {{\mathcal A}}
\newcommand{\Bcal}  {{\mathcal B}}
\newcommand{\Ccal}   {{\mathcal C }}

\newcommand{\Jcal}   {{\mathcal J }}

\newcommand{\Mcal}   {{\mathcal M }} 
\newcommand{\Ncal}   {{\mathcal N }} 
 
\newcommand{\Pcal}   {{\mathcal P }}

\newcommand{\Scal}   {{\mathcal S }}

\newcommand{\Vcal}   {{\mathcal V }} 
 
\newcommand{\Xcal}   {{\mathcal X }} 
\newcommand{\Ycal}   {\ensuremath{\mathcal Y }\xspace} 
\newcommand{\Zcal}   {{\mathcal Z }} 
\newcommand{\e}   {{\operatorname e }}

\endlocaldefs

\begin{document}

\begin{frontmatter}
\title{Spatial particle processes with coagulation:
Gibbs-measure approach, gelation and Smoluchowski equation}
\runtitle{Spatial coagulation processes}

\begin{aug}
\author[A]{\fnms{Luisa}~\snm{Andreis}\ead[label=e1]{luisa.andreis@polimi.it}},
\author[B]{\fnms{Wolfgang}~\snm{K\"onig}\ead[label=e2]{wolfgang.koenig@wias-berlin.de}},
\author[C]{\fnms{Heide}~\snm{Langhammer}\ead[label=e3]{heide.langhammer@wias-berlin.de}}
\and
\author[C]{\fnms{Robert~I.~A.}~\snm{Patterson}\ead[label=e4]{robert.patterson@wias-berlin.de}}

\address[A]{Politecnico di  Milano\printead[presep={,\ }]{e1}}

\address[B]{TU Berlin and WIAS Berlin\printead[presep={,\ }]{e2}}

\address[C]{WIAS Berlin\printead[presep={,\ }]{e3,e4}}

\end{aug}

\begin{abstract}
We study a spatial Markovian {\em particle system with pairwise coagulation}, a spatial version of the Marcus--Lushnikov process:  according to a {\em coagulation kernel} $K$,
particle pairs merge into a single particle, and their masses are united. 
We introduce a {\em statistical-mechanics approach} to the study of this process. We derive an explicit formula for the empirical process of the particle configuration at a given fixed time $T$ 
in terms of a reference Poisson point process, whose points are  trajectories that coagulate into one particle by time $T$. The non-coagulation between any two of them induces an exponential pair-interaction, which turns the description into 
a {\em many-body system with a Gibbsian pair-interaction}. 

Based on this, we first give a {\em large-deviation principle} for the joint distribution of the particle histories (conditioning on an upper bound for particle sizes), in the limit as the number $N$ of initial atoms diverges and the kernel scales as $\frac 1N K$.
We characterise the minimiser(s) of the rate function, 
we give criteria for its uniqueness and prove a {\em law of large numbers} (unconditioned). Furthermore, we use the unique minimiser to construct a solution of the {\em Smoluchowski equation} and give a criterion for the occurrence of a {\em gelation phase transition}.
\end{abstract}



\end{frontmatter}
\tableofcontents



\section{Introduction}\label{Intro}

In this paper, we investigate a Markovian process of spatially distributed particles that coagulate pairwise. Our main interest is in asymptotics and phase transitions in the limit of many particles. To investigate this, we adopt a typical statistical mechanics approach: we describe each particle at time $T$ as a trajectory in space-time that lead to the formation of such particle (via the multiple coagulation events) and we write the joint distribution of all the particles  as a many-body system with a pair-interaction. Starting from this, we apply the tools from large deviation theory to study the limits and we use the results to derive other properties of the process. This section serves as an introduction to the work. In particular, in Section~\ref{sec-Goals} we summarise the goals of this paper. In Section~\ref{sec-model} we introduce the model, a spatial Markovian coagulation process. Then we turn in Section~\ref{sec-treedecomp} to an explanation of our approach which lies in a decomposition of the process into distinct parts, each leading to one of the single particles that can be observed at a fixed time.

\subsection{Summary of our results}\label{sec-Goals}

\subsubsection*{The model} We consider 
a spatial version of the {\em Marcus--Lushnikov model} for coagulating particles. 
Each particle has a mass (an integer) and a location (a point in a Polish space $\Scal$). 
We start with a {\em monodispersed} configuration, where each particle is an atom with mass one. Coagulations of pairs occur independently over the pairs after exponential random holding times, whose parameters depend only on the locations and masses of the two particles and are assumed symmetric in the two particles. The collection of all these parameters forms what is called the {\em coagulation kernel}, $K\colon(\Scal\times \N)^2\to[0,\infty)$.  The mass of the new particle is the sum of the two old ones, its location  is picked randomly according to a {\em placement kernel}, $\Upsilon$.
This leads to a continuous time Markov process whose dynamics are determined by $K$ and $\Upsilon$. 


\subsubsection*{Basic properties and the process $\Xi$} The total amount of mass in the system is kept constant, but all the particle sizes are non-decreasing, and the number of particles decreases by one at each coagulation event. After each coagulation event the parameters of the exponential random times are updated to accommodate the pairwise interaction between the newly formed particle and all other ones. This process is a special case of what is called a {\it cluster coagulation process} in~\cite{Nor00}. One can write the configuration of the process at time $t$ as a point measure $\Xi_t$, which registers the statistics of the particles according to their positions and masses. Hence, $\Xi_t$ lies in the set $\Mcal_{\N_0}(\Scal \times\N)$ of all point measures on $\Scal\times \N$, and the natural state space of the process $(\Xi_t)_{t\in[0,T]}$ is the set $\Gamma_T$ of measure valued trajectories $[0,T]\to \Mcal_{\N_0}(\Scal \times\N)$ that are piecewise constant and are such that in each jump two particles are lost and one (larger) particle is gained. 


\subsubsection*{Initial condition} In the present paper, we restrict to an initial distribution of atoms (i.e., the  monodispersed situation of single-atom particles), taken as a Poisson process with intensity measure $N\mu$ for some probability measure $\mu$, as the answers that we find are particularly transparent; the study of the model under deterministic initial atom configurations is deferred to future work.


\subsubsection*{The overall goals} First, we are interested in the distribution of $(\Xi_t)_{t\in[0,T]}$, where $T$ is fixed. Then we study the behaviour of $(\frac 1N\Xi_t)_{t\in[0,T]}$ in the limit of large total mass $N$ in the system, when the coagulation kernel $K$ is replaced by $\frac 1N K$. This is often called the {\it hydrodynamic limit}. Here, we would like to understand the limiting distribution of $(\frac 1N\Xi_t)_{t\in[0,T]}$ and the question of a {\em gelation phase transition}, i.e., the emergence of large particles (i.e., with diverging size depending on $N$) containing in total a macroscopic amount of atoms. Furthermore, we want to derive a characteristic  partial differential equation for the limit.


\subsubsection*{Our new approach} While all the contributions to the analysis of the Marcus--Lushnikov model that we are aware of use the generator of the process, martingale arguments and characteristic equations like the Smoluchowski and the Flory equation for finding answers (see the literature survey in Section~\ref{sec-literature}), our approach is fundamentally different: it follows patterns that are known from statistical physics, large-deviations analysis and variational calculus. 


\subsubsection*{Our contributions} Our main contributions are the following:
\begin{enumerate}
\item[(1)] We derive an explicit formula for the distribution of $(\Xi_t)_{t\in[0,T]}$  in terms of an {\em interacting Poisson point process}, whose points are the histories (trajectories) of each particle present at time $T$.

\item [(2)]We prove a {\em large-deviation principle} for $ ( \frac 1N\Xi_t)_{t\in[0,T]}$ in the limit $N\to\infty$ of a large total mass $\asymp N$ in the system (for $K$ replaced by $\frac 1N K$).

\item[(3)] We give criteria for the occurrence or non-occurrence of large particles with macroscopic total atom mass  (i.e., a {\em gelation phase transition}) using properties of the large-deviation rate function under bounds on the coagulation kernel. In particular, we give criteria for the convergence of $ (\frac 1N \Xi_t)_{t\in[0,T]}$ as $N\to\infty$.

\item[(4)] As another by-product of our large-deviation analysis, we derive that in the subcritical regime, i.e., when there is no gelation, the process $(\frac 1N \Xi_t)_{t\in[0,T]}$ converges as  $N\to\infty$ to a solution of a {\em spatial version of the famous Smoluchowski equation}.

\end{enumerate}


\subsubsection*{Coagulation trajectories} Our ansatz for our goal (1) and for everything that follows is a decomposition of the entire configuration process on the time interval $[0,T]$ into the parts of the process (called {\em coagulation trajectories} or {\em history trees}) that coagulate by time $T$ into one of the particles that we see at that time. The main object of the  present paper is  the {\em normalised empirical measure}, $\Vcal_N^{\ssup T}$, of all these trajectories. This object is much more comprehensive and detailed than just the process $(\Xi_t)_{t\in[0,T]}$. All the randomness of the process (holding times and placement decisions) is attached to these one-particle trajectories (coagulation events that lead to the particle) and their mutual interaction (the non-coagulation between the trajectories). 


\subsubsection*{A crucial Poisson process} Indeed, in our first main result, Theorem~\ref{thm-distMi}, we introduce a Poisson point process whose points are in the space of coagulation trajectories. Its intensity measure describes the coagulation decisions within each trajectory and the family of these trajectories are independent and identically distributed. However, the non-coagulation between them is expressed in terms of a pairwise interaction between each pair of  trajectories. In other words, we represent the distribution of the configuration as an expectation over many independent coagulation trajectories with an exponential interaction term expressing the non-coagulations. This gives our representation the structure of a many-body system (better: a many-trajectory system): a {\em Gibbsian ensemble} of many independent trajectories with exponential pair interaction. This turns the Markovian coagulation process into  a static model of statistical mechanics with the underlying reference measure as the law of a  {\em Poisson point process} (PPP) on the set of coagulation trajectories. 


\subsubsection*{Other Poisson-process approaches} This Poisson-process description opens up a multitude of further research directions, in particular a comprehensive analysis of the entire trajectory configuration on the time-interval $[0,T]$.
It also makes the model amenable to an asymptotic analysis like hydro- or thermodynamic limits. Related Poisson representations have been carried out for a few models in statistical mechanics, like the interacting Bose gas \cite{ACK11, CJK23} and the Erd\H{o}s--R\'enyi graph with and without spatial component \cite{AKP21,AKLP23}, see also Section~\ref{sec-ERgraph} for an account on the latter subject.
Also in a recent paper \cite{Sun23}, the joint distribution of all the components of the Erd\H{o}s--R\'enyi graph has been identified in terms of a Poisson point process (here a compound one). It is used there to derive moderate-deviations results for three variables (size of largest component, number of components of a given fixed size and total number of components).


\subsubsection*{Large-deviation principle} Let us turn now to our goal (2). We use our Poissonian representation for deriving asymptotics of the empirical process $\Vcal_N^{\ssup T}$ for diverging system size $N$ and rescaled coagulation kernel $\frac 1N K$ with the help of {\em large-deviations analysis}, one of the ubiquitous approaches to Gibbsian systems. In our case a useful large-deviations principle (LDP) for the reference Poisson process is known, and it has an explicit rate function. Since the state space of $\Vcal_N^{\ssup T}$  is huge, we need to employ a conditioning on an event that generates compactness and continuity properties of certain functionals. This conditioning basically excludes the occurrence of large particles, which rules out the possibility to observe the gelation phase transition right away. In this conditional setting, we obtain in our second main result, Theorem~\ref{thm-LDP}, a full LDP for $\Vcal_N^{\ssup T}$ and an explicit formula for the rate function of the coagulation process. 


\subsubsection*{Consequences: Euler--Lagrange equation and  convergence} This LDP and its rate function now allows for a deeper analysis. We prove that every minimiser of the rate function (i.e., every possible accumulation point of $(\Vcal_N^{\ssup T})_{N\in\N}$ under the conditional measure) satisfies the Euler--Lagrange equation, and formulate assumptions under which the solutions in turn are unique, employing an argument known from Banach's fixed point theorem.
One obviously wishes to remove the conditioning, but, unfortunately, it cannot be removed on an exponential scale. Nevertheless, we succeed in formulating assumptions under which we can prove that the probability of the conditional event converges to one, and we derive tightness of the unconditioned distribution of $\Vcal_N^{\ssup T}$ and the convergence $\Vcal_N^{\ssup T}$ towards the minimiser of the large-deviation rate function. This is our third main result, Theorem~\ref{thm-Lossofmass}(1), together with Proposition~\ref{prop-criteria}(1). 


\subsubsection*{Our assumptions} Remarkably, all the preceding is obtained, for all sufficiently small $T$,
under the sole condition that the coagulation kernel $K$ is continuous and satisfies the bound in \eqref{AssK1}, that is,
$$
H:=\sup_{v,w\in\Mcal(\Scal\times \N)\colon \|v\|_1,\|w\|_1\leq 1}\langle v, Kw\rangle <\infty.
$$
(We will provide the necessary notation to understand the supremum in Section \ref{sec-results}.)
More precisely, the results about convergence of $\Vcal_N^\ssup{T}$ hold if $TH<\frac 1{\e^2}\frac\pi{1+\pi}$.
The criterion $H<\infty$ is in spirit of upper bounding the kernel against a product kernel of the form  $K((x,m),(x',m'))=\varphi(x,x') m m'$ with bounded $\varphi$, where $x,x'$ are the locations and $m,m'$ the masses of the two particles, but is much more general.  


\subsubsection*{Gelation and non-gelation} Turning to our goal (3), the macroscopic occurrence or non-occurrence of large particles, we also find that, under the above condition, there is no formation of a gel for all $TH<\frac 1{\e^2}\frac\pi{1+\pi}$. On the other hand, if additionally the \lq dual\rq\ condition
$$
\inf_{v,w\in\Mcal(\Scal\times \N)\colon \|v\|_1=\|w\|_1= 1}\langle v, Kw\rangle >0
$$
is satisfied (see \eqref{AssK2}), then we prove in Theorem~\ref{thm-Lossofmass}(2), jointly with Proposition~\ref{prop-criteria}(2), the occurrence of a gel for all sufficiently large $T$. 
We define gelation via the existence of large particles (with a diverging size that can be on any scale) that build a macroscopic part of the configuration (the gel), but do not specify the scale of the size of the large particles. The existence of such particles is identified by the fact that, when considering the total mass of atoms sitting in all the particles at time $T$ of sizes $\leq L$, after letting $N\to\infty$, followed by $L\to\infty$, we see that the mass is strictly less than the total mass at time $0$.
In future work, we plan to analyse also the macroscopic and mesoscopic part of the configuration to obtain a more detailed understanding of the gelation phase transition and the gel itself.


\subsubsection*{Results for the ML-process} So far, all our results are formulated in terms of the empirical trajectory measure, $\Vcal_N^{\ssup T}$, but one is also highly interested in the Marcus--Lushnikov process $\Xi$ itself. While  the process $(\Xi_t)_{t\in[0,T]}$ takes values in the set $\Gamma_T$ of trajectories $[0,T]\to  \Mcal_{\N_0}(\Scal\times\N)$, the empirical process $\Vcal_N^{\ssup T}$ is contained in a much more comprehensive and more abstract space.
Nevertheless, $(\frac 1N \Xi_t)_{t\in[0,T]}$ is a relatively simple functional of $\Vcal_N^{\ssup T}$, which turns out in Lemma~\ref{lem-Contrho} to be continuous in a sufficient sense, such that a great deal of our results about $\Vcal_N^{\ssup T}$ have a consequence for  $(\frac 1N \Xi_t)_{t\in[0,T]}$. 


\subsubsection*{LDP for the ML-process} Our first observation is in Corollary~\ref{cor-LDPXi} that, via the well-known contraction principle, $(\frac 1N \Xi_t)_{t\in[0,T]}$ also satisfies an LDP (under the same conditioning as $\Vcal_N^{\ssup T}$) with an explicit rate function. This nice fact opens up a completely new path to LDPs for the trajectories of the Marcus--Lushnikov process, which is totally different from and independent of the much-used  Freidlin--Wentzell theory that uses a toolbox from operator theory and stochastic processes (see the literature survey in Section~\ref{sec-literature}). It appears to have the great advantage to be successful also in the current setting of a pretty abstract and huge state space, in contrast with existing works using the mentioned toolbox.


\subsubsection*{The Smoluchowski equation} Furthermore, and now we are turning to our goal (4) of deriving a partial differential equation for the limiting objects, 
by the continuity of the map that maps $\Vcal_N^{\ssup T}$ onto $(\frac 1N \Xi_t)_{t\in[0,T]}$, 
we also obtain, under the mentioned assumptions, the convergence of the latter to a deterministic process, $\rho=(\rho_t)_{t\in[0,T]}$. From the Euler--Lagrange  equation for the minimisers of our rate function for $(\Vcal_N^{\ssup T})_{N\in \N}$, we derive in Lemma~\ref{lem-Smol} that $\rho$ is a solution to a natural spatial version of the famous Smoluchowski equation on $[0,T]$. In our approach, this equation is only a by-product and provides a nice additional information, but is not a vital part of the proof of convergence of $(\frac 1N \Xi_t)_{t\in[0,T]}$, as it is in many previous investigations of related processes.


\subsubsection*{Advantages and disadvantages} Let us briefly discuss advantages and disadvantages of our ansatz in comparison to earlier ansatzes; see Section \ref{sec-literature}. Convergence of the coagulation process was previously  established, under more restrictive conditions on the kernel $K$, using an entirely different route: tightness arguments that give convergent subsequences whose limits are solutions of the Smoluchowski equation that are combined with uniqueness arguments to identify the deterministic limit.  Previous criteria for convergence of the Marcus--Lushnikov process require $K$ to be conservative (or approximately conservative in a pretty narrow sense), i.e., to keep a certain quantity constant in each coagulation event.
The most natural example of this is a product kernel of the form $K((x,m),(x',m'))=\varphi(x,x') m m'$, with $\varphi$ a bounded function; this kernel is such that the sum of all the rates is preserved by a coagulation event. 
The good news is that our sufficient criterion for convergence is not restricted to such kernels, but consists just of the general upper bound \eqref{AssK1}, which does not require any information about the form of the kernel. The bad news is that our approach {\it a priori} does not work up to the gelation time (the first time at which large particles with macroscopic total mass arises), but only on a time interval $[0,T]$ with $T$ small enough. 


\subsubsection*{A uniqueness criterion} Let us finally mention, see Remark~\ref{Rem-convexity}, that our approach gives also uniqueness of the minimiser of our rate function under the condition that the kernel $K$ is nonnegative definite; however we believe that this criterion is not too helpful, since it is in general difficult to check this property.


\subsection{The spatial Marcus--Lushnikov process}\label{sec-model}

Let us enter into the details of our model. As mentioned, particles  live in a Polish  space $ \Scal$, i.e. a separable, complete metric space. Each particle carries a mass $m\in\N$ and sits at a site $x\in \Scal$; we  sometimes also say then that the particle sits at $(x,m)$. Initially, each particle has mass one, i.e., we consider what is usually called a monodisperse initial condition. The  unit-mass particles at time zero are also called {\em atoms}, since all later particles are composed out of them by merging and since they are never split anymore into smaller units. The particle process is a Markov process in continuous time. The dynamics of the process depends on a {\em coagulation kernel} and a {\em placement kernel},
\begin{equation}\label{kernels}
K\colon (\Scal\times\N)^2\to[0,\infty)\qquad\mbox{and}\qquad \Upsilon\colon (\Scal\times\N)^2\times \Bcal(\Scal)\to[0,1],
\end{equation}
where $\Bcal(\Scal)$ denotes the Borel-$\sigma$-field on $\Scal$.
We assume that $K$ is symmetric and measurable, and that $\Upsilon$ is a Markov kernel from $(\Scal\times\N)^2$ into $\Scal$, i.e., measurable in the argument $(\Scal\times\N)^2$ and a probability measure in the last argument. We also assume that $\Upsilon$ is symmetric in the first two arguments, i.e., $\Upsilon((x,m),(x',m'),\cdot)=\Upsilon((x',m'),(x,m),\cdot)$.

\medskip

{\bf The random mechanism:}  {\em At each time $t\in[0,\infty)$, for each unordered pair of particles in the current configuration, located at $x$ and $x'\in\Scal$ with masses $m$ and $m'$, respectively, there is an exponential random time with parameter $K((x,m),(x',m'))$ running. When it elapses, the pair is replaced by a single particle with mass $m+m'$, located at a random site that is picked according to $\Upsilon((x,m),(x',m'),\cdot)$. The exponential random times and the placement locations are independent over all particle pairs and over all times.}

\medskip

Hence, after each coagulation event, the parameters of all exponential random times involving any of the two coagulating particles are updated. Since we are starting with only finitely many particles, this excludes explosion and the total number of coagulation events during a time interval $[0,T]$ is finite. Hence, the entire process can be decomposed into finitely many time intervals during which the configuration remains constant. (The process that we study is the special case of what is called a {\it cluster coagulation process}  in \cite{Nor00}, with state space $E=\Scal\times\N $ and mass-preserving function $m((x,m))=m$.)

The above mechanism defines a Markov chain in continuous time
\begin{equation}\label{Xidef}
\Xi=(\Xi_t)_{t\in[0,\infty)},\qquad\mbox{where}\qquad\Xi_t=\sum_{i\in[n(t)]}\delta_{(X_i(t),M_i(t))}\in\Mcal_{\N_0}(\Scal\times\N),
\end{equation}
where $(X_i(t),M_i(t))\in\Scal\times\N$ is the location and the mass of the $i$-th particle at time $t$, and $n(t)$ is the number of particles at time $t$ (we abbreviate $[n]=\{1,\dots,n\}$). We denote by  $\Mcal_{\N_0}$ the set of  measures with values in $\N_0$, i.e., the set of finite point measures. We will treat elements $\phi\in \Mcal_{\N_0}(\Scal \times \N)$ as partially discrete measures and we will write for example $\phi(x,m)$ instead of $\phi(\{(x,m)\})$. We write also $\Xi_t(A,m)$ for the number of particles in $A\subset\Scal$  with mass $m\in\N$ at time $t$ in the configuration.

The process $(n(t))_{t\in[0,\infty)}$ is non-increasing in $t$; it actually decreases by one at every coagulation time. We include also the trivial case $n(t)=0$ for all $t$, in which case there is no atom and hence the process is empty. As usual for point processes, the index $i\in[n(t)]$ is arbitrary and does not specify the $i$-th particle. In particular, if there are multiple particles with the same location and mass, then $\Xi$ does not give information about which one of them are involved in a coagulation. 

The process $\Xi$ is a Markov process with jumps of type
\begin{equation}\label{steps}
\phi\mapsto \phi-\delta_{(x,m)}-\delta_{(x',m')}+\delta_{(z,m+m')},\qquad \phi\in\Mcal_{\N_0}(\Scal\times\N),
\end{equation}
(as long as the right-hand side is nonnegative) that happen with rate
\begin{equation}\label{Mkerneldef}
\begin{aligned}
{\bf K}_\phi\big((x,m),(x',m'), \d z\big) &= {\mathbf K}\big((x,m),(x',m'), \d z\big)\\
&\quad \times
\begin{cases}\phi(x,m)\phi(x',m'),&\mbox{if }(x,m)\not=(x',m'),\\
\phi(x,m)(\phi(x,m)-1)/2&\mbox{otherwise,}
\end{cases}
\end{aligned}
\end{equation}
and we abbreviated
\begin{equation}\label{z-kernels}
{\mathbf K}\big((x,m),(x',m'),\d z\big)=K\big((x,m),(x',m')\big)\, \Upsilon\big((x,m),(x',m'),\d z\big).
\end{equation}
Since $\Upsilon((x,m),(x',m'),\cdot)$ is a probability measure, the rate $K((x,m),(x',m'))$ is equal to the total mass of the measure ${\mathbf K}((x,m),(x',m'),\cdot)$. Analogously, we put \[K_\phi((x,m),(x',m'))={\bf K}_\phi((x,m),(x',m'), \Scal).\] Sometimes we call ${\mathbf K}$ also a rate, even if this makes strict sense only if ${\mathbf K}((x,m),(x',m'),\cdot)$ has a discrete support. The counting factor in the second line of \eqref{Mkerneldef} is the number of unordered  pairs of particles of types $(x,m)$ and $(x',m')$. In the case that $\Scal$ is uncountable, formally there are uncountably many $(x,x')$ that give rise to a step, but actually only finitely many have a positive rate, since $\phi$ has a finite support.  

We always start with a monodispersed situation, i.e., we fix $n(0)\in \N$, $M_i(0)=1$ for any $i\in [n(0)]$, and some configuration of $(X_i(0))_{i\in [n(0)]}$. Hence, $\Xi_0=\sum_{i\in[n(0)]}\delta_{(X_i(0),1)}$, and we often identify with  the configuration $\sum_{i\in[n(0)]}\delta_{(X_i(0))}\in\Mcal_{\N_0}(\Scal)$. The trajectories of $\Xi$ lie in the set
\begin{equation}\label{Gammadef}
\begin{aligned}
\Gamma&=\Big\{\xi=(\xi_t)_{t\in[0,\infty)}\in \D\big(\Mcal_{\N_0}(\Scal\times\N)\big)\colon  \xi_0 \text{ is concentrated on } \Mcal_{\N_0}(\Scal\times \{1\})\\
& \hspace{3cm} t\mapsto \xi_t \mbox{ is piecewise constant and makes steps as in \eqref{steps}} \Big\},
\end{aligned}
\end{equation}
where $\D(\mathcal X)$ denotes the set of càdlàg-paths $[0,\infty)\to \mathcal X$, i.e., paths that are right-continuous in $[0,\infty)$ and have left limits everywhere in $(0,\infty)$ (we give $\Mcal_{\N_0}(\Scal\times\N) $ the weak topology induced by integrals against continuous bounded test functions). In particular, the initial configuration of any $\xi\in\Gamma$ is finite, and the collected mass of all  particles at time $t$, which we indicate with $\|\xi_t\|_1=\sum_{m\in\N} m\,\xi_t(\Scal\times\{m\})$,  is a  constant in $t$.

For any $k\in \Mcal_{\N_0}(\Scal)$, we write $\P_k$ and $\E_k$ for probability and expectation with respect to the process $\Xi$ when started from $\Xi_0=k$. For $k=0$ we define $ \P_0$ as the measure that is concentrated on the constant zero point measure. 

\medskip

\subsubsection*{Initial configurations}
In this paper, we will formulate and prove all our results for the distribution, $\P_{\Poi_{N\mu}}$, of the process under poissonised initial conditions so that
\begin{equation}\label{PoissonInitial}
\P_{\Poi_{N\mu}}=\int_{\Mcal_{\N_0}(\Scal)}\Poi_{N\mu}(\d k)\, \P_k
\end{equation}
where $\Poi_{N\mu}$ is the law of a Poisson point process (PPP)  on $\Scal$ with intensity measure $N\mu$ for an arbitrary probability measure $\mu$ on $\Scal$.
Note that the number $n(0)$ of initial particles is not deterministic, but is $\Poi_{N}$-distributed and the empty configuration appears with  probability $\e^{-N}>0$. 
The analogous statements and proofs for the deterministic initial configuration $\P_{k_N}$ with $k_N\in\Mcal_{\N_0}(\Scal)$ satisfying $\frac 1N k_N  \to \mu$ for some $\mu\in\Mcal_1(\Scal)$ are deferred to future work.

\begin{rem}[Special choices]\label{Rem-placement}
A natural choice of $\Upsilon$ is the deterministic choice $z=\frac {xm+x'm'}{m+m'}$ (if $\Scal$ is convex), which keeps the centre of mass of the two particles (and hence also the centre of mass of the entire configuration) constant. Another one is the random one, where $z$ is put equal to $x$ or $x'$ with probability $\frac m{m+m'}$ and $\frac {m'}{m+m'}$ respectively; this choice keeps the centre of mass fixed on average.\hfill$\Diamond$
\end{rem}

\begin{rem}[Inhomogeneous Erd\H{o}s--R\'enyi graph]\label{Rem-ERgraph}
An important special case in the non-spatial setting (i.e., $\Scal$ is a singleton) is the product kernel $K(m,m')=H m m'$ for some $H\in(0,\infty)$. In this special case, the coagulation process can be mapped one-to-one onto the process of the family of connected component sizes of the well-known Erd\H{o}s--R\'enyi graph in its dynamic version.
 Indeed, this process can be seen as a simplified version of a coagulation process:  instead of a coagulation event of two particles, only a bond between two atoms in the two connected components is added, see Section~\ref{sec-ERgraph}. There we also explain the case of an {\em inhomogeneous} (i.e., spatial) version of this process and its relation to a spatial coagulation process.\hfill$\Diamond$
\end{rem}

\subsection{Decomposition into coagulation trajectories}\label{sec-treedecomp}

In this section, we will describe the distribution of $(\Xi_t)_{t\in[0,T]}$ for a fixed $T>0$ via a decomposition into subprocesses that coagulate into one particle by time $T$, which we will call  {\em coagulation trajectories} on the interval $[0,T]$. This section provides necessary notation.

Observe that the notation in \eqref{Xidef} only counts particles at a given site with a given size. It is not rich enough to differentiate between multiple particles that sit on the same site and can therefore not express information about which particles coagulate into which particle. As a consequence, it cannot grasp the full history of a particle. In order to express the evolution explicitly we need to introduce an alternative version of the coagulation process, that assigns a label to every atom at time $0$ such that we can follow its path through the coagulation process. To this end, we define a  Markov process $(P_t)_{t\in[0,\infty)}$ on the set of all partitions  of $[n(0)]$ starting from $P_0=\{\{i\}\colon i\in [n(0)]\}$, together with the locations of the particles (elements of the partition). 

As a preparation we need some notation. For any finite, non-empty set $A$ we denote the set of partitions of $A$ by
\begin{equation}\label{partion1}
\Pcal(A)=\Big\{\{C_j\}_{j\in J}\colon J\mbox{ an index set}, \dot\bigcup_{j\in J}C_j=A,\forall j\colon \emptyset\not= C_j\subset A\Big\}.
\end{equation}

We define a process $(P_t)_{t\in [0,\infty)}$ on $\Pcal(A)$ in such a way that $t\mapsto P_t$ makes discrete steps by joining two of the partition sets and is otherwise constant; in particular, $t\mapsto |P_t|$ decreases by one in each step.  We refer to each partition set $C\in P_t$ as a particle that exists at time $t$ and attach to it the site $X^\ssup{t}_C\in\Scal$ at which the particle sits. More precisely, we define the {\em labelled coagulation process} 
\begin{equation}\label{Zprocess}
Z=(Z_t)_{t\in[0,\infty)},\qquad \mbox{with }Z_t=(X_C^{\ssup t},C)_{C\in P_t},
\end{equation}
as a Markov process with the mechanism
\begin{equation}\label{mechanism1}
\big((X_C,C),(X_D,D)\big)\mapsto (X_{C\cup D},C\cup D)\mbox{ with  rate } {\mathbf K}((X_C,|C|),(X_D,|D|),\d  X_{C\cup D}),
\end{equation}  
where our short-hand notation $\big((X_C,C),(X_D,D)\big)\mapsto (X_{C\cup D},C\cup D)$ includes the fact that for all $\widetilde C \notin \{C,D\}$ 
the values $(X_{\widetilde{C}}, \widetilde{C})$ are unchanged in the transition. In contrast to $\Xi$, the process $Z$ contains for every particle $C$ present at time $t$ the information from which of the atoms at time $0$ it stems.
Given $\mathbf{x} = (x_i)_{i\in A} \in \Scal^A$ we denote by $\P_{\mathbf{x}}$ the distribution of this Markov process starting from $Z_0=(X_C^{\ssup 0},C)_{C\in P_0}=(x_i,\{i\})_{i\in A}$. We call $\mathbf x$ the {\it initial atom configuration} of the process.

Let us explain the connection between the process $Z$ and $\Xi$. Indeed, given $Z$, we can define 
\begin{equation}\label{XiZ}
\Xi(Z) = (\Xi_t(Z))_{t\in [0,\infty)}\qquad \text{ with } \Xi_t(Z) = \sum_{C\in P_t}\delta_{(X_C^{\ssup t},|C|)}.
\end{equation}
Assume that we fix $k\in \Mcal_{\N_0}(\Scal)$ and some vector $\mathbf{x}= (x_i)_{i=1}^{|k|}\in \Scal^A$ that is {\em compatible} with $k$, i.e. it satisfies $k=\sum_{i=1}^{|k|}\delta_{x_i}$. Then,
\begin{equation}\label{ZversusXi}
\P_k\circ\Xi^{-1}= \P_{\mathbf x}\circ \big( \Xi(Z)\big)^{-1}
\end{equation}
where we use the measure-theoretic notation $\mu\circ X^{-1}$ for the distribution of a random variable $X$ under a probability measure $\mu$.
In this way, we actively forget the precise history of each of the particles and just count them at any time with regard to their locations and sizes.

Now we fix a time $T\in(0,\infty)$. For each particle present at time $T$ we want to define a subprocess that describes the history of the particle. 
Now, let $(Z_t)_{t\in [0,T]}$ with $Z_t = (X_C^\ssup{t},C)_{C\in P_t}$ be the process defined in \eqref{Zprocess} on the time interval $[0,T]$. Fix a set $C$ such that either $C \in P_T$ (which is the case we consider most of the time) or such that $C$ can be written as a union of sets in $P_T$. Define
\begin{equation}\label{history}
\Xi^{\ssup {T,C}}=(\Xi^{\ssup C}_t)_{t\in[0,T]},\qquad\mbox{where}\qquad \Xi^{\ssup{T,C}}_t=\sum_{\widetilde C \in P_t\colon \widetilde C\subset C}\delta_{(X_{\widetilde C}^{\ssup t},|\widetilde C|)}\in\Mcal_{\N_0}(\Scal\times\N).
\end{equation}

\begin{figure}[!htpb]
	\centering
\includegraphics[scale=0.6]{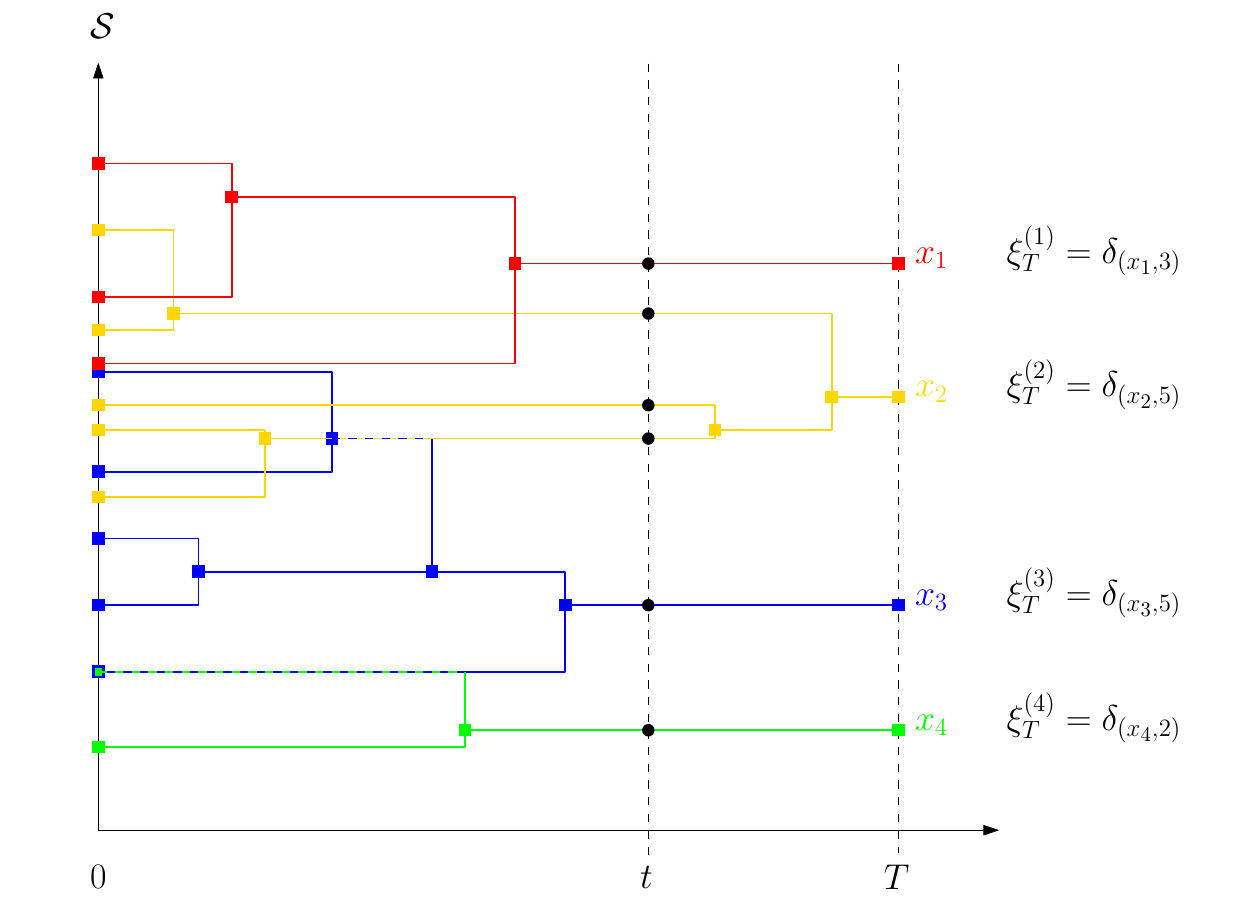}
	\caption{An illustration of the decomposition of $(\Xi_t)_{t\in [0,T]}$ into four subprocesses $(\xi^\ssup{i}_t)_{t\in[0,T]}$, $i=1, \ldots, 4$, that are distinguished by their colour. The process $(\Xi_t)_{t\in [0,T]}$ is started from $15$ atoms at time 0 and ends up with $4$ particles. However, each subprocess $(\xi^\ssup{i}_t)_{t\in[0,T]}$, $i=1, \ldots, 4$, ends in a single-particle configuration by definition. At time $t$, there are $6$ particles with masses $3, 2, 1, 2, 5$ and $2$. }
	\label{historytrees}
\end{figure}

One can easily check, that as a consequence of the mechanism \eqref{mechanism1} the sum on $\widetilde C\in P_t$ satisfying $\widetilde C \subset C$ is non-empty. $\Xi^\ssup{T,C}$ is the (sub-)process that only keeps track of the atoms with labels in $C$, i.e., $\Xi^\ssup{T,C}_t$ is the numer of particles with a given location and mass at time $t$ that have emerged from atoms with labels in $C$.  In particular, it holds that $\Xi^{\ssup {T,C}}_0=\sum_{i\in C}\delta_{(X_{\{i\}}^{\ssup 0},1)}$.
In the following we mainly choose $C\in P_T$ (i.e., $C$ is one particle at time $T$) and in that case we note that $\Xi^\ssup{T,C}$ is an element of the {\em set of (one-particle) coagulation trajectories},
\begin{equation}\label{Gamma1def}
\Gamma_T^{\ssup 1}=\{\xi\in\Gamma_T \colon \xi_T(\Scal\times \N)=1\},
\end{equation}
where 
\begin{equation}\label{GammaTdef}
\Gamma_T=\pi_{[0,T]}^{-1}(\Gamma), \qquad\mbox{where }\pi_{[0,T]}(\xi)=(\xi_t)_{t\in[0,T]},
\end{equation}
is the set of trajectories on the time interval $[0,T]$. 
Thus, if $C\in P_T$, then $\Xi^\ssup{T,C}$ tracks history of the particle $(X_C^{\ssup T},|C|)$, ending in a single-particle configuration $\delta_{(X_C^{\ssup T},|C|)}$ at time $T$.
In Figure~\ref{historytrees} we illustrate how the process $(\Xi_t)_{t\in [0,T]}$ decomposes into  the subprocesses $\Xi^\ssup{T,C}$, $C\in P_T$. 

Now we introduce the main object of this paper: the normalised empirical measure of  the coagulation trajectories of all the particles present at time $T$:
\begin{equation}\label{empmeas}
\Vcal_{N}^{\ssup T}=\frac{1}{N}\sum_{C\in P_T}\delta_{\Xi^{\ssup {T,C}}}\in\Mcal(\Gamma_T^{\ssup 1}).
\end{equation}

For the definition of $\Vcal_N^{\ssup T}$, we need the process $(Z_t)_{t\in [0,T]}$ and need to work {\it a priori} on the extended probability space with measure $\P_{\mathbf x}$. 
However, observe that $\Vcal_{N}^{\ssup T}$ does not depend on the labels, even though its definition  uses the labels in order to decompose the process into its subprocesses $\Xi^\ssup{T,C}$, $C\in P_T$. More precisely, if we fix $k\in \Mcal_{\N_0}(\Scal)$ 
and consider two vectors $\mathbf{x}= (x_i)_{i=1}^{|k|}$, $ \mathbf{\widetilde x}= ( \widetilde x_i)_{i=1}^{|k|}$ with $\sum_{i=1}^{|k|} \delta_{x_i} = \sum_{i=1}^{|k|} \delta_{\widetilde x_i} =k$, 
then ${\P_{\mathbf{x}}(\Vcal_N^\ssup{T}\in \cdot ) = \P_{ \mathbf{\widetilde x}}(\Vcal_N^\ssup{T}\in \cdot )}$. Hence, in a small abuse of notation, we will consider $\Vcal_{N}^{\ssup T}$ under the measure $\P_k$ and mean its distribution under $\P_{ \mathbf{\widetilde x}}$ for some $\mathbf{\widetilde x}$ that is compatible with $k$. Hence, we can assume that $\Xi$ and $\Vcal_{N}^{\ssup T}$ are defined on the same probability space with measure $\P_k$. In particular, we will speak of its distribution under the Poisson initial configuration, $\P_{\Poi_ {N\mu}}$, and mean the distribution with an initial configuration vector whose length is $\Poi_N$-distributed and has i.i.d.~$\mu$-distributed entries.

We decided to drop in $\Vcal_N^\ssup{T}$ the particle-by-particle history and to keep only the point-process information about the trajectories, but this is a matter of taste. However, some of our proofs will later need to work in the more detailed setting.

%
%

Let us explain how we can obtain the process $(\Xi_t(Z))_{t\in [0,T]}$ from $\Vcal_N^{\ssup T}$.
Let us denote the projection $\pi_t\colon \Gamma\to \Mcal_{\N_0}(\Scal \times \N)$ by $\pi_t(\xi)=\xi_t$, 
and let us write $\nu_t=\nu\circ \pi_t^{-1}$ for $\nu\in\Mcal(\Gamma)$ and $t\in[0,T]$. Then
\begin{equation}\label{Xi_as_function_of_nu}
\begin{aligned}
\frac 1N \Xi_t&=\frac 1N \sum_{\widetilde C\in P_t}\delta_{(X_{\widetilde C}^{\ssup t},|\widetilde C|)}
= \frac 1N \sum_{C\in P_T}\sum_{\widetilde C\in P_t\colon \widetilde C\subset C} \delta_{(X_{\widetilde C}^{\ssup t},|\widetilde C|)}\\
&
=\frac 1N \sum_{ C\in P_T} \Xi^\ssup{T, C}_t=
\int \Vcal_N^{\ssup T}(\d \xi)\, \xi_t.
\end{aligned}
\end{equation} 
This connection is crucial, since it shows that $(\frac 1N \Xi_t)_{t\in[0,T]}$ is a function of $\Vcal_N^\ssup{T}$ and allows us to understand the dynamics of $(\frac 1N \Xi_t)_{t\in[0,T]}$ by studying $\Vcal_N^\ssup{T}$.



The organisation of the remainder of this paper is as follows.
Section~\ref{sec-results} presents the main results of this paper, which we discuss and compare to the literature in Section~\ref{sec-discussion}. 
In Section~\ref{MiMaintroduction} we prove the first main result (identification of the distribution), in Section~\ref{sec-Prep} we prepare for the proofs of the others, which we complete in Sections~\ref{sec-LDPproof} (large-deviation principle) and  Section~\ref{sec-minimizerproof} (convergence and gelation criteria, and the Smoluchowski equation).

\section{Results}\label{sec-results}

In this section we formulate and comment our four main results: a representation of the empirical measure of coagulation trajectories in Section~\ref{sec-identdist}, a large-deviation principle for this empirical measure in Section~\ref{sec-LDP}, criteria for the occurrence of the gelation phase transition and convergence in Section~\ref{sec-minimizer}, and the validity of the Smoluchowski equation for the limit in Section~\ref{sec-Smol}.

First, let us introduce some notation that we use throughout the paper. For any measure $\mathfrak{m}$ on some measure space $\Xcal$ and any measurable function $f$ on $\Xcal$ we write $\langle \mathfrak{m}, f\rangle$ for the integral of $f$ with respect to the measure $\mathfrak{m}$. We write $|\mathfrak{m}| = \mathfrak{m}(\Xcal)$ for the total mass of $\mathfrak{m}$. We denote by $\Mcal_1(\Xcal)$ the set of probability measures on $\Xcal$.  For measures $v\in \Mcal(\Scal\times \N)$ we use the fact that $v$ is partially discrete and write $v(\d x,m)$ instead of $v(\d(x,m))$ and define $\|v\|_1 = \sum_{m\in \N}\int_\Scal v(\d x,m) m$. 

Considering Poisson point processes we use the following notation. For $\gamma\in(0,\infty)$ we denote with $\Poi_\gamma$  the Poisson distribution on $\N_0$ with parameter $\gamma$, and for a measure $\mu\in\Mcal(\Scal)$ we denote with  $\Poi_\mu\in \mathcal M_1(\Mcal_{\N_0}(\Scal))$ the  distribution of a Poisson point process on $\Scal$ with intensity measure $\mu$. By this, we mean a finite random collection $\sum_i \delta_{r_i}$ of points $r_i\in\Scal$ such that $\#\{i\colon r_i\in  A\}$ has the distribution $\Poi_{\mu(A)}$ for any measurable $A\subset \Scal$ and is independent over mutually disjoint sets $A$. Similarly, for any measure $M$ on $\Gamma_T^{\ssup 1}$, by $\Poi_{M}$ we denote a finite collection $\sum_i \delta_{ \Xi_i}$ of point measures on non-empty point configurations $\Xi_i\in \Gamma_T^{\ssup 1}$ such that $\#\{i\colon \Xi_i\in B\}$ is Poisson-distributed with parameter $M(B)$ for any measurable $B\subset \Gamma_T^{\ssup 1}$.

\subsection{First main result: identification of the distribution}\label{sec-identdist}

In this section, we present our first main result: an identification of the distribution of the empirical measure $\Vcal_N^{\ssup T}$ of the coagulation trajectories of the particles present at time $T$. This is in terms of a Poisson point process (PPP) on $\Gamma_T^{\ssup 1}$ and an exponential pair-interaction term. We fix $\mu\in\Mcal_1(\Scal)$ and will be considering the Marcus--Lushnikov process under the poissonised initial condition $\P_{\Poi_{N\mu}}$. 

In the following the reference measure $M^{\ssup T}_{\mu,N}\in \Mcal(\Gamma_T^{\ssup 1})$  will play an important role, which we define as
\begin{equation}\label{Mreferencemeasure}
M^{\ssup T}_{\mu,N}(\d \xi)=N^{|k|-1}\e\,\P_{\Poi_\mu}(\Xi\in\d \xi),\qquad \xi\in \Gamma^\ssup{1}_T, k = \xi_0,
\end{equation}
that is, the restriction of $\P_{\Poi_\mu}$ to $\Gamma^\ssup{1}_T$ with some factors that turned out convenient.
Note that $M^{\ssup T}_{\mu,N}$ does not weight empty configurations: while $\Poi_\mu$ gives positive measure to $k=0$, the restriction of $\P_0$ to $\Gamma_T^{\ssup 1}$ is the zero measure, since $\Gamma_T^{\ssup 1}$ does not contain the zero point measure trajectory. Furthermore, the total mass of $M^{\ssup T}_{\mu,N}$ is finite, since $\P_k$ has mass $\leq 1$ on $\Gamma_T^{\ssup 1}$ and $\Poi_\mu$ has exponential moments of all orders.

We now introduce a reference PPP $Y_N=\sum_{i\in I}\delta_{\Xi_i}$ on $\Gamma_T^{\ssup 1}$ with intensity measure $N M^{\ssup T}_{\mu,N}$. We write  expectation with respect to $Y_N$ as ${\tt E}_{N M^{\ssup T}_{\mu,N}}$.  The PPP of the initial configurations, $Y_{N,0}=\sum_i \delta_{\Xi_{i,0}}=Y_N\circ \pi_0^{-1}$, is also a PPP with intensity measure $N M^{\ssup T}_{\mu,N}\circ \pi_0^{-1}$.

Let us express the probability of non-coagulation between any two trees for a fixed time interval $[0,T]$ for fixed $T>0$. To do this, we need the labelled coagulation process $Z$. Fix disjoint finite label sets $A,B$ and vectors $\mathbf{x} \in \Scal^A$, $\mathbf{y}\in \Scal^B$. Suppressing the dependence on $T$ in the notation we consider the process $Z=(Z_t)_{t\in [0,T]}$ under $\P_{(\mathbf{x},\mathbf{y})}$, where the partition process $(P_t)_{t\in [0,T]}$ takes values in $\Pcal(A\cup B)$ with initial atom configuration $(\mathbf{x},\mathbf{y})\in\Scal^{A\cup B}$. We denote by $\{A\nleftrightarrow B\}$ the event that (up to time $T$) there is no coagulation between $A$ and $B$, i.e., no coagulation between any pair of particles (i.e., sets) $C,D$ with $C\subset A$ and $D\subset B$.  Now, we define
\begin{equation}\label{DefR}
R^{\ssup {T}}(\xi,\xi')= -\log \frac{\d \P_{(\mathbf{x},\mathbf{y})}(A\nleftrightarrow B, (\Xi^\ssup{T,A},\Xi^\ssup{T,B})\in \cdot)}{\d \P_{\mathbf{x}}\circ (\Xi(Z))^{-1}\otimes \P_{\mathbf{y}}\circ (\Xi(Z))^{-1}}(\xi,\xi^\prime),\qquad \xi, \xi' \in \Gamma_T.
\end{equation}
In words, $\e^{-R^{\ssup {T}}(\xi,\xi')}$ is the probability that the subprocesses $\xi=\Xi^\ssup{T,A}$ and $\xi'=\Xi^\ssup{T,B}$ do not coagulate with each other by time $T$. This probability is the probability that all the exponential holding times of pairs of particles between $\xi$ and $\xi'$ never elapse during $[0,T]$.  In Lemma \ref{lem-Rident} we prove the explicit formula
\begin{equation}\label{Rdef}
R^{\ssup{T}}(\xi,\xi')=\int_0^T \d t\,\langle \xi_t, K \xi'_t\rangle, \qquad \xi, \xi'\in \Gamma_T.
\end{equation} 
where we write $K\phi(x,m)=\int_\Scal\sum_{m\in\N} \phi(\d x',m')\, K((x,m),(x',m'))$ for any $(x,m) \in \Scal\times \N$ and any $\phi\in \Mcal_{\N_0}(\Scal\times \N)$.

In our first main result,  we assume nothing else than what we stated around \eqref{kernels}. We slightly generalise the definition \eqref{Mreferencemeasure} by defining for any $b\in (0,\infty)$ the measure
 \begin{equation}\label{Mbreferencemeasure}
M^{\ssup T}_{b\mu,N}(\d \xi)=N^{|k|-1}\e \,\P_{\Poi_{b\mu}}(\Xi\in\d \xi),\qquad \xi\in \Gamma^\ssup{1}_T, k = \xi_0,
\end{equation}
where the additional parameter $b$ will be used later to tune the total mass $|M^{\ssup T}_{b\mu,N}|$.

\begin{theorem}[Poissonian description of the empirical measure]\label{thm-distMi}
Fix $\mu\in\Mcal_1(\Scal)$ and $T>0$ and $N\in\N$ and  a measurable bounded function $f\colon\Mcal(\Gamma^\ssup{1}_T)\to[0,\infty)$.
Then, for any $b\in(0,\infty)$,
\begin{equation}\label{DistVcalPoiThm}
\begin{aligned}
\E_{\Poi_{N\mu}}\big(f(\Vcal_N^{\ssup T})\big)
&=
 {\tt E}_{NM^{\ssup T}_{b\mu,N}} \Big[\e^{-\frac 12\sum_{i,j\colon i\not=j}R^{\ssup T}(\Xi_i,\Xi_j)} \e^{(b-1)|Y_{N,0}|} b^{-\int Y_{N,0}(\d k) |k|} \, f(\smfrac {1}{N}Y_N)\Big]\\
 &\qquad \times \e^{N(|M^{\ssup T}_{b\mu,N} |-1)},
 \end{aligned}
 \end{equation}
 where $Y_N=\sum_i\delta_{\Xi_i}\sim\Poi_{N M^{\ssup T}_{b\mu,N}}$ is a PPP with intensity measure $N M^{\ssup T}_{b\mu,N}$.
\end{theorem}

The proof of Theorem \ref{thm-distMi} is in Section \ref{MiMaintroduction}. An easy interpretation of the formula can be given for $b=1$. In that case the distribution of $\Vcal_N^\ssup{T}$ is given via the distribution of $\frac{1}{N}Y_N$ with the additional density term $\exp\{-\frac{1}{2}\sum_{i\neq j} R^\ssup{T}(\Xi_i,\Xi_j)\}$ that takes into account the non-coagulation between any pair of coagulation trajectories.  
We derive the formula for an arbitrary $b\in(0,\infty)$ since we need later to ensure the finiteness of the limit of $|M^{\ssup T}_{b\mu,N} |$, when $K$ is replaced by $\frac 1N K$ and the limit $N\to\infty$ is taken. The correctional terms on the right-hands side of the formula involve the number of coagulation trajectories $|Y_{N,0}|$ and the sum of the number of atoms in each trajectory $\int Y_{N,0}(\d k) |k| = \sum_i |\Xi_{i,0}|$.

Actually, Theorem~\ref{thm-distMi} gives the coagulation process the meaning of a spatial {\it interacting many-body system} with Gibbsian interaction, see Remark~\ref{rem-manybody}. We underline that this is not a prerogative of the system when starting from a Poisson initial condition, as an analogue representation can be obtained for any initial condition. However, since in the following we focus only on Poisson initial condition, we omit it here and defer the treatment of such setting to future work.

\subsection{Second main result: a large-deviations principle}\label{sec-LDP}

\noindent Here is our second main result, which is about asymptotics as $N\to\infty$. We are mainly interested in laws of large numbers and gelation criteria, but in view of our Gibbsian representation in Theorem~\ref{thm-distMi}, this is most naturally approached in terms of a large-deviation principle for the empirical coagulation trajectory measure, as we will formulate in Theorem~\ref{thm-LDP}. We consider initial configurations distributed as $\Poi_{N\mu}$ for some $\mu\in\Mcal_1(\Scal)$, such that $N$ is the order of the number of particles at time zero. In order to see interesting phenomena, we replace now the kernel $K$ by $\frac 1N K$. We indicate this by adding an additional superscript $^{\ssup N}$. It is clear that the transition from $K$ to $\frac 1N K$ is mathematically equivalent to the transition from the time interval $[0,T]$ to $[0,TN]$. We choose to keep the time horizon as $[0,T]$ and  to rescale the rates of coagulation by $1/N$, since topologies on processes with fixed time horizons are handled in a more standard way than with diverging horizons. Having done this change, the formula for the particle distribution from Theorem~\ref{thm-distMi} receives a structure that is exponential in $N$, hence any kind of limiting assertions as $N\to\infty$ based on that formula will naturally involve a large-deviation principle.

For completeness, let us spell out the definition of a large-deviation principle (LDP): Let $\Ycal$ be a topological space that is equipped with the Borel-$\sigma$-algebra. One says that a sequence $(\mu_N)_{N\in\N}$ of probability measures on $\Ycal$ satisfies an LDP with speed $N$ and rate function $I\colon\Ycal\to[0,\infty]$ if, for any closed, respectively open, set $F, G\subset \Ycal$,
$$
\limsup_{N\to\infty}\smfrac 1N\log\mu_N(F)\leq -\inf_F I\qquad\mbox{and}\qquad \liminf_{N\to\infty}\smfrac 1N\log\mu_N(G)\geq -\inf_G I.
$$
This concept depends strongly on the topologies used, and we will be working on sets $\Ycal=\Mcal(\Xcal)$ of bounded Borel measures on a topological space $\Xcal$, so we emphasise that we will always equip them with the weak topology, which is induced by test integrals against bounded and continuous functions $\Xcal\to\R$. See Section~\ref{sec-topology} for specific remarks for the cases where $\Xcal$ is itself a space of measures.

For all limiting assertions in this paper, we will be under the following assumption:
\begin{equation}\label{AssK1}
H:=\sup_{v,w\in\Mcal(\Scal\times \N)\colon \|v\|_1,\|w\|_1\leq 1}\langle v, Kw\rangle <\infty,
\end{equation}
where we write $Kw(x,m)=\int_\Scal\sum_{m\in\N} w(\d x',m')\, K((x,m),(x',m'))$ for $(x,m)\in \Scal\times \N$ and $w\in \Mcal(\Scal\times \N)$.
Note that product kernels of the form $K((x,m),(x',m'))=\varphi(x,x') mm'$ satisfy \eqref{AssK1} if $\varphi$ is bounded. The assumption in \eqref{AssK1}, together with the conditioning on $\Acal_{f,\Aold}$ below, will imply several decisive boundedness and compactness properties in the sequel; see Remark~\ref{rem-UppboundExpla} below.

In preparation for formulating our main result, we formulate the convergence of the reference measure $M_{\mu,N}^{\ssup T}$ defined in \eqref{Mreferencemeasure} with $K$ replaced by $\frac 1N K$. In a slight abuse, we use the notation $\Q|_A$ for the restriction of a measure $\Q$ to a measurable set $A$, i.e., $\Q|_A(B)=\Q(A\cap B)$. Both measures that we are introducing in the next lemma are crucial for our LDP. 

\begin{lemma}[Convergence of $N^{|k|-1}\P_k^{\ssup {N}}(\Xi\in\cdot)|_{\Gamma_T^{\ssup 1}}$]\label{lem-Qlimit}
Assume that the kernel $K$ satisfies \eqref{AssK1}. Fix $\mu\in\Mcal_1(\Scal)$ and $T\in(0,\infty)$ and $k\in\Mcal_{\N_0}(\Scal)$. Replace the coagulation kernel $K$ by $\frac 1N K$ and write $\P_k^{\ssup{N}}$ for the probability measure.
 Then the following limiting measure exists in the weak sense:
\begin{equation}\label{Qlimit}
\Q_k^{\ssup{T}}(\cdot) =\lim_{N\to\infty}N^{|k|-1}\P_k^{\ssup{N}}(\Xi\in\cdot)|_{\Gamma_T^{\ssup 1}}\in\Mcal(\Gamma_T^{\ssup 1}).
\end{equation}
In particular, the  measure $M_{\mu,N}^{\ssup {T,N}}$ defined in \eqref{Mreferencemeasure} w.r.t. $\P^\ssup{N}$ converges towards the measure
\begin{equation}\label{Mmudef}
M^{\ssup T}_\mu =
\e\,\Poi_\mu\otimes \Q^{\ssup T} \in \Mcal(\Gamma_T^{\ssup 1}).
\end{equation}
\end{lemma}

In \eqref{Mmudef} we used the well-known measure-theoretic notation
\begin{equation}\label{productmeasure}
\Poi_{\mu}\otimes \Q^{\ssup T}\big(\d (k,\xi)\big)=\Poi_{\mu}(\d k)\, \Q^{\ssup T}_k(\d \xi),\qquad k\in\Mcal_{\N_0}(\Scal),\xi\in\Gamma_T^{\ssup 1} \text{\,s.t.\,}\xi_0=k.
\end{equation}   

The proof of Lemma~\ref{lem-Qlimit} is in Section~\ref{sec-Qproof}.  An explicit formula for $\Q_k^{\ssup{T}}$ and more information appear in Lemma~\ref{lem-limitconnection}. Informally speaking, $\Q_k^{\ssup T}$ assigns to a trajectory $\xi$ the product of  the rates over all coagulation events, but drops all terms coming from the exponential densities of the coagulation times. It is an important reference measure for our further analysis.

As we have already said in Section~\ref{sec-Goals}, in the present paper we study only the microscopic part of the process, i.e., those particles that are of finite-order size in $N$. We leave the study of other particles, with size growing in $N$, to future work. Note that the space of coagulation trajectories $\Gamma^\ssup{1}_T$ does include trajectories with arbitrarily large sizes. For this reason, we encounter the problem of lack of compactness in $\Mcal(\Gamma^\ssup{1}_T)$, the state space of $\Vcal_N^\ssup{T}$, and we are forced to condition on a set that induces compactness.  Fix any  $\Aold>0$ and some function $f \colon \N\to [0,\infty)$ that grows at infinity faster than linear, i.e., $f(r)/r\to\infty$ as $r\to\infty$, and satisfes that $f(r)\geq r$ for all $r\in \N$. Define
\begin{equation}\label{Acalfdef}
\Acal_{f,\Aold}=\Big\{\nu\in\Mcal(\Gamma^\ssup{1}_T)\colon \int_{\Gamma^\ssup{1}_T} \nu(\d \xi) \,f(|\xi_0|)\leq \Aold\Big\},
\end{equation}
and note that for any $\xi\in \Gamma^\ssup{1}_T$ it holds that its collected mass $\|\xi_t\|_1\in \N$ is constant in $t\in [0,T]$ and equal to $|\xi_0|$, the number of initial atoms of $\xi$. Hence, the condition in $\Acal_{f,\Aold}$ is an higher-moment integrability condition on the sizes/masses of the particles at time $T$. On the event $\{\Vcal_{N}^{\ssup{T}}\in \Acal_{f,\Aold}\}$  we will be able to derive compactness/continuity of important objects that allow for a smooth proof of the LDP (see Section \ref{sec-topology}).
See Remark~\ref{rem-Conditioning} below for an explanation that, on this event, the configuration cannot develop a macroscopically large particle by time $T$ in the limit $N\to\infty$ and so the conditioning is more than a purely technical step. Since $f(r)\geq 1$ for all $r\in\N$, $|\nu|\leq \Aold$ for all $\nu\in \Acal_{f,\Aold}$.
We will later mainly use $f(r)=r^2$.

We are going to introduce the rate function of our LDP. For two finite measures $m,p$ on a Polish space $\Xcal$, we denote the relative entropy of $m$ with respect to $p$ by 
\begin{equation}\label{entropy}
H(m\mid p)=\int_\Xcal m(\d x)\log\frac{\d m }{\d p}(x)+p(\Xcal)-m(\Xcal),
\end{equation}
if $m\ll p$, and $H(m\mid p)=\infty$ otherwise. Note that $H(\cdot\mid p)$ is convex and nonnegative with the only zero $p$. Furthermore, all its sublevel sets $\{m\colon H(m\mid p)\leq C\}$ for $C\in\R$ are compact in the weak topology (a short proof of which is given in Section~\ref{sec-LDPproof}).

Recall from \eqref{DefR} and \eqref{Rdef} the non-coagulation functional $R^{\ssup T}$ and introduce the  operator with kernel $R^{\ssup{T}}$,
\begin{equation}\label{ROperator}
\mathfrak R^{\ssup {T}} (\nu)(\xi)=\int_{\Gamma_T^{\ssup 1}} R^{\ssup{T}}(\xi,\xi')\,\nu(\d \xi'),\qquad \xi\in \Gamma_T^{\ssup 1},\nu\in\Mcal(\Gamma_T^{\ssup 1}).
\end{equation}

We fix $\mu\in\Mcal_1(\Scal)$ and we define the function $I^\ssup{T}_\mu\colon \Mcal(\Gamma_T^{\ssup 1})\to[0,\infty]$ by
\begin{equation}\label{Idef}
\begin{aligned}
I^{\ssup T}_\mu(\nu)
&=
\begin{cases}
\Big\langle \nu,\log \frac{\d \nu}{\d M^{\ssup T}_\mu}\Big\rangle+\frac 12 \langle \nu,\mathfrak R^{\ssup {T}}(\nu)\rangle + 1 -|\nu|,&\text{if } \nu \ll M^{\ssup T}_\mu\\
\infty&\mbox{otherwise,}
\end{cases}
\end{aligned}
\end{equation}

Notice that $M^{\ssup T}_\mu$ might not have a finite total mass, but we will show in Lemma~\ref{lem-refmeasure-finite} that, under the assumption in \eqref{AssK1}, $|M^{\ssup T}_{b\mu}|<\infty$ for sufficiently small $b\in(0,\infty)$, since $M^{\ssup T}_{b\mu}(\d\xi)=M^{\ssup T}_\mu(\d\xi)\, \e^{1-b}b^{|\xi_0|}$. As a consequence, for any $b\in(0,\infty)$ such that $|M^{\ssup T}_{b\mu}|<\infty$, if the density  $\frac{\d \nu}{\d M^{\ssup T}_{\mu}}$ exists (which exists if and only $\frac{\d \nu}{\d M^{\ssup T}_{b\mu}}$ exists), we have the alternative representation
\begin{equation}\label{Imuchar}
I_\mu^{\ssup T}(\nu)=1-|M^{\ssup T}_{b\mu}| +H(\nu|M^{\ssup T}_{b\mu})+\frac 12 \langle \nu,\mathfrak R^{\ssup {T}}(\nu)\rangle +\int \nu_0(\d k)\, |k|\, \log b +(1-b)|\nu|.
\end{equation}
(Recall that we write $\nu_0=\nu\circ \pi_0^{-1}\in\Mcal(\Mcal_{\N_0}(\Scal))$ for the projection of a measure $\nu\in\Mcal(\Gamma^\ssup{1}_T)$.) 

Here is our main result: an LDP for the collection $\Vcal_{N}^{\ssup {T}}$ of all the components of the coagulation process with $N$ particles and kernel $\frac 1N K$, restricted to the event that no infinite component appears in the limit. We will also be working with the measure $M^{\ssup{ T, \leq L}}_{\mu}$ for some fixed $L\in \N$, the restriction of $M^{\ssup T}_{\mu}$ to the set
\begin{equation}\label{GammaLdef}
\Gamma_{T,\leq L}^{\ssup{1}}= \{\xi\in\Gamma_T^{\ssup 1}\colon \xi_T\mbox{ is concentrated on }\Scal\times [L]\},\qquad \mbox{where }[L]=\{1,\dots,L\},
\end{equation}
of particle trajectories with particles of sizes $\leq L$. It is an easy consequence from Lemma \ref{lem-taubound} that $|M^{\ssup {T, \leq L}}_{\mu}|<\infty$ under the assumption \eqref{AssK1}.
In the following we often identify measures $\nu$ on $\Gamma_T^{\ssup 1}$ that satisfy $\nu(\Gamma^\ssup{1}_T \setminus \Gamma^\ssup{1}_{T, \leq L}) = 0$ with measures on   $\Gamma_{T,\leq L}^{\ssup 1}$.

\begin{theorem}[LDP for $\Vcal_{N}^{\ssup {T}}$]\label{thm-LDP}
Assume that the kernel $K$ is continuous  and satisfies the upper bound in \eqref{AssK1}. Replace the kernel $K$ by $\frac 1N K$. 
Pick $T\in(0,\infty)$ and $\mu\in\Mcal_1(\Scal)$.
\begin{enumerate}
\item[(1)] Pick $f\colon \N\to[0,\infty)$ satisfying $\lim_{r\to\infty}f(r)/r=\infty$ and $f(r)\geq r$ for any $r$. Then, for any $\Aold>0$, the distribution of $\Vcal_{N}^{\ssup {T}}$ under $\P^{\ssup{N}}_{\Poi_{N\mu}}(\,\cdot\,|\Vcal_{N}^{\ssup {T}}\in\Acal_{f,\Aold})$ satisfies an LDP on $\Mcal(\Gamma^\ssup{1}_T)$
 with speed $N$ and rate function 
\begin{equation}\label{ratefunctionA}
\Mcal(\Gamma^\ssup{1}_T)\to [0,\infty],\qquad \nu\mapsto \begin{cases} I^{\ssup T}_\mu(\nu)-\chi_\Aold & \text{ if } \nu \in \Acal_{f,\Aold}\\
\infty &\text{ otherwise,}
\end{cases}
\end{equation}
where  $\chi_\Aold=\inf_{\nu\in\Acal_{f,\Aold}}I^{\ssup T}_\mu(\nu)$.
The sublevel sets of this rate function are compact.

\item[(2)] For any $L\in\N$, the distribution of $\Vcal_N^{\ssup{T}}$ under $\P_{\Poi_\mu}(\,\cdot\,|\Vcal_N^{\ssup{T}}(\Gamma_{T}^{\ssup 1}\setminus \Gamma_{T,\leq L} ^{\ssup 1})=0)$ satisfies an LDP on $\Mcal(\Gamma^\ssup{1}_T)$ with speed $N$ and rate function 
\begin{equation}
\Mcal(\Gamma^\ssup{1}_T)\to [0,\infty],\qquad \nu\mapsto \begin{cases} I_\mu^{\ssup{T,\leq L}}-\inf_{\Mcal(\Gamma_{T,\leq L}^{\ssup{1}})}I_\mu^{\ssup{T,\leq L}} & \text{ if } \nu \in \Mcal(\Gamma^\ssup{1}_{T,\leq L})\\
\infty &\text{ otherwise,}
\end{cases}
\end{equation}
where
\begin{equation}\label{ILdef}
I_\mu^{\ssup{T,\leq L}}(\nu)
=H(\nu|M^{\ssup {T,\leq L}}_{\mu}) +1- |M_{\mu}^{\ssup {T,\leq L }}|+ \frac 12 \langle \nu,\mathfrak R^{\ssup T}(\nu)\rangle,\qquad \nu\in \Mcal(\Gamma_{T,\leq L}^{\ssup{1}}).
\end{equation}
The sublevel sets of this rate function are compact.
\end{enumerate}
\end{theorem}

The proof of Theorem~\ref{thm-LDP} is in Section~\ref{sec-LDPproof}. 

A standard conclusion from Theorem~\ref{thm-LDP} is the following.

\begin{cor}[Accumulation points]\label{cor-accu}
In the situation of Theorem~\ref{thm-LDP}, the distributions of $\Vcal_{N}^{\ssup {T}}$ both under $\P^{\ssup{N}}_{\Poi_{N\mu}}(\,\cdot\,|\Vcal_{N}^{\ssup {T}}\in\Acal_{f,\Aold})$ and under $\P^{\ssup{N}}_{\Poi_\mu}(\,\cdot\,|\Vcal_N^{\ssup{T}}(\Gamma_{T}^{\ssup 1}\setminus \Gamma_{T,\leq L} ^{\ssup 1})=0)$ are tight in $N$, and each limit point of this distribution along any subsequence is concentrated on the set of minimisers of $I_\mu^{\ssup T}|_{\Acal_{f,\Aold}}$, respectively of $I_\mu^{\ssup{T,\leq L}}$.
\end{cor}

We are not only interested in the empirical process $\Vcal_N^{\ssup T}$, but also in the Marcus--Lushnikov process $(\Xi_t)_{t\in[0,T]}$ itself. Because of \eqref{Xi_as_function_of_nu}, this is a function of $\Vcal_N^{\ssup T}$, more precisely $(\frac 1N\Xi_t)_{t\in[0,T]}=\rho(\Vcal_N^{\ssup T})$, where we define 
\begin{equation}\label{rhodef}
\rho\colon \Mcal(\Gamma_T^\ssup{1})\to\D_T(\Mcal(\Scal\times \N)),\qquad \rho(\nu)=(\rho_t(\nu))_{t\in[0,T]}=\Big(\int_{\Mcal(\Gamma^\ssup{1}_T)}\nu(\d \xi)\,\xi_t\Big)_{t\in[0,T]},
\end{equation}
where we write $\D_T$ for the set of  càdlàg functions on $[0,T]$.

The continuity of the map $\rho$ is handled in the following lemma. Details about the topologies are given in Section \ref{sec-topology}. 

\begin{lemma}[Continuity of $\nu\mapsto \rho(\nu)$]\label{lem-Contrho}
Fix any  $\Aold>0$ and some function $f \colon \N\to [0,\infty)$ that grows at infinity faster than linear, i.e., $f(r)/r\to\infty$ as $r\to\infty$. Let $(\nu_n)_{n\in\N}$ be a sequence in $\Acal_{f,\Aold}$ that converges towards some $\nu$ that has a finite entropy $H(\nu|M_{b\mu}^{\ssup T})$ with respect to $M_{b\mu}^{\ssup T}$ for some $b>0$. Then $\rho(\nu_n)\to \rho( \nu)$ as $n\to\infty$.
\end{lemma}

The proof of Lemma~\ref{lem-Contrho} is in Section~\ref{sec-proofrhocont}. As a consequence of this and the LDP of Theorem~\ref{thm-LDP}, we obtain also an LDP for the Marcus--Lushnikov process:

\begin{cor}[LDP for $(\frac 1N\Xi_t)_{t\in[0,T]}$]\label{cor-LDPXi}
In the situation of Theorem~\ref{thm-LDP}(1), the distribution of $(\frac 1N\Xi_t)_{t\in[0,T]}$ satisfies an LDP on $\D_T(\Mcal(\Scal\times \N))$ with rate function
$$
\rho\mapsto \inf\big\{I_\mu^{\ssup T}(\nu)-\chi_\Aold\colon \nu\in\Acal_{f,\Aold},\rho(\nu)=\rho\big\}.
$$
\end{cor}

This immediately follows from  the contraction principle, more precisely from Remark (c) on Theorem 4.2.1 in  \cite{DZ10}. 

Let us close this section with a remark on a handy criterion for uniqueness of minimisers for the rate function $I_\mu^{\ssup T}$.

\begin{rem}[Convexity of $I_\mu^{\ssup T}$ by nonnegative definiteness of $K$]\label{Rem-convexity}
{The map $\nu\mapsto \langle \nu,\mathfrak R^{\ssup {T}}(\nu)\rangle$ is {\it a priori} not convex. But under the additional assumption that $K$ be nonnegative definite, it is. Then  $I_\mu^{\ssup T}$ is strictly convex. Recall that 
\begin{equation}
\label{Knonnegdef}
K\mbox{ is nonnegative definite }\Longleftrightarrow\qquad\langle v, Kv\rangle\geq 0,\qquad v\in\Mcal_\R(\Scal\times \N),
\end{equation}
where we denote by $\Mcal_\R(\Scal\times \N)$ the set of signed finite measures on $\Scal\times\N$. The convexity of the map under the assumption of nonnegative definiteness is clear from the fact that, for positive measures $\nu_1,\nu_2$ on $\Gamma_T^{\ssup 1}$ and $\alpha\in(0,1)$, 
\begin{multline}
\big\langle \alpha\nu_1+(1-\alpha) \nu_2,\mathfrak R^{\ssup{T}}(\alpha\nu_1+(1-\alpha) \nu_2)\big\rangle\\
=\alpha\langle \nu_1,\mathfrak R^{\ssup{T}}(\nu_1)\rangle+ (1-\alpha)\langle \nu_2,\mathfrak R^{\ssup{T}}(\nu_2)\rangle
-\alpha(1-\alpha)\big\langle \nu_1-\nu_2,\mathfrak R^{\ssup{T}}(\nu_1-\nu_2)\big\rangle.
\end{multline}
Nonnegative definiteness yields a handy criterion of uniqueness of minimisers of $I_\mu^{\ssup T}$, but since it is difficult to check in practical examples, we will not rely on it.
\hfill$\Diamond$}
\end{rem}

\subsection{Third main result: criteria for convergence and gelation}\label{sec-minimizer}

One of the main questions in the Marcus--Lushnikov model is about the  existence or non-existence of a {\em gelation phase transition}. That is, the question about the existence of a deterministic critical time threshold $T_{\rm gel}\in(0,\infty)$ such that there are only microscopically sized particles (i.e., particles of size of finite order, not depending on $N$) before time $T_{\rm gel}$, and after this time, a positive fraction of the total mass (i.e.\ $\asymp N$) lies in large  particles (i.e., particles of $N$-depending diverging size). If this phenomenon occurs, then we call $T_{\rm gel}$ the {\em gelation time}, the group of all the macroscopic particles the {\em gel} and the coagulation kernel $K$ a {\em gelling kernel}. We stick to the convention that we use $\frac 1N K$ instead of $K$ (indicated by an additional superscript $^{\ssup N}$).

Let us  coin a rigorous definition of the occurrence of gelation.
We introduce the notation $\|v\|_{1,\leq L}=\int_\Scal\sum_{m=1}^L m\, v(\d x,m)$ for $v\in\Mcal(\Scal\times \N)$, then $\|\Xi_T\|_{1,\leq L}$ is the total  amount of mass in particles with size $\leq L$ at time $T$. Recall that $t\mapsto \|\Xi_t\|_1$ is constant under $\P^{\ssup N}_{\Poi_{N\mu}}$, and $ \|\frac 1N\Xi_t\|_1$ is equal to $\frac 1N $ times a $\Poi_N$-distributed random variable, i.e., it converges to one as $N\to\infty$ almost surely and in $L^1$-sense. One can call the difference $\|\frac 1N\Xi_T\|_1-\|\frac 1N\Xi_T\|_{1,\leq L}$ the $L$-{\it gel} of the process at time $T$. Then
\begin{equation}\label{nongel}
\begin{aligned}
\G^{\ssup \mu}_T&=\lim_{L\to\infty}\limsup_{N\to\infty}\E^{\ssup N}_{\Poi_{N\mu}}\Big[\|\smfrac 1N\Xi_T\|_{1,\leq L}\Big]\\
&=1- \lim_{L\to\infty}\limsup_{N\to\infty}\E^{\ssup N}_{\Poi_{N\mu}}\Big[\|\smfrac 1N \Xi_T\|_1-\|\smfrac 1N\Xi_T\|_{1,\leq L}\Big]
\end{aligned}
\end{equation}
is the limiting expected non-gel mass, i.e., the mass outside the gel. The map $T\mapsto \G^{\ssup \mu}_T$ is non-increasing with initial  value $\G^{\ssup \mu}_0=1$.  If $\G^{\ssup \mu}_T<1$, then we say that there is a gel at time $T$, and we define the {\em gelation time} by
\begin{equation}\label{gelationtime}
T_{\rm gel}^{\ssup \mu}=\inf\Big\{T\in(0,\infty)\colon \G^{\ssup \mu}_T<1\Big\}\in[0,\infty]
\end{equation}
This is the time at which the gelation phase transition occurs, if it is finite. If $T_{\rm gel}^{\ssup \mu}<\infty$, we   also speak of the phenomenon of {\em loss of mass} and say that {\em gelation occurs}. The interpretation is that some positive fraction of all the atoms sits in particles of sizes that depend on $N$ and diverge as $N\to\infty$, such that their total mass goes lost when looking only at particles of finite size, regardless how large this finite-size window is. We think this notion is (one of) the most natural notions of gelation and gelation times; see Section~\ref{sec-literature} for other notions of gelation  used in the literature on coagulation processes.

Note that the total mass of the $L$-gel can also be expressed in terms of our process $\Vcal_N^{\ssup{T}}$. Indeed, introduce the measure 
\begin{equation}\label{clambdadef}
c_\lambda(\cdot)=\int_{\Mcal_{\N_0}(\Scal)} \lambda(\d k)\, k(\cdot)\in\Mcal(\Scal),\qquad\lambda\in\Mcal(\Mcal_{\N_0}(\Scal)),
\end{equation}
and its $L$-restriction $c_\lambda^{\ssup{\leq L}}(\cdot)=\int \lambda(\d k)\,k(\cdot)\1\{|k|\leq L\}$. Recall that $(\Vcal_N^\ssup{T})_0$ is a measure concentrated on $\Mcal_{\N_0}(\Scal\times\{1\})$ and can hence be identified with a measure on $\Mcal_{\N_0}(\Scal)$ which allows us to write $c_{(\Vcal_N^\ssup{T})_0}$. Recalling the connection between $\Vcal_N^\ssup{T}$ and $\frac 1N \Xi_t$ from \eqref{Xi_as_function_of_nu} we can write the non-$L$-gel of the process at time $T$ as
\begin{equation}\label{totalmass-identity}
\Big|c^{\ssup{\leq L}}_{(\Vcal_N^{\ssup {T}})_0}\Big| = \int_{\Gamma^\ssup{1}_T} \Vcal_N^\ssup{T}(\d \xi)|\xi_0|\1\{|\xi_0|\leq L\} = \int_{\Gamma^\ssup{1}_T} \Vcal_N^\ssup{T}(\d \xi)\|\xi_T\|_{1,\leq L}= \|\smfrac 1N\Xi_T\|_{1,\leq L}.
\end{equation}
Especially from the second expression it is clearly seen that this quantity is a continuous functional of $\Vcal_N^{\ssup {T}}$ in the weak topology. This will be important in our proofs, since they deal with $\Vcal_N^{\ssup {T}}$.

Introduce
\begin{equation}\label{qTdef}
\begin{aligned}
q_\mu^{\ssup T}&=\limsup_{n\to\infty}\Big(M_\mu^{\ssup T}\big(\{\xi\in\Gamma_T^{\ssup 1}\colon |\xi_0|=n\}\big)\Big)^{1/n}\in(0,\infty).
\end{aligned}
\end{equation}
This quantity controls the finiteness of all the moments of $ |\xi_0|$ under the reference measure $M_\mu^{\ssup T}$; the threshold for all the finiteness of all moments is $q_\mu^{\ssup T}=1$.

Here is our main result on the existence and non-existence of gelation. For all our sufficient criteria for gelation we will  assume that the following lower bound for the kernel $K$ holds:
\begin{equation}\label{AssK2}
h:= \inf_{v,w \in \Mcal(\Scal\times \N)\colon \|v\|_1, \|w\|_1 = 1}\langle v, Kw \rangle >0.
\end{equation} 
Like our main upper bound in \eqref{AssK1}, one main example for a kernel that satisfies  \eqref{AssK2} is $K((x,m),(x', m'))=\varphi(x,x')mm'$ with $\varphi$ bounded away from zero, but \eqref{AssK2} is much more general and admits many kinds of spatial dependencies in $K$.

\begin{theorem}[Criteria for non-gelation and for gelation]\label{thm-Lossofmass} Fix $\mu\in\Mcal_1(\Scal)$ and $T\in(0,\infty)$ and assume that \eqref{AssK1} holds.
\begin{enumerate}
\item[(1)] {\em Criterion for non-gelation:} Assume that $q_\mu^{\ssup T}<1$. Then the following hold.
\begin{enumerate}
\item $I^{\ssup T}_\mu$ has compact sublevel sets and hence possesses minimisers.
\item ${\rm NG}^{\ssup\mu}_T=1$, i.e., there is no gelation at time $T$. 
\item Any minimising $\nu^{\ssup T}$ satisfies the Euler--Lagrange equation
\begin{equation}\label{ELpoissonized}
\nu(\d \xi)=M_\mu^{\ssup T}(\d \xi)\,\e^{-\mathfrak R(\nu)(\xi)},\qquad\xi\in\Gamma_T^{\ssup 1}.
\end{equation}
\item The distributions of $\Vcal_N^{\ssup {T}}$ and $c_{(\Vcal_N^{\ssup {T}})_0}$ under $\P^{\ssup N}_ {\Poi_{N\mu}}$ are tight in $N$.
\item Let $\P$ be a limit point of $\P^{\ssup N}_ {\Poi_{N\mu}}(\Vcal_N^\ssup{T}\in \,\cdot\,)$, $N\in \N$. Denote by $\Vcal$ a random variable with distribution $\P$. Then $|c_{\Vcal_0}|=1$ almost surely with respect to $\P$. 
\end{enumerate}

\item[(2)] {\em Criterion for gelation}: Additionally to  \eqref{AssK1} assume that \eqref{AssK2} holds and that ${\inf I_\mu^{\ssup T}>0}$.  
Then the following hold.
\begin{enumerate}
\item ${\rm NG}^\ssup{\mu}_T<1$, that is, gelation occurs.

\item For any $L\in\N$,  every minimiser $\nu^{\ssup{T,\leq L }}$ of $I_\mu^{\ssup{T,\leq L}}$ defined in \eqref{ILdef} satisfies  the Euler--Lagrange equation on  $\Mcal(\Gamma_{T,\leq L}^{ \ssup 1})$:
\begin{equation}\label{ELmitL}
\nu(\d \xi)=M_\mu^{\ssup T}(\d \xi)\,\e^{-\mathfrak R(\nu)(\xi)},\qquad\xi\in\Gamma_{T,\leq L}^{\ssup 1}.
\end{equation}
\end{enumerate}
\end{enumerate}
\end{theorem}

The proof of Theorem~\ref{thm-Lossofmass} consists of several steps and is spread over the entire Section~\ref{sec-minimizerproof}. The main argument for (1)(b) is in Section~\ref{sec-subcrphase} and for (2)(a) in Section~\ref{sec-PhaseTransition}.

\begin{rem}[Interpretation of the EL-equation]
{The EL-equations in \eqref{ELpoissonized} and  \eqref{ELmitL}, respectively, characterise minimisers $\nu$ of the rate function  in terms of a self-referencing equation for the (non-normalised) distribution $\nu$ on one-particle trajectories: they are equal to a characteristic reference distribution $\e\Poi_\mu\otimes \Q^{\ssup T}$, with a term that weights the sampled trajectory with the probability that it does not coagulate with another independent sample under the same distribution $\nu$. Via \eqref{Xi_as_function_of_nu}, this can also be turned into a characteristic equation for all accumulation points of $(\frac 1N\Xi_t)_{t\in [0,T]}$, see Section~\ref{sec-Smol}.
\hfill$\Diamond$
}
\end{rem}

Now we state simple estimates on the kernel $K$ that imply that the criteria from Theorem~\ref{thm-Lossofmass} are satisfied.

\begin{prop}[Bounds on $T$ that imply (non-)gelation or convergence]\label{prop-criteria} 
 Fix $\mu\in\Mcal_1(\Scal)$ and $T\in(0,\infty)$ and assume that \eqref{AssK1} holds.
\begin{enumerate}
\item[(1)] {\em Criterion for non-gelation:}
\begin{itemize}
\item[(a)] if $TH< \frac 1{\e^2}$, then $q_\mu^{\ssup T}<1$ and thus the statements from Theorem~\ref{thm-Lossofmass}(1) apply.

\item[(b)] If $TH<\frac 1{\e^2\,}\frac\pi{1+\pi}$, then \eqref{ELpoissonized} has at most one solution, and  the distribution of $(\frac 1N \Xi_t)_{t\in [0,T]}$ under $\P^{\ssup N}_ {\Poi_{N\mu}}$ converges to $(\int \nu^{\ssup T}(\d \xi)\,\xi_t)_{t\in[0,T]}$ with $\nu^{\ssup T}$ the unique solution to \eqref{ELpoissonized}.
\end{itemize}
\item[(2)] {\em Criterion for gelation:} Additionally to  \eqref{AssK1} assume that \eqref{AssK2} holds. Then
\begin{itemize}
\item[(a)] it holds that \begin{equation}\label{Ilowboundthm}
\inf_{\nu\in\Mcal(\Gamma_T^{ \ssup 1})} I^{\ssup T}_\mu(\nu)\geq 1- \frac 1{2T}\Big( \frac{\e}{\pi H} + \frac{(\log(2TH\e^2))^2}{h} \Big),
\end{equation}
that is, the criterion of Theorem~\ref{thm-Lossofmass}(2) applies for all $T$ such that the right-hand side is strictly positive.

\item[(b)] It holds that $|c ^{\ssup{\leq L}}_{\nu^{\ssup {T}}_0}| \leq \frac{2 \log T}{hT}$ for any $L\in\N$ and every minimiser $\nu^{\ssup{T,\leq L }}$ of $I_\mu^{\ssup{T,\leq L}}$ defined in \eqref{ILdef}.
\end{itemize}
\end{enumerate}
\end{prop}

The proof of assertion (1) of Proposition~\ref{prop-criteria} is at the end of Section~\ref{sec-subcrphase}, and the proof for (2) is in Section~\ref{sec-ELeq}. 

\begin{cor}[Bounds on the gelation time]
Under the assumption in \eqref{AssK1}, $T_{\rm gel}^{\ssup\mu}\geq \frac 1{H\e^2}$. If additionally \eqref{AssK2} is assumed, then \[T_{\rm gel}^{\ssup{\mu}}\leq\inf\{T\colon \frac 1{2T}( \frac{\e}{\pi H} + \frac{(\log(2TH\e^2))^2}{h})<1\} <\infty.\]
\end{cor}

The correct interpretation of the loss of mass phenomenon in terms of a measure $\nu$ on $\Gamma^{\ssup 1}_T$ is that the total mass of all the coagulation trajectories decays in the limit $N\to\infty$ as a function of $T$. Indeed $\nu(\d \xi)$ registers only those coagulation trajectories that end up in finite-size particles at time $T$. All the other (larger) particles that are present at time $T$ do not appear in $ \nu$ anymore, together with their coagulation trajectories. The reason is that the state space $\Gamma_T^{\ssup 1}$ cannot accomodate non-microscopic trees. In future work, we plan to extend the description of the coagulation trajectories by some enlarged space that is able to describe also the larger particles and their coagulation trajectories. It seems as if it is not possible to formulate (not to mention, prove) an LDP or a law of large numbers without detailed knowledge about the large particles, since each of them has a non-trivial influence on the dynamics of the entire coagulation process after they appear.

\subsection{Fourth main result: the Smoluchowski equation}\label{sec-Smol}

Let us now come to the fourth main result of this paper. To introduce it, let us recall that the rigorous analysis of coagulation processes started in 1916 with a formulation and analysis of the -- by now famous -- Smoluchowski equation \cite{Smol16b}. 
This is a partial differential equation for the evolution of the concentration of particles sitting in any site. It consists of a positive term that describes the formation of new particles in a certain site via coagulation of a pair of smaller particles and a negative term describing the loss of particles in that site due to coagulation with any other particle. In our situation, where we include space, we can formulate the equation in the form
\begin{equation}\label{SmolEQ}
\begin{aligned}
\partial_t \rho_t(\d x^*,m^*)&=\sum_{m,m'\in\N\colon m+m'=m^*}\int_\Scal\int_{\Scal}\rho_t(\d x,m)\rho_t(\d x',m') {\bf K}\big((x,m),(x',m'),\d x^*\big) \\
&\qquad -
\rho_t(\d x^*,m^*) \,  K \rho_t(x^*,m^*),\qquad x^*\in \Scal, m^*\in\N,
\end {aligned}
\end{equation}

This equation and its variants play a fundamental role in the investigation of the limiting behaviour of the  Marcus--Lushnikov process, as well as being of interest in its own right as a deterministic coagulation model. Indeed, in previous probabilistic investigations, proofs of the convergence of the process $(\frac 1N \Xi_t)_{t\in[0,\infty)}$ with kernel $\frac 1N K$ often (if not always) follow the route that (1) tightness arguments are employed, (2) it is shown that every accumulation point satisfies the Smoluchowski equation, and (3) criteria for the uniqueness of the solution are given, see Section~\ref{sec-literature}. If the kernel is such that gelation occurs, then one can only expect this convergence before the gelation time, whereas the full process $(\frac 1N \Xi_t)_{t\in[0,\infty)}$ has a limit that solves an extension of the Smoluchowski equation, called the Flory equation, which captures the influence of the gel on microscopic particles after the gelation time.

The main results of this paper on convergence and  gelation so far have nothing to do with the Smoluchowski equation. Nevertheless, this equation is so important that we decided to show that it is satisfied by all limit points of the process $(\frac 1N \Xi_t)_{t\in[0,T]}$, if $T$ is small enough such that we have no gelation. The novelty here is that we derive it from our main characterization, the Euler--Lagrange equation in \eqref{ELpoissonized}.

\begin{lemma}[The ML process converges to a solution to the Smoluchowski equation]\label{lem-Smol}
Assume that $\mu\in\Mcal_1(\Scal)$ and that $K$ satisfies \eqref{AssK1} and that $TH<\frac{1}{\e^2}\frac\pi{1+\pi}$. 
Then, under $\P_{\Poi_{N\mu}}^{\ssup N}$, the process $(\frac 1N \Xi_t)_{t\in[0,T]}$ converges to a solution $\rho$ of the Smoluchowski equation  in \eqref{SmolEQ}.
\end{lemma}

The proof is in Section~\ref{sec-Smolproof}. We use Proposition~\ref{prop-criteria} to ensure that we are in a regime, where gelation does not occur and to get convergence of $\Vcal_N^{\ssup T}$ to the (unique) solution of the Euler--Lagrange equation \eqref{ELpoissonized}. 

\section{Background discussion}\label{sec-discussion}

\subsection{Comments on the main results}\label{sec-comments}

In this section we outline our key heuristics and explain the nature of some crucial difficulties in the proofs.
We also highlight the benefits and shortcomings of our approach.

\begin{rem}[The coagulation process as a Gibbsian  many-body system]\label{rem-manybody}
In our first main theorem (Theorem~\ref{thm-distMi}), the expectation on the right-hand side of \eqref{DistVcalPoiThm} is with respect to a Poissonian reference measure on point measures $Y_N$ of coagulation trajectories on $[0,T]$.  The time marginal at time zero is the spatial distribution of the initial particles (atoms) that coagulate during $[0,T]$ into one particle, and the time marginal at time $T$ is the particle distribution at the end. By the Poisson nature, these trajectories are {\em a priori} mutually independent, but  they  are under an exponential pair interaction term that expresses that they do not coagulate by time $T$. Note that this interaction is mutually repellent.  Hence, the right-hand side is a many-body system of  points in $\Scal$ with marks (the mark being their coagulation trajectory) with Gibbsian pair interaction with a PPP as the underlying reference measure. Putting $f=1$, we can identify the last exponential term ($\e^{N(|M^{\ssup T}_{\mu,N}|-1)}$) in terms of the first one and obtain the formula (for $b=1$)
\begin{equation}
\label{DistVcalPoiThmGibbs}
\E_{\Poi_{N\mu}}\big(f(\Vcal_N^{\ssup T})\big)=\frac{
 {\tt E}_{NM^{\ssup T}_{\mu,N}} \Big[\e^{-\frac 12\sum_{i,j \colon i\not=j}R^{\ssup T}(\Xi_i,\Xi_j)}  \, f\big(\smfrac{1}{N}Y_N\big)\Big]}{
 {\tt E}_{NM^{\ssup T}_{\mu,N}} \Big[\e^{-\frac 12\sum_{i,j \colon i\not=j}R^{\ssup T}(\Xi_i,\Xi_j)}\Big] }.
\end{equation}
This shows that $\Vcal_N^{\ssup T}$ has the distribution of $\frac 1N Y_N$ under exponential tranformation with the density $\e^{-\frac 12\sum_{i\not=j}R^{\ssup T}(\Xi_i,\Xi_j)}$, properly normalised.\hfill$\Diamond$
\end{rem}

\begin{rem}[Plausibility of the LDP]\label{rem-heuristic}
On the basis of the Poissonian description in Theorem~\ref{thm-distMi}, one can easily guess that there might be an LDP valid for $\Vcal_N^{\ssup{T}}$ with rate function as in \eqref{Imuchar}, after replacing $K$ by $\frac 1N K$ and hence $R^{\ssup T}$ by $\frac 1N R^{\ssup T}$ in \eqref{DistVcalPoiThm}. The main point is that $\frac 1N Y_N$ satisfies an LDP under ${\tt E}_{N M_{\mu,N}^{\ssup {T,N}}}$ with rate function $H(\cdot| M_\mu^{\ssup T})$, since $M_{\mu,N}^{\ssup {T,N}}$ converges weakly towards $M_\mu^{\ssup T}$. The interaction term in the exponent (the double-sum on $i\not=j$) is approximated by the sum on  {\em all} $i,j$ and directly leads to the double-integral of $\mathfrak R^{ \ssup T}$ with respect to $\nu\otimes \nu$. Now collect all the exponential terms on the right-hand side of \eqref{DistVcalPoiThm} to see that they lead directly to the formula for  $I_\mu^{\ssup T}$ in \eqref{Imuchar}.
\hfill$\Diamond$
\end{rem}

\begin{rem}[Difficulties]\label{rem-difficulties}
The heuristics of Remark~\ref{rem-heuristic} suggest that $I_\mu^{\ssup T}$ governs an LDP for $(\Vcal_N^{\ssup{T}})_{N\in\N}$ under $\P_{\Poi_{N\mu}}^{\ssup N }$ without any conditioning, but this is not true in general:
A couple of problems arise when attempting to prove an unconditional LDP and they turn out to be substantial issues, not merely technical difficulties. Indeed, these problems have a lot to do with the gelation phase transition that we discuss in Section~\ref{sec-minimizer}. 
One problem is that the approximation of the sum on $i\not= j$ in the interaction term by the sum on all $i,j$ (i.e., the addition of all the self-interactions) fails if the particles are too large, more precisely, if some of them are of a size proportional to $N$.
Another problem is that we do not see an argument for compact sublevels sets of $I_\mu^{\ssup T}$ (not even for lower semi-continuity) on the entire set $\Mcal(\Gamma_T^{\ssup 1})$ if $|M_\mu^{\ssup T}|=\infty$, hence existence of minimisers is not certain. In this case, one would pick $b\in(0,1)$ appropriately in order to make $|M_{b\mu}^{\ssup T}|$ finite, but then the term $\int \nu_0(\d k)\, |k|$, which is lower-semicontinuous in $\nu$, but in general not continuous) has a negative prefactor in $I_\mu^{\ssup T}(\nu)$. A third problem is that we see {\em a  priori} no argument for the fact that the infimum of $I_\mu^{\ssup T}$ over $\Mcal(\Gamma_T^{\ssup 1})$ should be equal to zero; in fact we disprove it in Section~\ref{sec-minimizer} under certain assumptions on $K$ and $T$. 

As mentioned, the above issues are not only technical. The reason why the unconditioned process $(\Vcal_N^{\ssup{T}})_{N\in\N}$  fails to satisfy an LDP  with rate function $I_\mu^{\ssup T}$ is the possible emergence of a macroscopic particle by time $T$. One can see that, under the assumptions in Theorem~\ref{thm-Lossofmass}(2), actually a different scenario arises, the emergence of a gel, and this makes an LDP with rate function $I_\mu^{\ssup T}$ impossible. Let us recall that, in the special case of an inhomogeneous Erd\H{o}s--R\'enyi graph (see Remark~\ref{Rem-ERgraph} and Section~\ref{sec-ERgraph}), it turned out that a non-trivial contribution to the true rate function without conditioning comes from  the macroscopic part of the configuration, which we neglect in Theorem~\ref{thm-LDP}. We plan to incorporate this part into our analysis in forthcoming work.
\hfill$\Diamond$
\end{rem}

\begin{rem}[The role of conditioning on $\Acal_{f,\Aold}$]
\label{rem-Conditioning} 
In Theorem~\ref{thm-LDP} we use two strategies to avoid the problems deriving an LDP for $(\Vcal_N^{\ssup{T}})_{N\in\N}$  described in Remark~\ref{rem-difficulties}.  In part (1) of the theorem we condition on $\Acal_{f,\Aold}$, while in part (2) we condition on only having particles with bounded size.
It is simple to see that the map $\nu\mapsto \int \nu(\d \xi)\, |\xi_0|$ is continuous and bounded on the set $ \Acal_{f,\Aold}$. Furthermore, this enables us to show that $I_\mu^{\ssup T}$ has compact sublevel sets on it, and the diagonal sum (i.e., the sum on $i=j$) in the exponent can be shown to be small (also using the assumption in \eqref{AssK1}); see Lemma~\ref{lem:diag-term}.
This makes the proof of the LDP of Theorem~\ref{thm-LDP} practical. 

The compactness of the level sets of the rate function indeed hinges on the restriction to the set $\Acal_{f,\Aold}$. In \eqref{Imuchar}, at least for large $T$, $b$ needs to be taken small, and then the one-but-last term is not lower-semicontinuous in $\nu$. However, on $\Acal_{f,\Aold}$, it is easily seen to be even continuous. Then the compactness of the sublevel sets of the rate function follows easily from the lower semicontinuity of $\nu\mapsto \langle \nu,\mathfrak R^{\ssup {T}}(\nu)\rangle$ (see Lemma~\ref{Lem-RfrakProp}) and the compactness of the level sets of the entropy.

However, on the event $\{\Vcal_{N}^{\ssup{T}}\in \Acal_{f,\Aold}\}$, there can be no loss of mass in the limit $N\to\infty$ for the random coagulation process $\Vcal_{N}^{\ssup{T}}$, since the total mass of atoms present at time $T$ in particles larger than $L$ is bounded  for any $L\in(0,\infty)$. Indeed, using the majorizing function $f(r)=r^2$ implies that
$$
\int (\Vcal_{N}^{\ssup{T}})( \d \xi)\, \|\xi_0\|_1\1\{\|\xi_0\|_1>L\}\leq  \frac \Aold L.
$$
This means that the conditioning rules out any occurrence of mesoscopic or macroscopic particles with total mass of order $\asymp N$. In this way, $f$ acts like a majorant that induces tightness.  
In our arguments in Section~\ref{sec-minimizerproof}, we will use only $f(r)=r^2$ and large $\Aold$, since this is easy to handle.

One might think that, instead of conditioning on $\{\Vcal_{N}^{\ssup{T}}\in \Acal_{f,\Aold}\}$, in the proof of the LDP, one can use a decomposition into this set and its complement and try to show that the probability of the complement is exponentially small in $N$ with a very large rate if $\Aold$ is large. This argument, if it could be carried through, would indeed lead to a proof of the LDP without conditioning. However, it is not successful, since the probability of such complementary event is exponentially small only when {gelation} is super exponentially unlikely to occur before time $T$ and therefore it requires an a priori knowledge about the occurrence of gelation.
\hfill$\Diamond$
\end{rem}

\begin{rem}[The role of the upper bound in \eqref{AssK1}]
\label{rem-UppboundExpla} 
From a technical perspective, the assumption in \eqref{AssK1}, in combination with the conditioning on $\Acal_{f,\Aold}$, implies boundedness, continuity and compactness for a number of crucial terms and sets that play decisive roles in the proof of the LDP of Theorem~\ref{thm-LDP}.

With regard to contents, any upper bound for $K$ obviously puts a lower bound for the inter-coagulation times, i.e., makes the early occurrence of large particles improbable. One might regard \eqref{AssK1} as a criterion under which there is a non-trivial non-gelation phase.  However, this connection is not easy to be understood; not even using comparison techniques like pathwise couplings to random-graph models like in \cite{AnIyMa23}. 

There is one special case where the value of the quantity in \eqref{AssK1} stands in a direct connection to the gelation time, namely for the well-known product kernel $K_H(m,m')=H m m'$ in the spaceless setting. For this process  a lot is known about the existence of the gelation phase transition and even the value of the gelation time, which is equal to $1/H$. See for example \cite{AKP21} and Section~\ref{sec-ERgraph} for the analysis of this  model from a large deviations point of view. \hfill$\Diamond$
\end{rem}

\begin{rem}[Difficulties in proving (non-)gelation]\label{Rem-DifficultGelation}
If $\Vcal_N^{\ssup{T}}$ were to satisfy an LDP under $\P_{\Poi_{N\mu}}^{\ssup N}$ with rate function $I_\mu^{\ssup T}$ as discussed in Remark~\ref{rem-heuristic}, then one would expect that it converges (at least along subsequences) to a minimiser $\nu^{\ssup T}$ of $I_\mu^{\ssup T}$, which then satisfies the EL-equations in \eqref{ELpoissonized}.
In this setting one might hope that the definition of gelation via \eqref{nongel} would be equivalent to mass loss in the minimiser, that is, $\int \nu^{\ssup T}(\d\xi)\,|\xi_0|<1$.

However, as we already pointed out in Remark~\ref{rem-difficulties}, the truth is much more delicate.
Nevertheless, we succeeded in proving that there is no mass loss in either sense under the assumption $q_\mu^{\ssup T}<1$,  implying no gelation.
However, the condition $q_\mu^{\ssup T}<1$ is only sufficient, but not necessary; we dropped the exponential term in \eqref{ELpoissonized} when estimating in Lemma~\ref{lem-estiM}. Actually, the criterion $q_\mu^{\ssup T}<1$ for non-gelation (and convergence and characterization of limits) is good only for $T$ small enough, depending on the upper bound on the kernel $K$ in \eqref{AssK1}.
Unfortunately, even for kernels that are known never to produce a gel, we are currently not able to derive $q_\mu^{\ssup T}<1$ from our criterion; for example, for the constant kernel $K\equiv H$ we can currently deduce only that $q_\mu^{\ssup T}\leq \frac 12 HT$. This shows that our non-gelation criterion is only sufficient and may be far from sharp.

On the other hand, under \eqref{AssK2}, we did prove that $\int \nu^{\ssup T}(\d\xi)\,|\xi_0|<1$ for $T$ sufficiently large, however, we were not able to use this to prove gelation. The reason is that we found no argument to show that $\Vcal_N^{\ssup{T}}$ converges to $\nu^{\ssup T}$ as $N\to\infty$ under any measure that satisfies an LDP with rate function $I_\mu^{\ssup T}$.  Instead, the main point in our proof of gelation for large $T$ in Section~\ref{sec-PhaseTransition} is the fact that $I_\mu^{\ssup T}$ is bounded away from zero for large $T$. 
\hfill$\Diamond$
\end{rem}

\subsection{Literature survey}\label{sec-literature}
In this section we wish to give an overview on related literature.
 
\subsubsection*{Coagulation models}
{\em Spaceless} models for coagulation have been studied for decades, and there are a number of works that derive criteria for the occurrence of gelation.  A review by Aldous~\cite{Ald99} gives a general overview, covering deterministic and  stochastic points of view.
He also suggests many open questions, several of which have since been resolved.

Mathematical modelling of coagulation began with Smoluchowski~\cite{Smol16b} in connection with his work on diffusion.
He wrote down a deterministic model in the form of a coupled set of ODEs, (together known as Smoluchowski equation), which he informally derived from an underlying stochastic model of Brownian particles.
The original Smoluchowski equation is the spaceless version of  \eqref{SmolEQ}.
A very natural stochastic (Markovian) model, which may be viewed as spatially homogenised limit of the stochastic model used by Smoluchowski, has been introduced by Marcus~\cite{Mar68}, and again by Gillespie~\cite{Gil72} and later studied by Lushnikov~\cite{Lus78b}. It is called the Marcus--Lushnikov model and we study a spatial version of it in this paper.

The stochastic setting is well-connected with the deterministic one. Indeed, replacing the coagulation kernel $K$ by $\frac 1N K$, several later authors prove, in the spaceless setting, the convergence of the normalised Marcus--Lushnikov process (written $\frac 1N \Xi$ in our setting) towards the solution of the Smoluchowski equation, under various assumptions on the kernel and on the initial condition. Some authors \cite{Nor99, FoGi04} prove convergence towards a more general version of the Smoluchowski version, the Flory equation~\cite{Flo41}, which characterises the evolution also in the presence of a gelation phase transition after the gelation time.

Much literature focuses on the deterministic setting, proving existence, uniqueness and other properties of solutions of the Smoluchowski equation and its variants, under more and more general assumptions~\cite{BaCa90, DuSt96, BoNiVe19, FeLuNoVe21}. In this setting, the phenomenon of gelation was initially interpreted as an explosion of moments of solutions at a finite time, later as the existence of a time at which the solution loses mass, i.e., the first moment strictly decreases.
The first rigorous treatment of gelation in the stochastic model comes with Jeon's work~\cite{Jeo98} (see also \cite{Re13} for extensions and generalisations), where several notions of gelation are considered. One  main notion is in terms of some kind of boundedness in $N$ (either in probability or in expectation) of the first time at which the total mass of the configuration that sits in particles larger than $\psi_N$ is larger than $\delta N$ for some $\delta >0$ and for some scale function $\psi_N$. The choice of $\psi_N\asymp N$ (i.e., the appearance of a macroscopic particle) leads to what is called a strong gelation. It is clear that all these notions {\it a priori} depend on the parameter $\delta$ and the scale function  $\psi_N$ and on the sense of boundedness that is required. The notion that we chose refers to $\psi_N\equiv L$ and an additional limit as $L\to\infty$, and we take expectation of the total mass, a notion that is not handled in \cite{Jeo98} nor in \cite{Re13}. Our gelation time defined in \eqref{gelationtime}
is -- so to speak -- the \lq earliest\rq\ of all possible gelation times.

It is intuitively clear that sufficient criteria for gelation to occur should consist of lower bounds for the coagulation kernel $K$. The most interesting yet simple example is the product kernel, $K(m,m')=C mm'$, which admits a one-to-one map of the coagulation process onto a natural growing sparse Erd\H{o}s--R\'enyi graph (see Section \ref{sec-ERgraph}), where gelation has a natural analogue with the famous emergence of a giant component. In this special case, one can identify the gelation time, and it turns out that at that time a macroscopic particle arises, hence it is a strong gelation time in the above notion.

In~\cite{Jeo98} Jeon proved that a sufficient condition for gelation is that $K(m,m')\geq \epsilon (m m')^q$ for all $m,m'\in\N$, with some positive $\epsilon$ and $q>\frac 12$.
Later Rezakhanlou~\cite{Re13} proved that kernels satisfying $K(m,n)\geq m^q +n^q$ with $q>1$ produce instantaneous gelation, i.e., the time at which a giant particle appears tends to zero. 
Let us mention that it is in general believed in the applied mathematics areas that all homogeneous kernels, i.e., kernels that satisfy $K(cm,cm')=c^{\gamma}K(m,m')$ for all $m,m'\in\N$, are gelling for $\gamma>1$, but as far as we know there is not a proof for this statement.

The first attempt to include space or particle features other than mass is made by Norris~\cite{Nor00}. He introduces what he calles a \emph{cluster coagulation model}, where each particle is called a cluster  and is attached to some point in some measurable space $E$. With some mass function $m\colon E\to (0,\infty)$, the coagulation of a pair of particles is a replacement of two with one single particle such that the sum of the masses before and after coagulation is preserved. Choosing $E=\Scal\times \N$ makes our model a special case of this. 
Norris focuses on the deterministic side, proving existence of solutions to the cluster version of the Smoluchowski equation, in our setting \eqref{SmolEQ}, and conservation of mass in some cases, i.e., when $K(x,y)\leq \phi(m(x))\phi(m(y))$ for all $x,y\in E$ and $\phi$ is at most linear. Moreover Norris calls the kernel $K$ \emph{approximately multiplicative} if
\[
\epsilon m(x)m(y)\leq K(x,y)\leq M(1+m(x))(1+m(y)),\qquad x,y\in E,\]
for some $\epsilon>0$ and $M<\infty$. In Theorem~2.2 he  proves that equation \eqref{SmolEQ} in this case has a unique mass-preserving solution up to a certain time $T$, which can be upper and lower bounded by functions of the initial condition. This time $T$ can be seen as another definition of gelation time. This gives upper and lower bounds on the gelation time under the above conditions, which can be seen as special cases of our assumptions~\eqref{AssK1} and~\eqref{AssK2}.  In~Section~4 Norris  studies convergence of  $\smfrac 1N\Xi$ (in our notation) to the solution of \eqref{SmolEQ} before gelation, under the assumption that the space $E$ is countable. In one special case, which he calls the (eventually) multiplicative kernel, he proves convergence to the solution of (the equivalent of) the Flory equation, even beyond gelation. Let us underline that Norris proves also that the convergence holds at least exponentially fast in $N$, which is a first indication that a large deviation principle might hold in this case.
Recently, in~\cite{AnIyMa23}, the authors extend Norris' convergence results  to convergence towards the solution of a generalised Flory equation, when $E$ is  a $\sigma$-compact metric space.
Like in all mentioned papers, the proof methods work on the generator of the coagulation process and employ martingale arguments.


The uniqueness of solutions to the Smoluchowski equation is quite a delicate question.
Under relatively strong assumptions, uniqueness has been established both before and after gelation \cite{NoZa11}.
However, see~\cite{Nor99} for a (quite involved) example of a spaceless kernel for which multiple solutions of the Smoluchowski equation are possible.
It is not clear whether the introduction of spatial positions for the particles influences uniqueness in a positive or in a negative way. Our work does not shed any light on this question, since our uniqueness assertions are only for the minimiser of the rate function.

Let us mention that there are few other models of coagulation, where the coagulating particles move in Euclidean space: \cite{Lang80,HaRe07} treat the case where the particles move as independent Brownian motions in $\R^d$ and \cite{Gui01} where the Brownian motions are replaced by random walks.
The works \cite{Bel03,Wel06} deal with diffusive particles with very general interactions of which coagulation is only a special case.
All these investigations proceed under assumptions that exclude gelation, whereas \cite{SiWa06, Wag06} restrict to particles moving as independent Markov chains on some finite state space in order to investigate gelation in detail.

\subsubsection*{Large deviations for jump Markov processes} We turn now to large deviation results. In particular, we would like to draw comparisons between our Corollary~\ref{cor-LDPXi} and existing results of  dynamical LDPs, i.e.,~LDPsx for the path of the empirical measure of weakly interacting jump Markov processes over a compact time interval $[0,T]$. 
Results of this type have been proved in the case that the empirical measure takes value in a finite dimensional space, starting from \cite{SW05} and then, under increasingly general assumptions, see for example~\cite{ADE18, AAPR22} and references therein. The classical approach follows  Freidlin-Wentzell theory; it consists in a tilting argument, for which a law of large numbers for the transformed dynamics is needed. Alternative approaches have appeared recently  (using fluxes \cite{PR19}, Hamilton-Jacobi equations \cite{Kra21}, weak convergence methods~\cite{DRW16}), however they all seem to be restricted to the finite dimensional case at the moment. The rate function of the dynamical LDP proved in these ways has a specific form: it is the time integral of an action functional, which involves a minimisation problem and the generator of the Markov process. It is worth mentioning that such a form of the rate function plays a key role in the understanding of the model in terms of gradient flows, see~\cite{MPR14}. 

The case of a LDP for infinite dimensional empirical measure has been under investigation mainly in kinetic theory, in relation to large deviations for Kac type of particle systems (stochastic microscopic models for Boltzmann equation). In this case particles are characterised by a spatial location and a velocity, they interact pairwise by changing their velocities and preserving kinetic energy. Here, the empirical measure is a measure on the space of locations and velocities, therefore it is a truly infinite dimensional object. For this type of models a rate function of the classical Freidlin-Wentzell type was suggested by L\'eonard~\cite{Leo95} via a large deviation upper bound. However matching  lower bounds were proved only when restricting to classes of sufficiently good paths~\cite{Hey23, BBBO21, Sun21}, the difficulty lies in the lack of a law of large numbers for the perturbed dynamics.
In particular, Heydecker in~\cite{Hey23} shows that this ansatz is not always successful, indeed even if the particle system almost surely preserves energy, the large deviation behaviour involves also paths that do not preserve energy and this happens with a rate that is not consistent with L\'eonard's rate function.
In~\cite{BBBC22, BBCB22} two models of Kac type are studied and a new rate function is suggested, which assigns a non-trivial rate to paths with increasing energy. The matching lower bound is then proved for a (larger) class of paths, but not yet for all paths.
Let us mention that upper and lower bounds do coincide in some infinite dimensional case, as they do in the finite dimensional one. This has been recently proved by Sun in~\cite{Sun21}, where L\'eonard's ansatz is extended to a more general setting. Sun focuses on a dynamical LDP for interacting particles with infinite dimensional empirical measure (including Kac type of model and coagulation fragmentation ones) on a restricted class of paths. Moreover matching upper and lower bounds for all paths are proved for coagulation-fragmentation models of Becker-D\"oring type under certain assumptions on the rates of interaction. 

 How is our Corollary~\ref{cor-LDPXi} linked to the above models? The role of conservation of energy for Kac type of models is played in our coagulation case by the conservation of the mass. Moreover, we are proving  a conditional LDP, i.e.,~we prove matching upper and lower bounds only for paths in a certain subset of the space (the image of the subset $\Acal_{f,\Aold}$ under $\rho$). Exactly as predicted in~\cite{Hey23, BBBC22} for Kac type of models, we expect the true rate function (without ruling out gelation) to be different but we defer this to future work.

\subsection{Comparison to inhomogeneous Erd\H{o}s--R\'enyi graph}\label{sec-ERgraph}

\noindent There is a special case in which our spatial coagulation process  can be realised as the connected component size process of an inhomogeneous Erd\H{o}s--R\'enyi graph, which was analysed in  \cite{AKLP23}.
The graph dynamics differs from the coagulation mechanism studied here only in that, for the graph, we add an edge between atoms (and leave their locations untouched) instead of replacing two particles by with a single particle at a new location. 
The connected components of the graph then play the role of the particles of the coagulation system.



Here is the special form of the coagulation and placement kernels that admit a representation as a graph process provided that $\Scal$ is a convex subset of a linear space.  Consider
\begin{equation}\label{specialkernel}
K\big((x,m),(x',m')\big)= \kappa(x,x')\, m m',\qquad x,x'\in\Scal, m,m'\in\N,
\end{equation}
with some symmetric and bilinear function $\kappa\colon \Scal\times\Scal\to[0,\infty)$, together with the deterministic placement kernel 
\begin{equation}
\Upsilon \big((x,m),(x',m'),\cdot\big) = \delta_{\frac{mx+m^\prime x^\prime}{m+m^\prime}}.
\end{equation}
Recall from Remark~\ref{Rem-placement} that this kernel  determines the new location of the particle in such a way that the center of mass is preserved. This implies that the center of mass of a component (the convex combination of its vertices) is equal to the location of the particle in the coagulation process, for any component (respectively particle), at any time. 

Since $\kappa$ is linear in each of its arguments
$$
|I| \,|I'|\, \kappa\Big(\frac 1{|I|} \sum_{i\in I} x_i,\frac 1{|I'|}\sum_{i\in I'} x_{i}\Big)
= \sum_{i\in I} \sum_{i\in I'}\kappa(x_i,x_{i}),\qquad I,I'\subset [N].
$$
Assuming that $I$ and $I'$ are disjoint, on the right-hand side, we see the rate (in the graph model) of putting a bond between any two vertices of the groups $\{x_i\colon i\in I\}$ and $\{x_i\colon i\in I'\}$, which corresponds, in the coagulation model, to the rate of the coagulation between two particles at the locations $\frac 1{|I|} \sum_{i\in I} x_i$ and $\frac 1{|I'|}\sum_{i\in I'} x_{i}$ with  mass $|I|$ respectively $|I'|$, that we see on the left-hand side. 


The rate function that we derived in \cite{AKLP23} shows how the macroscopic and the mesoscopic part of the configuration influence the large deviations of the microscopic part.
We strongly expect to see corresponding effects for the more general coagulation process studied in this work, but this issue is deferred to future work.
Our present work introduces a major improvement to the proof technique from \cite{AKLP23} by introducing in Theorem~\ref{thm-distMi} a Poisson point process that enables us to prove large-deviations results without the projective limits that are at the core of the proof in \cite{AKLP23}.

\section{Proof of Theorem~\ref{thm-distMi}: distribution of the configuration}\label{MiMaintroduction}

\noindent In this section, we prove our first main result, Theorem~\ref{thm-distMi}. That is, given an arbitrary probability measure $\mu$ on $\Scal$, we identify the distribution of the empirical measure $\Vcal_{N}^{\ssup T}$ defined in \eqref{empmeas} for any fixed $T\in(0,\infty)$ and $N\in\N$ under the poissonised initial condition $\P_{\Poi_{N\mu}}$ given in \eqref{PoissonInitial}. We assume that the coagulation kernel $K$ from Section~\ref{sec-model} is just measurable.

The main idea is the decomposition of the process  into all the subprocesses that coagulate into one single particle by time $T$. This is done via a decomposition in terms of all the coagulation trajectories, which we will also call {\em coagulation trees}, i.e., the parts of the coagulation process that coagulate into one single particle by time $T$. We are going to rewrite the joint distribution of all these trees in terms of  a self-interacting Poisson point process. 

Recall the partition process $Z= (Z_t)_{t\in[0,T]}$ from definition \eqref{Zprocess}, i.e., $Z_t = (X_C^\ssup{t}, C)_{C\in P_t}$, which consists of a process of partitions $P= (P_t)_{t\in [0,T]}\in \Pcal(A)$ for some $A\subset\N$ with a location $X^\ssup{t}_C\in\Scal$ for every particle $C\in P_t$ at every $t\in [0,T]$. The  Markovian dynamics is given in \eqref{mechanism1} and the initial state is $Z_0=(x_i,\{i\})_{i\in A}$ for some vector $\mathbf x\in \Scal^A$.  We denote the state space of $Z$ by $\Zcal_A$ and we indicate with $\P_{\mathbf x}$ the law of such process. For two disjoint finite sets $A,B\subset C$, we  denote by $\{A\nleftrightarrow B\}\subset \Zcal_{C}$ the event that there is no coagulation between any subset of $A$ and any subset of $B$ up to time $T$ in a coagulation process that starts from some superset $C$.

\begin{lemma}[Identification of $R^{\ssup{T}}$]\label{lem-Rident}
For any $\xi,\xi'\in \Gamma_T$, the density $R^{\ssup T}(\xi,\xi')$ in \eqref{DefR} exists, and the formula \eqref{Rdef} holds.
\end{lemma}
From formula \eqref{Rdef} one sees that $R^{\ssup{T}}$ does not depend on the placement kernel $\Upsilon$ and is symmetric, since $K$ is. Since for any $t\in [0,T]$ we have that $\xi_t$ and $\xi'_t$ are discrete in space, the integrals $\langle \xi_t, K\xi'_t\rangle$ are indeed sums.

\begin{proof}As for the definition \eqref{DefR} for $R^{\ssup T}(\xi,\xi')$, let $A$ and $B$ be disjoint finite sets and fix $\mathbf{x} = (x_i)_{i\in A}\in \Scal^A$ and $\mathbf{y} = (y_i)_{i\in B}\in \Scal^B$ such that $\xi_0=\sum_i\delta_{x_i}$ and $\xi'_0=\sum_i\delta_{y_i}$. By $\P_{(\mathbf{x},\mathbf{y})}$ we denote the distribution of the process $Z$ as above with partition process in $\Pcal(A\cup B)$ and initial configuration $(\mathbf x,\mathbf y)\in\Scal ^{A\cup B }$. 
On the event $\{A\nleftrightarrow B\}$, for any $t\in [0,T]$, $P_t$ can be decomposed into a partition of $A$ and a partition of $B$. In particular, we can decompose any state $z\in \{A\nleftrightarrow B\}$ of $Z$ as $z = (z_A,z_B)$, where $z_A\in \Zcal_A$, $z_B\in \Zcal_B$. On the other hand, to any pair $z = (z_A,z_B)$ with $z_A\in \Zcal_A$, $z_B\in \Zcal_B$ we can associate a state  $z\in \{A\nleftrightarrow B\}$.  Then
\begin{multline}\label{eq:conditioning}
\int_{\Zcal_{A\cup B}} \P_{(\mathbf x, \mathbf y)}(A \nleftrightarrow B, Z_A\in \d z_A, Z_B \in \d z_B) =\\
 \int_{\Zcal_A \times \Zcal_B} \P_{\mathbf{x}}(Z \in \d z_A) \P_{\mathbf{y}}(Z \in \d z_B)  p(z_A,z_B)
\end{multline}
where $p(z_A,z_B)=\P_{(\mathbf x, \mathbf y)}(z_A\nleftrightarrow z_B)$ is equal to the probability that no exponential time elapsed that would result in a coagulation between subsets of $A$ and subsets of $B$. More precisely,
given any two states $z_A\in \Zcal_A$, $z_B\in \Zcal_B$ we can iteratively, for $i=1, 2, \ldots$ construct maximal intervals $I_i = [t_{i-1}, t_i)$ such that the corresponding partition processes $(P^\ssup{A}_t)_{t\in [0,T]}$ and $(P_t^\ssup{B})_{t\in [0,T]}$ are both constant on $[t_{i-1}, t_i)$. Then, we indicate with $\{z_A\nleftrightarrow z_B\}$ the event that for all $i$ and for all $C\in P^\ssup{A}_{t_{i-1}}$, $D\in P^\ssup{B}_{t_{i-1}}$ the exponentially distributed times with parameter $K((X_C, |C|), (X_D, |D|))$ are larger that $t_i-t_{i-1}$. Since for any $i$ we can get the number of sets/particles $C\in P^\ssup{A}_{t_{i-1}}$ sitting in $(x,m)\in \Scal\times \N$ via the empirical measure $\Xi_{t_{i-1}}(z_A)(\d(x,m))$ as defined in \eqref{XiZ} (and the same for $B$) we have that  
\begin{align*}
p(z_A,z_B) = \prod_i \e^{-(t_i-t_{i-1})\langle \Xi_{t_{i-1}}(z_A),K \Xi_{t_{i-1}}(z_B)\rangle} = \exp\Big(-\int_0^T \langle \Xi_t(z_A), K \Xi_t(z_B)\rangle \d t \Big),
\end{align*}
which only depends on $\Xi(z_A), \Xi(z_B)$. Given that these exponential clocks did not ring, by definition the exponential clocks running between elements of $\Pcal(A)$ are independent from the ones running  between elements of $\Pcal(B)$. Therefore we have
\[\P_{(\mathbf x, \mathbf y)}(Z_A\in \d z_A, Z_B \in \d z_B|z_A\nleftrightarrow z_B)= \P_{\mathbf{x}}(Z \in \d z_A) \P_{\mathbf{y}}(Z \in \d z_B), \]
which justifies \eqref{eq:conditioning}. 
Now going to the distributions of $\Xi(z_A), \Xi(z_B)$ under $\P_{\mathbf{x}}$ and $\P_{\mathbf{y}}$ gives the result. 
\end{proof}


Now we express the probability of the event that the coagulation process decomposes into two pieces, not necessarily terminating with precisely one particle at time $T$.

\begin{lemma}\label{lem:Z-decomp}
Let $A$ and $B$ be disjoint finite sets and fix $\mathbf{x} = (x_i)_{i\in A}\in \Scal^A$ and $\mathbf{y} = (y_i)_{i\in B}\in \Scal^B$. By $\P_{(\mathbf{x},\mathbf{y})}$ we denote the distribution of the process $Z$ as above with partition process in $\Pcal(A\cup B)$ and initial configuration $(\mathbf x,\mathbf y)\in\Scal ^{A\cup B }$. Then, for any measurable $\Ncal\subset \Mcal_{\N_0}(\Gamma_T^\ssup{1})$,
\begin{multline}\label{eq:Z-decomp}
\P_{(\mathbf{x}, \mathbf{y})}\big(A\nleftrightarrow B, N\Vcal_N^\ssup{T} \in \Ncal\big) 
=\\
 \Big\langle \P_{\mathbf{x}}\otimes \P_{\mathbf{y}}, \e^{ -R^{\ssup T}(\Xi^\ssup{T,A}, \Xi^\ssup{T,B}) } \1\big\{ N[\Vcal_N^\ssup{T}(Z_A)+ \Vcal_N^\ssup{T}(Z_B)] \in \Ncal\big\}\Big\rangle,
\end{multline}
where we recall definition \eqref{history} and  \eqref{empmeas}. (We write the random processes under $\P_{\mathbf{x}}$ and $\P_{\mathbf{y}}$ as $Z_A$ respectively $Z_B$ and used \eqref{empmeas} for each separately.)
\end{lemma}

\begin{proof}
On the event $\{A\nleftrightarrow B\}$, for any $t\in [0,T]$, $P_t$ can be decomposed into a partition of $A$ and a partition of $B$. In particular, we can decompose any state $z\in \{A\nleftrightarrow B\}$ of $Z$ as $z = (z_A,z_B)$, where $z_A\in \Zcal_A$, $z_B\in \Zcal_B$. The measure $\Vcal_N^\ssup{T}$ decomposes accordingly into the sum of $ \Vcal_N^\ssup{T}(Z_A)$ and $\Vcal_N^\ssup{T}(Z_B) $.  This gives us 
\begin{align*}
\P_{\mathbf{x}, \mathbf{y}}\big(&A\nleftrightarrow B, N\Vcal_N^\ssup{T} \in \Ncal\big)=\\
  &\int_{\Zcal_A}\int_{\Zcal_B} \frac{\d \P_{\mathbf{x}, \mathbf{y}}(A\nleftrightarrow B, (Z_A, Z_B) \in \cdot)}{\d \P_{\mathbf{x}}\otimes \P_{\mathbf{y}}}(z_A,z_B) \1\{N\Vcal_N^\ssup{T} \in \Ncal\}\,\P_{\mathbf{x}}(\d z_A) \P_{\mathbf{y}}(\d z_B).
\end{align*} 
Hence, we can insert the density which is equal to $\e^{-R^{\ssup T}(\Xi^\ssup{T,A}(z_A), \Xi^\ssup{T,B}(z_B))}$ and arrive at \eqref{eq:Z-decomp}.
\end{proof}

Now we formulate the decomposition of the coagulation process into pieces that end up with just one particle each at time $T$. We denote the restriction of $\P_{\mathbf x}$ to $\{\Xi(Z)\in\Gamma_T ^{\ssup 1}\}$ by $\P^{\ssup 1}_{\mathbf x}$ and note that this is a sub-probability measure.

\begin{lemma}\label{lem:partition-decomp}
Fix $M\in\N$ and  let $P=\{C_j\colon j\in[m]\}\in \Pcal([M])$. Fix $\mathbf{x} = (x_i)_{i=1}^M\in \Scal^M$ and denote $\mathbf{x}^\ssup{j} = (x_i)_{i\in C_j}$. Then, for any measurable $\Ncal\subset \Mcal_{\N_0}(\Gamma_T^\ssup{1})$ 
\begin{equation}\label{eq:partition-decomp}
\P_{\mathbf{x}}\Big(P_T = P, N\Vcal_N^\ssup{T} \in \Ncal\Big) 
= \Big\langle \bigotimes_{j=1}^m \P^\ssup{1}_{\mathbf{x}^\ssup{j}}, \e^{ -\frac{1}{2}\sum_{j,j' \colon j\neq j'}R^{\ssup T}(\Xi_j, \Xi_{j'}) } \1\big\{
\textstyle\sum_{j=1}^m \delta_{\Xi_j} \in \Ncal\big\}\Big\rangle 
\end{equation}
where for each $j=1, \ldots, m$, we denote the random variables under $\P_{\mathbf{x}^\ssup{j}}$ by $Z_j$ and put $\Xi_j=\Xi(Z_j)$.
\end{lemma}

\begin{proof}
Note that $\{P_T = P\}$ implies that $\{C_1 \nleftrightarrow \bigcup_{j=2}^m C_j\}$ and use Lemma \ref{lem:Z-decomp}. Iterating this argument gives formula \eqref{eq:partition-decomp}, since $\sum_{C,\tilde C \in P\colon C \neq \tilde C}R^\ssup{T}(\Xi^\ssup{T,C},\Xi^\ssup{T,\tilde C}) = \frac{1}{2}\sum_{j,j'\colon j\neq j'}R^\ssup{T}(\Xi_j, \Xi_{j'})$ 
\end{proof}

\begin{proof}[Proof of Theorem \ref{thm-distMi}]
Without loss of generality it suffices to show the statement of the Theorem for  functions of the form $f(\nu)= \1_{\Ncal}(N\nu)$ for some measurable set $\Ncal\in \Mcal_{\N_0}(\Gamma^\ssup{1}_T)$. Denote by $L$ the set of collections of numbers $\ell = (\ell_n)_{n\in \N}\in \N_0^{\N}$. For any point measure $\nu= \sum_{i\in I} \delta_{\xi^\ssup{i}} \in \Mcal_{\N_0}(\Gamma^\ssup{1}_T)$ there exists an $\ell \in L$ such that $\ell_n = \#\{i\in I \colon |\xi^\ssup{i}_0| = n \}$, for any $n\in\N$, i.e., $\nu$ consists of exactly $\ell_n$ (coagulation) trees of of size $n$, for each $n\in \N$ (where the size of a tree $\xi\in \Gamma^\ssup{1}_T$ is given as the number of atoms, i.e. by $|\xi_0|$). In that case we say that the tree sizes of $\nu$ are given by $\ell$. We can decompose any measurable set $\Ncal\subset \Mcal_{\N_0}(\Gamma^\ssup{1}_T)$ as $\Ncal = \bigcup_{\ell \in L} \Ncal(\ell)$ with 
\[
\Ncal(\ell) = \{\nu \in \Ncal \colon \text{the tree sizes of }\nu \text{ are given by }\ell\}, \quad \text{ for any } \ell \in L,
\]
where the sets $\Ncal(\ell)$, for $\ell \in L$, are disjoint. In the following we will assume  that $\Ncal = \Ncal(\ell)$ for some $\ell\in L$. 
Recall that $N\Vcal_N^\ssup{T} =  \sum \delta_{\Xi^\ssup{T,C}}$, where the sum extends over $C\in P_T$, which we will leave out in the notation. 
We want study the event $\{\sum \delta_{\Xi^\ssup{T,C}} \in \Ncal \}$.  We abbreviate $M = \sum_n \ell_n n$ and $m= \sum_n \ell_n$ and note that on the event $\{\sum \delta_{\Xi^\ssup{T,C}} \in \Ncal \}$ it holds that $M$ is the total number of atoms and $m$ is the total number of coagulation trees. Recall that under $\P_{\Poi_ {N\mu}}$ we want to consider the empirical measure $\Vcal_N^\ssup{T}$ under $\P_{\mathbf{x}}$, where the initial condition $\mathbf{x}=(x_i)_{i\in I}$ is such that $\sum_{i\in I}\delta_{x_i}$ is $\Poi_{N\mu}$ distributed. More precisely, we choose the length  $|I|$ of the initial vector $\mathbf{x}$ as $\Poi_N$-distributed and sample the $(x_i)_{i\in I}$ i.i.d. and with distribution $\mu$. Therefore, we have that
\begin{multline*}
\P_{\Poi_{N\mu}}\big(\sum \delta_{\Xi^\ssup{T,C}} \in \Ncal \big) =\\
 \sum_{M' =0}^\infty \Poi_N(M') \int_{\Scal^{M'}}\mu^{\otimes M'}(\d(x_1, \ldots, x_{M'})) \,\P_{\mathbf{x}} \big( \sum \delta_{\Xi^\ssup{T,C}} \in \Ncal\big).
\end{multline*} 
The sum reduces to the summand $M' = M$, since otherwise the probability of the event is $0$.
Now, 
\begin{align*}
 \P_\mathbf{x}\big({\textstyle \sum \delta_{\Xi^\ssup{T,C}} \in \Ncal }\big)  = \sum_{\substack{P \in \Pcal([M]) \\ P \text{ compatible with }\ell}}\P_\mathbf{x}\big(P_T= P, \,{\textstyle\sum \delta_{\Xi^\ssup{T,C}} \in \Ncal }\big),
\end{align*}
where we say that a partition $P\in \Pcal([M])$ is compatible with $\ell$ if for each $n\in \N$ we have that $\#\{C\in P\colon |C|=n\} = \ell_n$. Using the formula from Lemma \ref{lem:partition-decomp} we get that 
\begin{equation}\label{ProofStep1}
\begin{aligned}
&\P_{\Poi_{N\mu}}\big(\sum \delta_{\Xi^\ssup{T,C}} \in \Ncal \big) = \Poi_N(M) \\
& \quad\times \sum_{\substack{P \in \Pcal(M) \\ P \text{ compatible with }\ell}}  \int_{\Scal^{M}}\mu^{\otimes M}(\d\mathbf x)\, \Big\langle \bigotimes_{j=1}^m \P^\ssup{1}_{\mathbf{x}^\ssup{j}}, \e^{ -\frac{1}{2}\sum_{j,j'\colon j\neq j'}R(\Xi_j, \Xi_{j'}) } \1\big\{\textstyle\sum_{j=1}^m \delta_{\Xi_j} \in \Ncal \big\}\Big\rangle.
\end{aligned}
\end{equation}
For $n\in \N$ we denote $\mu^{\otimes n} \otimes \P^\ssup{1}(\d (\mathbf x,Z)) = \mu^{\otimes n}(\d \mathbf x) \otimes \P^\ssup{1}_{\mathbf x}(\d Z)$ (i.e., we conceive $\P^{\ssup 1}$ as the kernel $(\mathbf x,A)\mapsto \P^\ssup{1}_{\mathbf x}(A)$). Then for any $P\in \Pcal(M)$ that is compatible with $\ell$ we have that
{ \begin{equation}\label{eq:Pell}
\begin{aligned}
\begin{array}{ll}
\text{each summand}\\
\text{in 2nd line of \eqref{ProofStep1}}
\end{array}= \Big\langle \bigotimes_{j=1}^m \Big(\mu^{\otimes |C_j|}\otimes \P^\ssup{1}\Big), \e^{ -\frac{1}{2}\sum_{j,j'\colon j\neq j'}R^{\ssup T}(\Xi_j, \Xi_{j'}) } \1\big\{\textstyle\sum_{j=1}^m \delta_{\Xi_j} \in \Ncal \big\}\Big\rangle,
 \end{aligned}
\end{equation}}where we put  $P=\{C_j\colon j\in[m]\}$. The right-hand side of \eqref{eq:Pell} depends only on the cardinalities $n_j$ of the partition sets $C_j$. Given $\ell$ we can uniquely fix a collection of numbers $n_1,\dots,n_m\in\N$ such that $\#\{j\colon n_j = n\} = \ell_n$ for each $n\in \N$. Hence, we need to find the number of partitions $P=\{C_j\colon j\in[m]\}$ such that the cardinalities $n_j$ of the partition sets satisfy this, in other words, such that $P$ is compatible with $\ell$.
Indeed, we observe that
\begin{equation}\label{eq:countingpartitions}
\#\big\{P\in \Pcal(M) \colon P \text{ compatible with } \ell\big\}  = \frac{M!}{m!\prod_{j=1}^m n_j!}= \frac{M!}{m!\prod_{n} n!^{\ell_n}}.
\end{equation}
To see that recall that the multinomial factor $\frac{M!}{\prod_{j=1}^m n_j!}$ is the number of possibilities of putting the indices from $[M]$ into boxes labelled by $j= 1, \ldots, m$, such that box $j$ has exactly $n_j$ many indices. Each of those possibilities gives us an ordered partition of $[M]$ that is compatible with $(n_j)_{j=1}^m$ (in the sense that its $j$-th set has $n_j$ many elements). Any (non-ordered) partition of $[M]$ that is compatible with $\ell$ corresponds to precisely $m!$ many ordered partitions that are compatible with $(n_j)_{j=1}^m$. Thus, formula \eqref{eq:countingpartitions} holds.

Summarizing the previous steps, we have shown that
\begin{equation}\label{poiiniteq2}
\begin{aligned}
&\P_{\Poi_{N\mu}}\big(\sum \delta_{\Xi^\ssup{C}} \in \Ncal \big)
= \Poi_N(M)\,  \frac{1}{m!}\, \frac{M!}{\prod_{j=1}^m n_j!}\, \\
& \hspace{3cm} \Big\langle \bigotimes_{j=1}^m \big(\mu^{\otimes n_j}\otimes\P^\ssup{1}\big), \e^{-\frac{1}{2}\sum_{j,j'\colon j\neq j^\prime}R^{\ssup T}(\Xi_j,\Xi_{j^\prime})} \1\big\{\textstyle\sum_j \delta_{\Xi_j}\in \Ncal\big\}\Big\rangle,
\end{aligned}
\end{equation}
where we recall that $\Ncal=\Ncal(\ell)$ and $\#\{j\colon n_j = n\} = \ell_n$ for each $n\in \N$ and that $\Xi_j$ is short for $\Xi(Z_j)$ (see \eqref{XiZ}), where $Z_j$ is the random                                                       variable under $\mu^{\otimes n_j}\otimes\P^\ssup{1}$. 

Now we make the connection with the reference measure $M_{b\mu,N}^{\ssup T}$ defined in \eqref{Mbreferencemeasure}.
Note that on $\{\xi\in \Gamma^\ssup{1}_T \colon |\xi_0| = n_j\}$ we have 
\[
\frac{1}{n_j!}\,\Big(\mu^{\otimes n_j}\otimes N^{n_j-1}\P^\ssup{1}\Big)\circ\Xi^{-1}(\d \xi) =  b^{-n_j}\e^{b-1} M_{b\mu,N}^{\ssup T}(\d \xi),\qquad b\in(0,\infty),
\]
where we conceive $\Xi$ as the map defined in \eqref{XiZ}.

Inserting this into \eqref{poiiniteq2} and using that $\sum_j n_j = M$ (and hence $\prod_j N^{n_j-1}=N^{M-m}$) we get
\begin{align*}
&\P^{\Poi_{N\mu}}\big(\sum \delta_{\Xi^\ssup{C}} \in \Ncal \big)\\
&\quad =  \Poi_N(m) \, \Big\langle \bigotimes_{j=1}^m \Big(\frac1{n_j!}\mu^{\otimes n_j}\otimes N^{n_j-1}\P^\ssup{1}\Big), \e^{-\frac{1}{2}\sum_{j,j'\colon j\neq j^\prime}R^{\ssup T}(\Xi_j,\Xi_{j^\prime})} \1\big\{\textstyle\sum_j \delta_{\Xi_j}\in \Ncal\big\}\Big\rangle\\
&\quad =  \Poi_N(m)\,\e^{(b-1)m}\,b^{-M} \Big\langle (M_{b\mu,N}^{\ssup T})^{\otimes m},\e^{-\frac{1}{2}\sum_{j,j'\colon j\neq j^\prime}R^{\ssup T}(\Xi_j,\Xi_{j^\prime})} \1\big\{\textstyle\sum_j \delta_{\Xi_j}\in \Ncal\big\}\Big\rangle\\
&\quad = \e^{(b-1)m}\,b^{-M}\, \e^{N(|M^{\ssup T}_{b\mu,N}|-1)}\,{\tt E}_{NM^{\ssup T}_{b\mu,N}}\Big[\e^{-\frac{1}{2}\sum_{j,j'\colon j\neq j^\prime}R^{\ssup T}(\Xi_j,\Xi_{j^\prime})} \1\big\{\textstyle\sum_j \delta_{\Xi_j}\in\Ncal\big\}\Big],
\end{align*}
where the $\Xi_j$ denote in the second line $\Xi(Z_j)$, in the third line the random variables under $M_{b\mu,N}^{\ssup T}$ (even though this is not normalised) and in the fourth line the points of the PPP.

We have arrived at the assertion.
\end{proof}

\section{Preparations}\label{sec-Prep}

\noindent In this section we prepare for the proof of our remaining main  results by doing the following.
\begin{itemize}
\item We make some technical remarks about topologies and metrics and prove that $\Gamma_T ^{\ssup 1}$ is closed in Section~\ref{sec-topology}.

\item We investigate properties of the non-coagulation operator $\mathfrak R^{\ssup T}$ in Section~\ref{sec-noncoagrates}.

\item We prove the convergence of $\Q_k^{\ssup{T,N,N}}$ and introduce the measure  $\Q_k^\ssup{T}$ in Section~\ref{sec-Qproof}.

\item We prove the continuity of the map $\nu\mapsto \rho(\nu)$ in Section~\ref{sec-proofrhocont}.

\end{itemize}

\subsection{Topologies and metrics}\label{sec-topology}

For our further developments, we need to make some comments on the topologies used. We need to do this on various levels. Recall that we assumed $\Scal$ to be a Polish space. Let $d$ be a metric on $\Scal\times\N$ such that $(\Scal\times \N,d)$ is a complete space.

On the space $\Mcal(\Scal \times \N)$ we consider the weak topology, i.e., $\phi^\ssup{n} \to \phi$ as $n\to \infty$ if $\langle \phi^\ssup{n}, f\rangle \to \langle \phi, f\rangle$ as $n\to \infty$ for any continuous and bounded test functional $f\in C_{\rm b}(\Scal \times \N)$. Note that $\Mcal(\Scal\times \N)$ is  Polish.
Weak convergence on $\Mcal(\Scal\times \N)$ is induced by the  L{\'e}vy-Prohorov metric ${\mathfrak d}$, which is defined as follows, see \cite{Bil13}. For all $\phi,\phi'\in \Mcal(\Scal \times \N)$,
\begin{equation}\label{Prorohovmetric}
{\mathfrak d}(\phi,  \phi^\prime ) = \inf\{\eps>0\colon \phi(A)\leq \phi'(A^{\eps})+\eps,\phi'(A)\leq \phi(A^{\eps})+\eps, \text{ for any mb.\ } A\subset\Scal\times\N\},
\end{equation}
where $A^{\eps}=\{x\in \Scal\times \N\colon {d}(x,y)<\eps\text{ for some } y\in A\}$ denotes the open $\eps$-neighbourhood of $A\subset\Scal\times\N$.

Eventually, we want to consider paths of measures. Denote by $\Mcal$ either $\Mcal(\Scal\times \N)$ or $\Mcal_{\N_0}(\Scal\times \N)$ (the space of point measure on $\Scal\times \N$) and equip it with the metric $\mathfrak{d}$. Then $\D_T=\D_T(\Mcal)$ denotes the space of c\`adl\`ag functions $[0,T]\to\Mcal $. We endow $\D_T$ with the Skorohod $J_1$-topology, which is induced by a certain metric ${\rm d}_T$ such that the space $(\D_T, {\rm d}_T)$ is separable and complete (see \cite{Kal97}, Thm.~A2.2), i.e., it is a Polish space. Convergence in $(\D_T, {\rm d}_T)$ can be characterised via time-changes on $[0,T]$, which are  strictly increasing bijections $\lambda \colon [0,T]\to [0,T]$ (which are necessarily continuous  with $\lambda(0) = 0$ and $\lambda(T) = T$) as follows. For $ \xi,\xi^\ssup1,\xi^\ssup2,\dots \in \D_T$ it holds that $\xi^\ssup{n} \to \xi$ as $n\to \infty$ if and only if
\begin{equation}\label{tree-conv-timechange}
\sup_{t\in[0,T]} |\lambda_n(t) - t| + \sup_{t\in[0,T]} \mathfrak d\big(\xi^\ssup{n}_{\lambda_n(t)},\xi_t\big) \to 0 \text{ as } n\to \infty,
\end{equation}
for some time-changes $\lambda_n$ on $[0,T]$.

Recall that the set of one-particle coagulation trajectories is 
\[
\begin{aligned}
\Gamma^\ssup{1}_T = \big\{ \xi = &\,(\xi_t)_{t\in[0,T]} \in \D_T(\Mcal_{\N_0}(\Scal\times \N)) \colon \xi_0 \text{ is concentrated on } \Scal\times \{1\},\\
&  t\mapsto\xi_t \text{ is piecewise constant and makes steps as in } \eqref{steps} \text{ and }\xi_T(\Scal\times \N) =1\big\},
\end{aligned}
\]
Since $\D_T(\Mcal_{\N_0}(\Scal\times \N))$ is separable, the same is true for $(\Gamma^\ssup{1}_T,{\rm d}_T)$. In Lemma \ref{lem-Gammacomplete} we will show that $(\Gamma^\ssup{1}_T,{\rm d}_T)$ is also closed and hence complete as a closed subset of the complete space $\D_T(\Mcal_{\N_0}(\Scal\times \N))$. With other words, it is itself a Polish space.

Recall that the state space of our process $\Vcal_N^\ssup{T}$ is equal to the set $\Mcal(\Gamma^\ssup{1}_T)$ of positive measures on $\Gamma^\ssup{1}_T$. We equip $\Mcal(\Gamma^\ssup{1}_T)$ with the weak topology, i.e., $\nu^\ssup{n} \to \nu$ as $n\to \infty$ if $\langle \nu^\ssup{n}, f\rangle \to \langle \nu, f\rangle$ as $n\to \infty$ for any continuous and bounded test functional $f\in C_{\rm b}(\Gamma^\ssup{1}_T)$.

Now we show that the one-particle coagulation trajectory set $\Gamma_T^{\ssup 1}$ defined in \eqref{Gammadef} and \eqref{Gamma1def} is closed. 

\begin{lemma}\label{lem-Gammacomplete}
The set $\Gamma^\ssup{1}_T$ is a closed subset of $\D_T(\Mcal_{\N_0}(\Scal\times \N))$. Consequently, $(\Gamma^\ssup{1}_T, {\rm d}_T)$ is a Polish space.
\end{lemma}

\begin{proof}
We abbreviate $\D_T$ for $\D_T(\Mcal_{\N_0}(\Scal\times\N))$. Let $(\xi^\ssup{n})_{n\in \N}$ be a sequence in $\Gamma^\ssup{1}_T$ that converges to $\xi\in \D_T$ as $n\to \infty$. We will argue that $\xi\in \Gamma^\ssup{1}_T$. There exist time changes $\lambda_n$, $n\in \N$, such that \eqref{tree-conv-timechange} holds. First, we need to show that $\xi_0$ is concentrated on $\Scal\times\{1\}$ and that $\xi_T(\Scal\times \N)=1$.  But this is easy. Indeed, since $\lambda_n(0)=0$ and $\lambda_n(T)=T$, we have that $\xi_0 = \lim_{n\to \infty} \xi^\ssup{n}_0$ is concentrated on $\Scal\times \{1\}$ and $\xi_T(\Scal\times \N) = \lim_{n\to \infty} \xi^\ssup{n}_T(\Scal\times \N) =1$. 

Now we need to show that $\xi$ is piecewise constant and makes only steps as in \eqref{steps}. For this we need a uniform lower bound for the holding times of $\xi^\ssup{n}$, $n\in \N$. Otherwise one might have a succession of jumps that happen increasingly fast and add up in the limit to a larger jump that is no longer of the form \eqref{steps}. For any $\tilde\xi\in \D_T$ let $\tilde t_1 < \tilde t_2 < \ldots$ be the discontinuity points of $\tilde \xi$ and  define 
\[
s_*(\tilde\xi) = \inf_{j\geq 1}  (\tilde t_{j+1} -\tilde t_{j}) \wedge T.
\]
If $\tilde \xi$ is piecewise constant (which is the case for $\tilde \xi \in \Gamma^\ssup{1}_T$) then $s_*(\tilde \xi)$ is the infimum over all holding times of $\tilde\xi$. We also define the modulus of continuity for $\tilde \xi \in \D_T$ with spacing $h>0$ as
\begin{equation}\label{modulcont}
w(\tilde\xi,h) = \inf_{(I_i)_i}\, \max_i \,\sup_{u,v\in I_i} \mathfrak d(\tilde\xi_u,\tilde\xi_v),
\end{equation}
where the infimum extends over all partitions of the interval $[0,T]$ into subintervals $I_i = [s_i,t_i)$ such that $|I_i| := t_i-s_i>h$ for all $i$. 
Note that if $\tilde\xi$ is piecewise constant, then $w(\tilde\xi, h) = 0$ for all $h< s_*(\tilde\xi)$. 

By Theorem A2.2 in \cite{Kal97} the convergence of $\xi^{\ssup n}$ towards $\xi$ implies that  \[\lim_{h\to 0}\sup_{n\in \N} w(\xi^\ssup{n},h)=0.\] We want to show that $\inf_{n\in\N} s_*(\xi^\ssup{n})>0$. We argue via contraposition, i.e., we show that $\inf_{n\in\N} s_*(\xi^\ssup{n})=0$ implies that $\lim_{h\to 0}\, \sup_{n\in\N} \,w(\xi^\ssup{n},h)>0$. Fix any $h>0$ and let $n$ be such that $s_*(\xi^\ssup{n}) <h$. Then any partition $(I_i)_i$ of $[0,T]$ satisfying $|I_i| >h$ for all $i$ contains an interval $I_j$ that contains a discontinuity point $t_n$ of $\xi^\ssup{n}$ and thus
\[
\sup_{u,v\in I_j} \mathfrak d(\xi^\ssup{n}_u,\xi^\ssup{n}_v) \geq  \xi^\ssup{n}_{t_n-}(\Scal\times \N) -  \xi^\ssup{n}_{t_n}(\Scal\times \N) =1.
\]
This implies that $\lim_{h\to 0}\, \sup_{n\in\N} \,w(\xi^\ssup{n},h)\geq 1>0$.
Hence, we have shown that \[\inf_{n\in\N} s_*(\xi^\ssup{n})>0.\]

Now, we argue that $\xi$ is piecewise constant. Fix $h>0$ with $h< \inf_{n\in\N} s_*(\xi^\ssup{n})$ and note that $\sup_n w(\xi^\ssup{n},h) =0$. Using the characterization of convergence from \eqref{tree-conv-timechange} one can show that for  any $\eps\in(0,h)$ one has that $w(\xi, h) \leq w(\xi^\ssup{n},h-\eps) + \eps = \eps$, if $n$ is large enough. Hence $w(\xi, h)=0$, which implies that $\xi$ is piecewise constant.

Now, we show that $\xi$ makes steps as in \eqref{steps}. Fix a discontinuity point $t$ of $\xi$. Take $\eps>0$ small enough such that $3 \eps < \inf_{n\in\N} s_*(\xi^\ssup{n})\wedge t$. We have that $\xi_t \neq \xi_{t-\eps}$ and
\[
\xi_t -\xi_{t-\eps} = \lim_{n\to \infty}\big( \xi^\ssup{n}_{\lambda_n(t)} - \xi^\ssup{n}_{\lambda_n(t-\eps)}\big).
\]
We choose $n$ large enough such that $\mathfrak d(\xi^\ssup{n}_{\lambda_n(t-\eps)}, \xi^\ssup{n}_{\lambda_n(t)}) >0$ and $\sup_{u}|\lambda_n(u) - u|\leq \eps$. Then $\mathfrak d(\xi^\ssup{n}_{\lambda_n(t-\eps)}, \xi^\ssup{n}_{\lambda_n(t)}) >0$ implies that the interval $[\lambda_n(t-\eps), \lambda_n(t))$ contains at least one discontinuity point of $\xi^\ssup{n}$. On the other hand we have that ${\lambda_n(t)- \lambda_n(t-\eps)\leq 3\eps < \inf_{n\in\N} s_*(\xi^\ssup{n})}$ and this implies that the interval $[\lambda_n(t-\eps), \lambda_n(t))$ contains at most one discontinuity point of $\xi^\ssup{n}$. This gives us that 
\[
\xi_t -\xi_{t-\eps} = \lim_{n\to \infty}\big( \xi^\ssup{n}_{\lambda_n(t)} - \xi^\ssup{n}_{\lambda_n(t-\eps)}\big) = \lim_{n\to \infty}\big( -\delta_{(x_n,m_n)} -\delta_{(x'_n,m'_n)} + \delta_{(z_n,m_n+m'_n)}\big)
\]
for some $(x_n,m_n),(x'_n,m'_n)\in \Scal\times \N$ and $z_n\in \Scal$. 
It is not hard to argue that the left-hand side does not depend on $\eps$ and that the convergence of the measures on the right-hand side implies the convergence of the atoms and hence the right-hand side is equal to $-\delta_{(x,m)} -\delta_{(x',m')} + \delta_{(z,m+m')}$, which finishes the proof. 
\end{proof}

\subsection{The rates of non-coagulation}\label{sec-noncoagrates}

In this section, we study properties of the  operator $\mathfrak{R}^\ssup{T}$ defined in \eqref{ROperator} with kernel equal to the non-coagulation probability, $R^{\ssup{T}}$, defined in \eqref{DefR} and identified in \eqref{Rdef}.

\begin{lemma}[Properties of $\mathfrak R^{\ssup{T}}$]\label{Lem-RfrakProp}
Assume that $K\colon (\Scal\times\N)^2\to[0,\infty)$ is continuous and symmetric, and fix $T\in(0,\infty)$. Then the following holds: 
\begin{enumerate}
\item[(1)] The mapping $(\xi, \xi') \mapsto R^\ssup{T}(\xi, \xi')$, is continuous on $\Gamma^\ssup{1}_T\times \Gamma^\ssup{1}_T$.
\item[(2)]  The map $\nu\mapsto \langle \nu,\mathfrak R^{\ssup{T}}(\nu)\rangle$ is lower semicontinuous with respect to the weak topology.
\item[(3)]
Assume that $K$  satisfies \eqref{AssK1}, and fix a majorizing function $f$ satisfying $f(r)/r\to\infty$ as $r\to\infty$ and $f(r)\geq r$ for any $r\in\N$. Then  $\nu\mapsto \langle \nu,\mathfrak R^{\ssup{T}}(\nu)\rangle$ is bounded and continuous on $\Acal_{f,\Aold}$ for any $\Aold\in(0,\infty)$.
\end{enumerate}
\end{lemma}

\begin{proof} 

We start by showing (1) Recall the Skorohod $J_1$-topology introduced at the beginning of Subsection \ref{sec-topology}. Fix $\xi' \in \Gamma^\ssup{1}_T$ and let $\xi^\ssup{n}, \xi\in \Gamma^\ssup{1}_T$ be such that $\xi^\ssup{n}\to \xi$ and let $\lambda_n$ be  time-changes, such that \eqref{tree-conv-timechange} holds. We can assume that $|\xi^\ssup{n}_0| = |\xi_0|$ for all $n$. For $i=1, \ldots,|\xi'_0|-1$, let $\phi'_i$ be the value of the path $\xi'$ on $[t_{i-1},t_i)$. Continuity of $K$ implies that  the mappings $(x, m) \mapsto K\phi_i'(x,m)$ are continuous. They are also bounded on $\Scal \times \{1, \ldots, |\xi_0|\}$ and hence bounded on the support of $\bigcup_{n,t}\xi^\ssup{n}_{\lambda_n(t)}$. Hence, $\sup_{t\in [0,T]}\mathfrak{d}(\xi^\ssup{n}_{\lambda_n(t)},\xi_t) \to 0$, as $n\to \infty$, implies that 
\[
\lim_{n\to \infty}R^\ssup{T}(\xi^\ssup{n}_{\lambda_n}, \xi') = \lim_{n\to \infty}\sum_{i=1}^{|\xi'_0|-1}\int_{t_{i-1}}^{t_i}\langle \xi^\ssup{n}_{\lambda_n(t)}, K\phi'_i\rangle \d t = \sum_{i=1}^{|\xi'_0|-1}\int_{t_{i-1}}^{t_i}\langle \xi_{t}, K\phi'_i\rangle \d t= R^\ssup{T}(\xi, \xi').
\]
(The fact that we can control the continuity of the mappings $K\xi_t'$ uniformly for $t\in [0,T]$ is implied more generally by the fact that $\lim_{h\to 0} w(\xi',h) =0$, recalling \eqref{modulcont}). Also, using that the jumps are of the form \eqref{steps} one has that
\[
\big|R^\ssup{T}(\xi^\ssup{n}, \xi') - R^\ssup{T}(\xi^\ssup{n}_{\lambda_n}, \xi') \big| \leq 3T \sup_{i,x,m\colon m\leq |\xi_0|}K\phi_i' (x,m) \, \sup_{t\in [0,T]}|\lambda_n(t) -t|\to 0  \text{ as } n\to \infty 
\]
 Altogether we proved that $R^\ssup{T}(\xi^\ssup{n}, \xi')\to R^\ssup{T}(\xi, \xi')$, as $n\to \infty$. Continuity of the mapping $(\xi, \xi') \mapsto R^\ssup{T}(\xi, \xi')$ is then implied by symmetry.

We continue with a proof for (2).
If $\nu_n\to\nu$ as $n\to\infty$ weakly, then also $\nu_n\otimes \nu_n\to\nu\otimes\nu$. Using the continuity of $R^\ssup{T}$ that we established in point (1), Lemma~\ref{lem-Rident} and Fatou's lemma, one easily sees that then $\liminf_{n\to\infty}\langle \nu_n,\mathfrak R^{\ssup{T}}(\nu_n)\rangle\geq \langle \nu,\mathfrak R^{\ssup{T}}(\nu)\rangle$, which shows lower semicontinuity.

Let us now show (3). We show only the upper semi-continuity. 
Because of \eqref{AssK1}, we have for any $\xi,\xi'\in\Gamma^\ssup{1}_T$ the upper bound
\begin{equation}\label{Resti}
\begin{aligned}
R^{\ssup T}(\xi,\xi')&=\int_0^T\d t\, \langle \xi_t, K\xi_t'\rangle 
\leq H\int_0^T\d t\, \|\xi_t\|_1\,\|\xi_t'\|_1 = H T \|\xi_0\|_1\,\|\xi_0'\|_1=HT |\xi_0|\,|\xi'_0|,
\end{aligned}
\end{equation}
where we used that $t\mapsto \|\xi_t\|_1=\sum_{x,m}\xi_t(x,m)m$ is constant for any $\xi\in \Gamma^\ssup{1}_T$. For $t=0$, $\|\xi_0\|_1=|\xi_0|$ is the total mass of $\xi_0$. 

Hence, we may split, for any $L>0$,
\begin{equation}\label{eq:Rfracdecomp}
\begin{aligned}
\langle \nu,\mathfrak{R}^{\ssup{T}}(\nu)\rangle
&=\int_{|\xi_0|\leq L} \nu(\d \xi)\int_{|\xi'_0|\leq L}\nu(\d \xi')\, R^{\ssup T}(\xi,\xi')\\
&\qquad+\int \nu(\d \xi)\int\nu(\d \xi')\, R^{\ssup T}(\xi,\xi')\1\{|\xi_0|> L\mbox{ or }|\xi'_0|> L\}\\
&\leq \int_{|\xi_0|\leq L} \nu(\d \xi)\int_{|\xi'_0|\leq L}\nu(\d \xi')\, R^{\ssup T}(\xi,\xi')
+2H T |c_{\nu_0}|\,\int_{|\xi_0|> L}\nu(\d \xi)\,|\xi_0|,
\end{aligned}
\end{equation}
where we recall the definition of $c_\lambda$ from \eqref{clambdadef} and note that $|c_{\nu_0}|= \int\nu(\d \xi)\,|\xi_0|$ is bounded by $\Aold$ for all  $\nu\in\Acal_{f,\Aold}$.
Consider the last term on the right-hand side of \eqref{eq:Rfracdecomp}. Note that, for $r>L$, we can estimate $r\leq f(r)\eps_L$ with $\eps_L=\sup_{r>L}r/f(r)$, which vanishes as $L\to\infty$. Hence,
\[
\int_{|\xi_0|> L} \nu(\d \xi)\, |\xi_0| \leq \eps_L \int_{|\xi_0|> L} \nu(\d \xi)\, f(|\xi_0|)\leq \Aold\eps_L. 
\]
So, the last term on the right-hand side of \eqref{eq:Rfracdecomp} vanishes, as $L\to \infty$.  
Consider the first term on the right-hand side of \eqref{eq:Rfracdecomp}. Observe that due to the point (1) the function $(\xi,\xi') \mapsto R(\xi, \xi')\1\{|\xi_0|\leq L\}\1\{|\xi'_0|\leq L\}$ is dominated by a bounded and continuous function.
Together with the fact that $\nu\mapsto \nu\otimes \nu$ is continuous in the weak topology this implies that $\nu \mapsto \int_{|\xi_0|\leq L} \nu(\d \xi)\int_{|\xi'_0|\leq L}\nu(\d \xi')\, R(\xi,\xi')$ is continuous with respect to the weak topology.
\end{proof}

Lower semicontinuity of $\langle \nu,\mathfrak R^{\ssup{T}}(\nu)\rangle$ is part of the proof that the rate function $I_\mu^{\ssup{T}}$ defined in \eqref{Imuchar} has compact sublevel sets under certain assumptions (since $H(\cdot|M^\ssup{T}_{b\mu})$ has the property); see Section~\ref{sec-minimizerproof}.

Now let us turn to an issue that will arise in the proof of the LDP when we write the sum on $R$ in the exponent on the right-hand side of \eqref{DistVcalPoiThm} in terms of $\frac 1N Y_N$. First, let us write

\begin{equation}\label{RintermsoffrakR}
\sum_{i\not= j} R^{\ssup {T,N}}(\Xi_i,\Xi_j)=N\langle \smfrac 1N Y_N, \mathfrak R^{\ssup{T}}(\smfrac 1N Y_N)\rangle -\frac 1N\sum_i R^{\ssup T}(\Xi_i,\Xi_i),
\end{equation}
where $Y_N=\sum_i \delta_{\Xi_i}$ is the Poisson process on $\Gamma_T^{\ssup 1}$ with intensity measure $N M_{\mu,N}^{\ssup T}$ (note that $R^{\ssup{T,N}}=\frac 1N R^{\ssup T}$ if $K$ is replaced by $\frac 1N K$, since $R^{\ssup T}$ is linear in $K$).

Now, for the proof of the upper bound in the LDP in Section~\ref{sec-LDPproof} we will need that the last term on the right-hand side of \eqref {RintermsoffrakR} is small. This is provided in the next lemma. We write $Y_{N,0}=(Y_N)_0$ and $(\Xi_i)_0=\Xi_i(0)$ for the projection of $Y_N$ respectively $\Xi_i$ on the time marginal at time 0. We fix a majorizing function $f\colon (0,\infty)\to (0,\infty)$ satisfying $f(r)/r\to\infty$ as $r\to\infty$ and $f(r)\geq r$ for any $r$.

\begin{lemma}[Upper bound for the diagonal term]\label{lem:diag-term}
Fix a point process $Y_N=\sum_i \delta_{\Xi_i}$ on $\Gamma_T^{\ssup 1}$ and assume the upper bound on $K$ from \eqref{AssK1}.  Then, for any  $L,N\in \N$,
\begin{equation}\label{upbd-diag}
\frac 1N\sum_i R(\Xi_i,\Xi_i) \leq HTL\Aold   + HT N \Aold^2\eps_L^2,\qquad \mbox{with }\eps_L=\sup_{r>L}\frac r{f(r)},
\end{equation} 
holds on the event $\{\smfrac 1N Y_N\in \Acal_{f,\Aold}\}$.
\end{lemma}

Hence, the diagonal term is $o(N)$. If a gel would be present, i.e., if $|\Xi_{i,0}|$ was of order $N$ for some $i$, then it would be $\asymp N$.

\begin{proof}
We use the  estimate \eqref{Resti} that holds under \eqref{AssK1} and the fact that $\|\xi_0\|_1 = |\xi_0|$ for $\xi\in \Gamma_T^\ssup{1}$ to get
\begin{equation}\label{diagonaltermsmall}
\begin{aligned}
\frac 1N\sum_i R^\ssup{T}(\Xi_i,\Xi_i) &= \int_{\Gamma^\ssup{1}_T} \frac 1N Y_N(\d \xi)\,  R^\ssup{T}(\xi, \xi)\leq H T \int \frac{1}{N}Y_N(\d \xi)\, |\xi_0|^2\\
&\leq HT  L\int_{{|\xi_0|\leq L}} \frac{1}{N}Y_N(\d \xi)\, |\xi_0|  +HT \int_{|\xi_0|> L}\frac{1}{N}Y_N(\d \xi)\, |\xi_0|^2 
\end{aligned}
\end{equation}
Now, if $\{\smfrac 1N Y_N\in \Acal_{f,\Aold}\}$ we have that the first integral is smaller than $\Aold$ and for the second integral we have that
\[
\int_{|\xi_0|> L}\frac{1}{N}Y_N(\d \xi)\, |\xi_0|^2 \leq N \Big(\int_{|\xi_0|> L}\frac{1}{N}Y_N(\d \xi)\, |\xi_0|\Big)^2 \leq N \Aold^2\eps_L^2.
\]

\end{proof}

\subsection{The measure $\Q_k^{\ssup T}$}\label{sec-Qproof}

In this section, we study the limit of $\Q^\ssup{T,N,N}_k := N^{|k|-1}\P_k^\ssup{N}(\Xi \in \cdot \,)|_{\Gamma^\ssup{1}_T}$ as $N\to\infty$ for $ k\in \Mcal_{\N_0}(\Scal)$, i.e., the restriction of the distribution of the coagulation process to the space $\Gamma^\ssup{1}_T$ of one-particle coagulation trajectories. Recall that the (additional) super-index \lq$N$\rq\, indicates that we have replaced the kernel $K$ by $\frac 1N K$. The limiting measure $\Q_k^\ssup{T}$ will play an important role in the description of the limiting behaviour of the coagulation process. We will identify $\Q^\ssup{T}_k$ in terms of an explicit formula, using a kind of chart, i.e., a push-forward measure under some explicit measure on the space $[(\Scal\times \N)^2\times \Scal\times (0,\infty)]^{|k|-1}$, that carries all the data needed to understand which transition is happening in each of the $|k|-1$ coagulation steps.

Fix an initial configuration $k\in \Mcal_{\N_0}(\Scal)$ and a coagulation kernel $K$ as in Section~\ref{sec-model}. Let us identify the distribution of the coagulation process $\Xi=(\Xi_t)_{t\in[0,T]}$ under $\P_k$ on ${\Gamma^\ssup{1}_T}$, more precisely, on the set
\begin{equation}\label{Gamma1k}
\Gamma_{T,k}^{\ssup 1}=\big\{\xi\in\Gamma^\ssup{1}_{T}\colon \xi_0=k\big\}.
\end{equation}
Recall that coagulation trajectories $\xi\in \Gamma^\ssup{1}_{T,k}$ are only allowed to perform steps as in \eqref{steps}. For a pair of particles at $(x,m),(x^\prime, m^\prime)\in \Scal\times \N$ that coagulates into a particle with location $z\in \Scal$, this step is given by the addition of the signed measure
\begin{equation}\label{Wdef}
W^{\ssup z}_{(x,m),(x',m')}=-\delta_{(x, m)}-\delta_{(x', m')} +\delta_{(z,m+m')}.
\end{equation}
For a fixed tuple $(y_i,y_i',z_i)_{i=1, \ldots, |k|-1}$ with $y_i,y_i'\in \Scal\times \N$ and $z_i\in\Scal$,  we define
\begin{equation}\label{defphi}
\phi_0(\cdot,m)=k(\cdot)\delta_1(m)\mbox{ and }\phi_i=\phi_0+\sum_{j=1}^i W^{\ssup{z_j}}_{y_j,y_j'}, \text{ for } i=1, \ldots, |k|-1,
\end{equation}
and say that $((y_i,y_i',z_i)_{i=1, \ldots, |k|-1})$ is compatible with $k$ if
\begin{equation}\label{compatibility}
\phi_i\geq 0,\quad \mbox{ for }i=0,\dots,|k|-1.
\end{equation}
Let ${\mathfrak X}_k\subset [(\Scal\times\N)^2 \times \Scal]^{|k|-1}$ denote the set of such tuples. Furthermore, introduce the set of admissable time tuples,
$$
{\mathfrak F}_k=\Big\{(s_1,\dots,s_{|k|-1})\in [0,\infty)^{|k|-1}\colon\sum_{i=1}^{|k|-1}s_i\leq T\Big\}
$$ 
and define 
$$
\Psi_k\colon {\mathfrak X}_k \times {\mathfrak F}_k\to \Gamma_{T,k}^{\ssup 1},\qquad 
\big(y_i,y_i',z_i, s_i\big)_{i=1, \ldots, |k|-1} \mapsto \xi=(\xi_t)_{t\in[0,T]},
$$
by 
\begin{equation}\label{xitidenit}
\xi_t=\sum_{i=1}^{|k|} \1\{t\in I_i\}\phi_{i-1},\qquad\mbox{where }I_i=
\begin{cases}
\Big[\sum_{j=1}^{i-1}s_j, \sum_{j=1}^{i}s_j\Big),&\mbox{for } i<|k|,\\
&\\
\Big[\sum_{j=1}^{|k|-1}s_j,T\Big],&\mbox{for }i=|k|,
\end{cases}
\end{equation}
where $\phi= (\phi_i)_{i=1,\ldots, |k|-1}$ is defined in \eqref{defphi}.  
In words, if, starting from $\phi_0$, for $i=1,\dots,|k|-1$, iteratively after a time elapsure of $s_i$ time units, two particles at $y_i=(x_i,m_i)$ and $y_i'=(x_i',m_i')$ coagulate into a particle at $(z_i,m_i+m_i')$, then at each time $t\in[0,T]$, the configuration is equal to $\xi_t$. Note that $\xi_T$ is a delta-measure, i.e., after time $\sum_{j=1}^{|k|-1} s_j$ there is no coagulation possible anymore.

It is not hard to see that the mapping $\Psi_k$ is a bijection. We will leave the details to the reader. 

Now we describe the distribution of $\phi$. Recall the definition of $K_\phi$ and ${\bf K}_\phi$ from around \eqref{Mkerneldef} and recall that we conceive ${\bf K}_\phi$ as a measure on $(\Scal\times \N)^2\times\Scal$ and $K_\phi$ as its marginal on $(\Scal\times \N)^2$; in the first four arguments $K_\phi$ and ${\bf K}_\phi$ are indeed point measures. With a slight abuse of notation, one can conceive ${\bf K}_{\phi_{i-1}}(\d (x,m), \d (x',m'),\d z)$  as a Markov kernel, since $ \varphi_i=\varphi_{i-1}+W^{\ssup{z_i}}_{y_i,y_i'}$.  Their product $\bigotimes_{i=1}^{|k|-1} {\bf K}_{\phi_{i-1}}$ is concentrated on $\mathfrak X_k$. 

\begin{lemma}[The distribution of a one-particle coagulation tree]\label{xidecomp}
Fix $k\in \Mcal_{\N_0}(\Scal)$, $k\neq 0$. Then, we have, for any measurable bounded test function $f\colon \Gamma_{T,k}^{\ssup 1}\to\R$,
\begin{equation}\label{Xidistribution}
\begin{aligned}
\E_k\big(f(\Xi)\1\{\Xi \in\Gamma_{T,k}^{\ssup 1}\}\big)
&=\int_{{\mathfrak X}_k\times {\mathfrak F}_k}\Big(\bigotimes_{i=1}^{|k|-1} {\bf K}_{\phi_{i-1}}\otimes\bigotimes_{i=1}^{|k|-1} \d s_i\Big)(\d  \Theta)
\,f(\Psi_k(\Theta))\,\e^{- \widetilde\varphi_k(\Theta)},
\end{aligned}
\end{equation}
where
\begin{equation}\label{varphitildedef}
\widetilde \varphi_k(\Theta)=\frac 12\sum_{i=1}^{|k|-1} s_i\Big[\langle \phi_{i-1}, K \phi_{i-1}\rangle - \langle \phi_{i-1}, K^{\ssup{\rm diag}}\rangle\Big],\qquad \Theta=(y_i,y_i',z_i, s_i)_{i=1, \ldots, |k|-1},
\end{equation}
where we introduced $K^{\ssup{\rm diag}}(y)=K(y,y)$ for $y\in\Scal\times \N$.
\end{lemma}

The right-hand side of \eqref{Xidistribution} is the image measure under $\Psi_k$ of the restriction of the measure  $\bigotimes_{i=1}^{|k|-1} ({\bf K}_{\phi_{i-1}}\otimes \d s_i)$ to its support $\mathfrak X_k\times \mathfrak F_k$ with density $\e^{- \widetilde \varphi_k}$.

\begin{proof}
Fix $i\in\{1,\dots,|k|-1\}$.  During the time interval $I_i$, for any unordered pair $\{y,y'\}$ of elements of $\Scal\times \N$, there are $\phi_{i-1}(y) \phi_{i-1}(y')$ if $y\not=y'$, respectively $\frac 12\phi_{i-1}(y)(\phi_{i-1}(y)-1)$ if $y=y'$ independent exponential holding times running with parameter $K(y,y')$, one of which elapses at the respective coagulation time, namely one with unordered pair $\{y_i,y_i\}$. There are only finitely many such exponential clocks  involved, since $\phi_{i-1}$ has a finite support. The Lebesgue density for the event that the unordered pair of particles $\{y_i, y_i'\}$ coagulates at time $s_i$ as the first pair in the configuration is
\begin{equation}\label{exponentials-ith-coag}
\begin{aligned}
s_i&\mapsto K_{\phi_{i-1}}(y_i, y_i')\e^{-s_i K_{\phi_{i-1}}(y_i, y_i')}\Big(\prod_{\{y,y'\}\not=\{y_i, y_i'\}\colon y\not= y'}\e^{-s_i\phi_{i-1}(y)K(y,y')\phi_{i-1}(y')} \Big)\\
&\qquad\times\prod_{y\colon \{y\}\not=\{y_i, y_i'\}}\e^{-s_i\frac 12 \phi_{i-1}(y)K(y,y)(\phi_{i-1}(y)-1)}  \\
&= K_{\phi_{i-1}}(y_i, y_i')\e^{-s_i \eta^{\ssup i}(y_i,y_i')},
\end{aligned}
\end{equation}
where, for $y_i\not= y_i'$,
$$
\begin{aligned}
&\eta^{\ssup i}(y_i,y_i')= \frac 12 \sum_{y,y'}\phi_ {i-1}(y)K(y,y')\phi_{i-1}(y')-\frac 12 \sum_{y}\phi_{i-1}(y) K(y,y)\\
&\quad =\sum_{\{y,y'\}\colon y\not= y'}\phi_ {i-1}(y)K(y,y')\phi_{i-1}(y')+\sum_{y\colon \{y\}\not=\{y_i,y_i'\}} \frac 12 \phi_{i-1}(y)K(y,y)(\phi_{i-1}(y)-1),
\end{aligned}
$$
and for $y_i=y_i'$
$$
\eta^{\ssup i}(y_i,y_i')=\sum_{\{y,y'\}\colon y\not= y'}\phi_ {i-1}(y)K(y,y')\phi_{i-1}(y')+\sum_{y} \frac 12 \phi_{i-1}(y)K(y,y)(\phi_{i-1}(y)-1),
$$
i.e., the same formula. 

The probability that the new particle is  placed at $z_i$ is expressed by multiplying \eqref{exponentials-ith-coag} with  $\Upsilon((x_i, m_i),(x'_i,m'_i),\d z_i)$, which turns the first factor into ${\mathbf K}_{\phi_{i-1}}(y_i,y_i',\d z_i)$.

Because of the Markov property of the coagulation process, the probability $\P_k(\Xi\in\d\xi)$ is equal to the product over $i=1, \ldots, |k|-1$ of \eqref{exponentials-ith-coag}. Noting that $\widetilde\varphi_k(\Theta)= \sum_{i=1}^{|k|} \eta^{\ssup i}(y_i,y_i')$, this implies \eqref{Xidistribution}.
\end{proof}

Now we introduce the measure 
\begin{equation}\label{QTNdef}
\Q^\ssup{T,N}_k(\cdot)=N^{|k |-1} \P_k(\Xi\in\cdot)|_{\Gamma_T^{\ssup 1}}\in\Mcal(\Gamma_T^{\ssup 1}).
\end{equation}
The next question that we consider is the identification of the limit of $\Q^\ssup{T,N,N}_k$ as $N\to\infty$ for fixed $k\in \Mcal_{\N_0}(\Scal)$, where we recall that we add a superindex \lq$N$\rq\ when we replace the kernel $K$ by $\frac 1N K$. Using the description of $\xi$ that was given in Lemma \ref{xidecomp} we define a measure $\Q_k^{\ssup{T}}$ on $\Gamma_{T.k}^{\ssup 1}$ by dropping just the density in \eqref{Xidistribution}:
\begin{equation}\label{Qdef}
\Q_k^{\ssup{T}}(\d \xi)=\Big(\bigotimes_{i=1}^{|k|-1} {\bf K}_{\phi_{i-1}}\otimes \bigotimes_{i=1}^{|k|-1}\d s_i\Big)\circ \Psi_k^{-1}(\d \xi)=\e^{\varphi_k(\xi)}\, \P_k\big(\xi\in \Gamma_{T,k}^{\ssup 1};\Xi\in\d\xi\big),
\end{equation}
where
\begin{equation}\label{varphidef}
\varphi_{k}(\xi)=\frac 12\int_0^T \Big[\langle \xi_t, K\xi_t\rangle-\langle \xi_t, K^{\ssup{\rm diag}}\rangle\Big] \,\d t,\qquad k=\xi_0.
\end{equation}
Note that $\varphi_k(\Psi_k(\Theta))=\widetilde \varphi_k(\Theta)$ for $\Theta\in {\mathfrak X}_k$. 

\begin{lemma}[Limit of $\P_k^{\ssup{N}}(\Xi\in\cdot)|_{\Gamma_T^{\ssup 1}}$]\label{lem-limitconnection}
Fix any $k\in \Mcal_{\N_0}(\Scal)$, $k\neq 0$. Then 
\begin{equation}\label{Qdensity}
\frac{\d \Q_k^{\ssup{T,N,N}}}{\d \Q_k^{\ssup{T}}}(\xi)=\e^{-\frac 1N\varphi_k(\xi)},\qquad \xi\in\Gamma_{T,k}^{\ssup 1}.
\end{equation} 
If $K$ satisfies \eqref{AssK1}, then the density satisfies
\begin{equation}\label{density-bounds}
0\leq \varphi_{k}(\xi)\leq \frac 12 HT |k|^2, \qquad \xi\in\Gamma_{T,k}^{\ssup 1},
\end{equation}
and in particular,
\begin{equation}\label{eq:QNconvergence}
\Q_k^{\ssup{T}}=\lim_{N\to\infty}\Q^\ssup{T,N,N}_k,
\end{equation}
weakly and in total variation.
\end{lemma}

\begin{proof} We apply formula \eqref{Xidistribution} to the coagulation model with kernel $\frac 1N K$ (instead of $K$). Hence, the $(|k|-1)$-fold product of the $K$-terms receives a prefactor $N^{-(|k|-1)}$, which is compensated by the prefactor $N^{|k|-1}$ in \eqref{QTNdef}. Hence, on $\Gamma^\ssup{1}_{T,k}$ we have that 
\begin{align*}
\Q^\ssup{T,N,N}(\d\xi) &= \Big(\bigotimes_{i=1}^{|k|-1} {\bf K}_{\phi_{i-1}}\otimes \bigotimes_{i=1}^{|k|-1} \d s_i\Big)\circ \Psi_k^{-1}(\d\xi)\,\exp\Big(-\frac1{N}\varphi_k(\Theta^{-1}(\xi))  \Big)\\
&= \exp\Big(-\frac{1}{N}\varphi_k(\xi)\Big)\,\Q^\ssup{T}_k(\d\xi).
\end{align*}
This shows \eqref{Qdensity}.

To show \eqref{density-bounds} we use that $\langle \phi, K \phi\rangle \leq H \|\phi\|_1^2$ for any non-trivial $\phi\in \Mcal_{\N_0}(\Scal\times \N)$ and that  $\|\xi_t\|_1 = \|\xi_0\|_1 = |k|$ for all $t\in [0,T]$ and for any $\xi\in \Gamma^\ssup{1}_{T,k}$. Therefore, we get 
\begin{equation}\label{density-bd-details}
\varphi_{k}(\xi)\leq \frac{1}{2}\int_0^T \langle \xi_t, K\xi_t\rangle \,\d t \leq \frac{1}{2}HT|k|^2.
\end{equation}

The bound \eqref{density-bounds} obviously implies convergence in total variation as well as weak convergence.
\end{proof}

In particular, we can identify the true $N$-dependence of the reference measure $M_{b\mu,N}^{\ssup T}$ defined in \eqref{Mmudef}, when replacing the kernel $K$ by $\frac 1N K$: 

\begin{cor} Fix $\mu\in\Mcal_1(\Scal)$ and $b\in(0,\infty)$ and $N\in \N$ and replace $K$ by $\frac 1N K$ (adding a superscript $^{\ssup N}$). Then
\begin{equation}\label{MNident}
M_{b\mu,N}^{\ssup{T,N}}(\d \xi)=M_{b\mu}^{\ssup{T}}(\d \xi)\,\e^{-\frac 1N \varphi_{\xi_0}(\xi)}.
\end{equation}
\end{cor}

Now we can rewrite the representation in Theorem~\ref{thm-LDP} for $K$ replaced by $\frac 1N K$ in such a way that the intensity measure of the reference PPP does not depend on $N$ (up to the prefactor $N$): we just carry out a change of measure from $N M_{b\mu,N}^{\ssup{T,N}}$ to $N M_{b\mu}^{\ssup{T}}$ in the Poissonian expectation.

\begin{cor}\label{DistVcal-rescaledkernel}
Replace the kernel $K$ by $\frac{1}{N} K$ (which is denoted via the superscript $^\ssup{N}$). Then, for any $b\in(0,\infty)$ such that $|M^\ssup{T}_{b\mu}|<\infty$ and for any bounded and continuous test function $f\colon \Mcal(\Gamma_T^{\ssup 1})\to\R$,
\begin{equation*}
\begin{aligned}
\E^\ssup{N}_{\Poi_{N\mu}}\big(f(\Vcal_N^{\ssup T})\big)&=
 {\tt E}_{N M^\ssup{T}_{b\mu}} \Big[\left({\e^{b-1}}\right)^{|Y_N|}b^{-\sum_{i}|\Xi_{i,0}|}\exp\Big\{-\frac {1}{2N}\sum_{i,j \colon i\not=j}R^{\ssup T}(\Xi_i,\Xi_j)\Big\} \\
 &\quad\, f\big(\smfrac{1}{N}Y_N\big)\e^{-\frac 1N \int\varphi_{\xi_0}(\xi)\,Y_N(\d\xi)}\Big]\, \e^{N(|M^\ssup{T}_{b\mu}|-1)}.
 \end{aligned}
\end{equation*}
\end{cor}

Note that the assumption $|M^\ssup{T}_{b\mu}|<\infty$ holds for all sufficiently small $b$ under the assumption \eqref{AssK1}, according to Lemma \ref{lem-refmeasure-finite}.

\begin{lemma}\label{lem-tauesti}
We define $\sigma\colon \Mcal_{\N_0}(\Scal\times \N)\to [0,\infty)$ recursively via $\sigma(\delta_{(x,m)}) = 1$ for any  $ (x,m)\in \Scal\times \N$ and
\begin{equation*}
\sigma(\phi) = \sum_{(x,m),(x^\prime, m^\prime)} \int_\Scal  {\bf K}_{\phi}\big( (x,m),(x^\prime, m^\prime),\d z\big) \,\sigma\big( \phi - \delta_{(x,m)} -\delta_{(x^\prime, m^\prime)} + \delta_{(z, m+ m^\prime)}\big),
\end{equation*}
where the sum is taken over the support of $\phi$.
Then, the measure $\Q^\ssup{T}_k$ has total mass equal to
\begin{equation}\label{eq:expression_tau}
\Q^{\ssup{T}}_k(\Gamma_{T,k}^{\ssup 1})=\frac{T^{|k|-1}}{(|k|-1)!}\sigma(\phi_0),\qquad\mbox{where }\phi_0(\d x,m) = k(\d x) \delta_1(m).
\end{equation}
\end{lemma}
\begin{proof}
Recall definition \eqref{Qdef} and the fact that $\Psi_k \colon {\mathfrak X}_k \times {\mathfrak F}_k\to \Gamma_{T,k}^{\ssup 1}$ is a bijection. Then 
\[
\Q^{\ssup{T}}_k(\Gamma_{T,k}^{\ssup 1}) = \Big( \bigotimes_{i=1}^{|k|-1} {\bf K}_{\phi_{i-1}}\Big)({\mathfrak X}_k) \Big(\bigotimes_{i=1}^{|k|-1} \d s_i \Big)({\mathfrak F}_k).
\]
Note that
\[
\Big(\bigotimes_{i=1}^{|k|-1} \d s_i \Big)({\mathfrak F}_k) = \int_{[0,T]^{|k|-1}}(\d (s_1, \ldots, s_{|k|-1}))\,\1\Big\{\sum_{i=1}^{|k|-1}s_i\leq T\Big\}=\frac {T^{|k|-1}}{(|k|-1)!}.
\]
We now generalise the definition of ${\mathfrak X}_k$ in order to write down the recursion. We fix a point measure $\phi \in \Mcal_{\N_0}(\Scal \times \N)$ and define for any tuple $((y_i,y'_i,z_i)_{i=1, \ldots, |\phi|-1}$ the finite sequence
\[
 \phi_i=\phi+\sum_{j=1}^i W^{\ssup{z_j}}_{y_j,y_j'}, \text{ for } i=1, \ldots, |\phi|-1,
\]
where we recall the definition of $W^\ssup{z}_{(y,y')}$ given in \eqref{Wdef}. We denote by ${\mathfrak X}_\phi$ the set of all tuples $((y_i,y'_i,z_i)_{i=1, \ldots, |\phi|-1}$ that are compatible, i.e. that $\phi_i \geq 0$ for all $i=1, \ldots, |\phi|-1$. Observe that 
\[
{\mathfrak X}_\phi = \bigcup_{(y,y',z)\in \supp(\phi)^2 \times \Scal} \{(y,y',z)\}\times {\mathfrak X}_{\phi - \delta_{y} -\delta_{y'} + \delta_{(z, m(y)+ m(y^\prime))}}
\]
where for any $y= (x,m)\in \Scal \times \N$ we abbreviated $m(y)= m$. This implies that
\begin{align*}
\sigma(\phi) &:= \Big( \bigotimes_{i=1}^{|\phi|-1} {\bf K}_{\phi_{i-1}}\Big)({\mathfrak X}_\phi) = \sum_{(x,m),(x^\prime, m^\prime)} \int_\Scal  {\bf K}_{\phi}\big( (x,m),(x^\prime, m^\prime),\d z\big) \,\\
&\hspace{4cm}\Big( \bigotimes_{i=2}^{|\phi|-1} {\bf K}_{\phi_{i-1}}\Big)({\mathfrak X}_{\phi - \delta_{(x,m)} -\delta_{(x^\prime, m^\prime)} + \delta_{(z, m+ m^\prime)}}),
\end{align*}
which is the claimed recursion and implies the result.

\end{proof}

\begin{lemma}[Bounds on the total mass of $\Q^\ssup{T}$]\label{lem-taubound}
Assume that the kernel $K$ satisfies \eqref{AssK1} with constant $H$. Then, for any $k\in \Mcal_{\N_0}(\Scal)$, $k\neq 0$, we have that 
\begin{equation}\label{tauestiprecise}
\Q^{\ssup{T}}_k(\Gamma_{T,k}^{\ssup 1})\leq \frac{(TH)^{|k|-1}}{(|k|-1)!} |k|^{2(|k|-1)}.
\end{equation}
\end{lemma}

\begin{proof}
We use Lemma \ref{lem-tauesti} and show via induction over $|\phi|$ that 
\begin{equation}\label{sigma-esti}
\sigma(\phi) \leq H^{|\phi|-1} \|\phi\|_1^{2(|\phi|-1)}.
\end{equation}
 For $|\phi| = 1$,  both sides of \eqref{sigma-esti} are equal to $1$. Fix any $\phi\in \Mcal(\Scal\times \N)$ with $|\phi|\geq 2$. 
  We use the recursion from Lemma \ref{lem-tauesti}, the induction hypothesis and the fact that for any $\widetilde \phi = \phi - \delta_{(x,m)} -\delta_{(x^\prime, m^\prime)} + \delta_{(z, m + m')}$ we have that $|\widetilde\phi| = |\phi|-1$ and $\|\widetilde \phi\|_1 = \| \phi\|_1$. Also recall from \eqref{AssK1} that the definition of $H$ implies that $\langle \phi, K \phi\rangle \leq H \|\phi\|_1^2$. Then
\begin{align*}
\sigma(\phi) &\leq H^{|\phi|-2}\| \phi\|_1^{2(|\phi|-2)} \sum_{(x,m),(x^\prime, m^\prime)}  K_{\phi}\big( (x,m),(x^\prime, m^\prime)\big) 
 \\
& \leq H^{|\phi|-2} \|\phi\|_1^{2(|\phi|-2)} \langle \phi, K \phi\rangle \leq H^{|\phi|-1} \|\phi\|_1^{2(|\phi|-1)}
\end{align*}
This implies \eqref{tauestiprecise}.


\end{proof}

\begin{rem}[Examples of explicit expressions for $\Q^{\ssup{T}}_k(\Gamma_{T,k}^{\ssup 1})$]
In some specific cases the recursive definition of $\Q^{\ssup{T}}_k(\Gamma_{T,k}^{\ssup 1})$ in \eqref{eq:expression_tau} can be solved and the quantity can be expressed explicitly for every $k\in \Mcal_{\N_0}(\Scal)\setminus\{0\}$. Notable cases are the non-spatial kernels of multiplicative and of additive type. When $K((x,m), (x',m'))=mm'$, then $\Q^{\ssup{T}}_k(\Gamma_{T,k}^{\ssup 1})=\frac{T^{|k|-1}}{(|k|-1)!}|k|^{2(|k|-1)}$. When $K((x,m), (x',m'))=m+m'$, then $\Q^{\ssup{T}}_k(\Gamma_{T,k}^{\ssup 1})=\frac{T^{|k|-1}}{2^{|k|-1}}|k|!$.\hfill$\Diamond$
\end{rem}

\subsection{Proof of Lemma~\ref{lem-Contrho}}\label{sec-proofrhocont}

In this section, we prove Lemma~\ref{lem-Contrho}, i.e., the continuity of the map $\nu\mapsto \rho(\nu)$ defined in \eqref{rhodef}.
First we show the continuity of every marginal:

\begin{lemma}[Continuity of $\nu\mapsto \rho_t(\nu)$]\label{lem-Xicontinuous} Fix any  $\Aold>0$ and some function $f \colon \N\to [0,\infty)$ that grows at infinity faster than linear, i.e., $f(r)/r\to\infty$ as $r\to\infty$. Let $(\nu_n)_{n\in\N}$ be a sequence in $\Acal_{f,\Aold}$ that converges towards some $\nu$ that has a density with respect to $M_{b\mu}^{\ssup T}$ for some $b>0$. Then, for any $t\in[0,T]$, $\rho_t(\nu_n)\to\rho_t(\nu)$ as $n\to\infty$.
\end{lemma}
\begin{proof}
Fix $A\subset \Scal\times\N$ with $\rho_t(\nu)(\partial A)=0$. It suffices to show that $\rho_t(\nu_n)(A)\to \rho_t(\nu)(A)$ as $n\to\infty$. Note that, for any $L\in\N$,
$$
\rho_t(\nu_n)(A)= \int \nu_n(\d \xi)\,\xi_t(A)\1\{|\xi_0|\leq L\}
+\int \nu_n(\d \xi)\,\xi_t(A)\1\{|\xi_0|> L\}.
$$
The last term vanishes uniformly in $A$ and $n$ as $L\to\infty$, since
\begin{equation}\label{Lcutting}
\begin{aligned}
\int \nu_n(\d \xi)\,\xi_t(A)\1\{|\xi_0|> L\}
&\leq \int \nu_n(\d \xi)\,|\xi_t|\1\{|\xi_0|> L\}=\int \nu_n(\d \xi)\,|\xi_0|\1\{|\xi_0|> L\} \\
&\leq \eps_L \int \nu_n(\d \xi)\,f(|\xi_0|)\leq \eps_L \beta,
\end{aligned}
\end{equation}
where $\eps_L=\inf_{r>L}r/f(r)$ vanishes as $L\to\infty$.

Concerning the first term, we now show that the map $\xi\mapsto \xi_t(A)\1\{|\xi_0|\leq L\}$ is continuous in each $\xi$ that satisfies $\xi_t(\partial A)=0$ and does not jump in $t$. Let $ \xi$ be such a point, and pick a sequence $(\xi^{\ssup n})_{n\in\N }$ that converges to $\xi$. Then $|\xi_0^{\ssup n}|\to |\xi_0|$ as $n\to\infty$. If $|\xi_0|>L$, then we have $\lim_{n\to\infty}\xi^{\ssup n}_t(A)\1\{|\xi^{\ssup n}_0|\leq L\}\to0= \xi_t(A)\1\{|\xi_0|\leq L\}$. Otherwise, for any sufficiently large $n$ (recall that $|\xi_s|\in\N$ for any $s$)  we have $\xi^{\ssup n}_t(A)\1\{|\xi^{\ssup n}_0|\leq L\}=\xi^{\ssup n}_t(A)\to \xi_t(A)=\xi_t(A)\1\{|\xi_0|\leq L\}$ because $\xi^{\ssup n}_t\to \xi_t$ weakly (since $\xi$ is continuous in $t$) since $\xi_t(\partial A)=0$.

Finally we need to show that the set of considered $\xi$ exhausts all $ \xi$, i.e., that $\nu(\{\xi\colon \xi_t(\partial A)>0\})=0$ and $\nu(\{\xi\colon \xi\mbox{  jumps at }t\})=0$. The first holds since $\xi_t(\partial A)$ is $\N_0$-valued and hence $\nu(\{\xi\colon \xi_t(\partial A)>0\})\leq \int \nu(\d \xi)\, \xi_t(\partial A)=\rho_t(\nu)(\partial A)=0$. The second holds since $\nu$ has a density with respect to $M_{b\mu}^{\ssup T}$, and the latter has a density with respect to the distribution of the Marcus--Lushnikov process, which does not jump with positive probability at time $t$.
\end{proof}

Now we prove the continuity of the map $\nu\mapsto \rho(\nu)$ defined in \eqref{rhodef}.

\begin{proof}[Proof of Lemma~\ref{lem-Contrho}.]

We are going to show that $\lim_{n\to\infty}\sup_{t\in[0,T]}\mathfrak{d}(\rho_t(\nu_n),\rho_t(\nu))=0$, where $\mathfrak{d}$ is the L{\'e}vy-Prohorov metric on $\Mcal(\Scal\times \N)$ defined in \eqref{Prorohovmetric}. This implies convergence of $\rho(\nu_n)$ towards $\rho(\nu)$ with respect to the $J_1$-topology on $\D_T(\Mcal(\Scal\times \N))$. With a small parameter $\eps>0$, we decompose $[0,T]$ into pieces $I_i^{\ssup\eps}=[t_{i-1},t_i]$ of length $\leq \eps$ and use the triangle inequality to estimate
\begin{equation}\label{triangle}
\begin{aligned}
\sup_{t\in[0,T]}\mathfrak{d}(\rho_t(\nu_n),\rho_t(\nu))\leq \max_{i}\big[\mathfrak{d}&(\rho_{t_i}(\nu_n),\rho_{t_i}(\nu))\\&+\sup_{t\in I_i^{\ssup\eps}}\mathfrak{d}(\rho_t(\nu_n),\rho_{t_i}(\nu_n))+\sup_{t\in I_i^{\ssup\eps}}\mathfrak{d}(\rho_t(\nu),\rho_{t_i}(\nu))\big].
\end{aligned}
\end{equation}
The first of the three terms on the right vanishes as $n\to\infty$, according to Lemma~\ref{lem-Xicontinuous}. We estimate now the second. More generally, we give an upper bound for
\begin{equation}
\sup_{s,t\in[0,T]\colon |s-t|\leq \eps}\mathfrak{d}(\rho_t(\nu_n), \rho_s(\nu_n)).
\end{equation}
Let $s,t\in [0,T]$ be such that $s<t$ and $|s-t|\leq \eps$. According to the definition of the L{\'e}vy-Prohorov metric, if we can find some $\eta>0$ such that 
\begin{align*}
&\rho_t(\nu_n)(A)  \leq \rho_s(\nu_n)(A^\eta) +\eta\quad \text{ and }\quad 
  \rho_s(\nu_n)(A)  \leq \rho_t(\nu_n)(A^\eta) +\eta 
\end{align*}
holds for all measurable $A\subset \Scal\times \N$, 
then $\mathfrak{d}(\rho_t(\nu_n),\rho_s(\nu_n)) \leq \eta$. For any $\eta>0$ we have the estimate 
\[
\rho_t(\nu_n)(A) \leq \rho_t(\nu_n)(A^\eta) \leq \rho_s(\nu_n)(A^\eta) + \int \nu_n(\d \xi) \big| \xi_t(A^\eta) - \xi_s(A^\eta)\big|
\]
and the same for $s$ and $t$ exchanged. For the latter term we can estimate
\[
\int \nu_n(\d \xi) \big|\xi_t(A^\eta) - \xi_s(A^\eta)\big| \leq 2 \int \nu_n(\d \xi) J_{[s,t]}(\xi),
\]
where $J_I(\xi)$ is the number of jumps of $\xi$ in the interval $I\subset [0,T]$. This is true since all steps of $\xi$ are of the form of adding $\delta_{(z,m+m')}-\delta_{(x,m)}-\delta_{(y,m')}$ for some $x,y,z\in \Scal$ and $m,m'\in\N$.

Now, for any interval $I\subset [0,T]$ and $L\in\N$, we define 
\begin{equation}\label{def-J_function}
\Jcal_{I}^\ssup{L}(\nu)=\int\nu(\d \xi)\,J_{I}(\xi)\1\{|\xi_0|\leq L\},\qquad \nu\in \Mcal(\Gamma^\ssup{1}_T).
\end{equation}
Splitting the integrals into the cases $|\xi_0|\leq L$ and $|\xi_0|> L$ and proceeding as in \eqref{Lcutting}, we get  for any $\eps,\eta>0$,
$$
\begin{aligned}
\sup_{s,t\in[0,T]\colon |s-t|\leq \eps}\mathfrak{d}(\rho_t(\nu_n),\rho_s(\nu_n))
&\leq 2 \sup_{s,t\in[0,T]\colon |s-t|\leq \eps}\Jcal^\ssup{L}_{[s,t]}(\nu_n)+2\beta \eps_L\\
&=2\, \Jcal^\ssup{L}_{[s_n,t_n]}(\nu_n)+2\beta \eps_L,
\end{aligned}
$$
where $(s_n,t_n)$ are picked as a maximising pair in $[0,T]^2$ with $|s_n-t_n|\leq \eps$.

Now, along subsequences, we may assume that $(s_n,t_n)\to (s,t)\in[0,T]^2$ such that $|s-t|\leq\eps$. For a given $\delta>0$ and all sufficiently large $n$ in this subsequence, we have $[s_n,t_n]\subset [s-\delta,t+\delta]\cap[0,T]$. Furthermore, observe that $\xi\mapsto J_{[s-\delta,t+\delta]\cap[0,T]}(\xi)$ is upper semi-continuous. Indeed, if $\xi_n\to\xi$ in $\Gamma^\ssup{1}_T$, each jump of $\xi_n$ converges to a jump of $\xi$ and since each element is right-continuous, we have that $J_{[s-\delta,t+\delta]\cap[0,T]}(\xi_n)\leq J_{[s-\delta,t+\delta]\cap[0,T]}(\xi)$.  Hence we have, by \cite[Thm.~D.12]{DZ10},
\begin{equation}\label{ProofLemma(i)}
\limsup_{n\to\infty}\Jcal^{\ssup L}_{[s_n,t_n]}(\nu_n)\leq \Jcal^{\ssup L}_{[s-\delta,t+\delta]\cap[0,T]}(\nu).
\end{equation}

We show now that the right-hand side is not larger than $C_L (\eps+2\delta)$ for some $C_L>0$ that does not depend on $s$ nor on $t$. For doing this, we note that
\begin{equation}\label{nuboundM}
M_{b\mu}^\ssup{T}(A)<\e^{-1}\qquad\Longrightarrow\qquad \nu(A)\leq \frac{H(\nu| M_{b\mu}^\ssup{T})+\e^{-1}}{-1-\log M_{b\mu}^\ssup{T}(A) },\qquad A\subset \Acal_{f,\Aold}.
\end{equation}
Indeed, note that $x\mapsto x\log x+1-x$ is nonnegative and convex in $(0,\infty)$ and therefore 
\begin{align*}
H(\nu| M_{b\mu}^\ssup{T})& \geq  \int \Big (\frac{\d \nu}{\d M_{b\mu}^\ssup{T}}\log \frac{\d \nu}{\d M_{b\mu}^\ssup{T}}+1-\frac{\d \nu}{\d M_{b\mu}^\ssup{T}}\Big)\mathbf{1}_{A}\, \d M_{b\mu}^\ssup{T}\\
&\geq M_{b\mu}^\ssup{T}(A) \Big (\frac{\nu(A)}{M_{b\mu}^\ssup{T}(A)}\log \frac{\nu(A)}{M_{b\mu}^\ssup{T}(A)}+1-\frac{\nu(A)}{ M_{b\mu}^\ssup{T}(A)}\Big)\\
&\geq \nu(A)\log\nu(A)+\nu(A)(-1-\log M_{b\mu}^\ssup{T}(A))\\
&\geq  -\e^{-1}+\nu(A)(-1-\log M_{b\mu}^\ssup{T}(A)),
\end{align*}
where the second inequality is obtained thanks to Jensen's inequality, the third because of $M_{b\mu}^\ssup{T}(A)\geq 0$ and the last one thanks to the fact that $x\log x\geq -\e^{-1}$. We see that, when $M_{b\mu}^\ssup{T}(A)<\e^{-1}$, the bracket on the last line is positive and we obtain \eqref{nuboundM}.

Now 
note that the right-hand side of \eqref{ProofLemma(i)} may be bounded as 
\begin{equation}\label{eq:bound_on_J}
\Jcal_{[s-\delta,t+\delta]\cap[0,T]}^{\ssup L}(\nu)\leq \sum_{j=1}^L \nu(\{\xi\colon J_{[s-\delta,t+\delta]\cap[0,T]}(\xi)\geq j,\, |\xi_0|\leq L\}).
\end{equation}
Before bounding the right-hand side,  we are going to argue now that 
\begin{equation} \label{eq:bound_on_jumps_M}
M_{b\mu}^\ssup{T}(\{\xi\colon J_{[s-\delta,t+\delta]\cap[0,T]}(\xi)\geq j,\, |\xi_0|\leq L\})\leq \widetilde C_L (\eps+2\delta)^j
\end{equation} 
for all $j\in \{1, \ldots,L\}$ and some $\widetilde C_L>0$ that does not depend on $s$ nor on $t$. 
Recall that $M_{b\mu}^\ssup{T} (\d \xi) = \e^{2-b}\, \Poi_\mu \otimes \Q^\ssup{T}(\d \xi) \, b^{|\xi_0|}$ 
and the formula for $\Q^\ssup{T}$ from \eqref{Qdef}. For $\xi\in \Gamma^\ssup{1}_T$ let $k=\xi_0$ and let $(y_i, y^\prime_i,z_i,s_i)_{i=1, \ldots, |k|-1} \in \mathfrak{X}_k\times \mathfrak{F}_k$ be such that $\xi = \Psi_{k}((y_i, y^\prime_i,z_i,s_i)_{i=1, \ldots, |k|-1})$. Abbreviate $I= [s-\delta, t+\delta]\cap [0,T]$ and note that the value of $J_I(\xi)$ only depends on the inter-coagulation times $(s_i)_{i=1, \ldots, |k|-1}\in \mathfrak{F}_k$, which are i.i.d.\ with distribution given by the Lebesgue measure (see \eqref{Qdef}). More precisely, we can recover the  jump times of $\xi$ via $t_i= s_1 + \ldots + s_i$, for $i=1, \ldots, |k|-1$, and have that $J_I(\xi) = \#\{i \colon t_i\in I\}$. Then,
\begin{align*}
\int_{\mathfrak{F}_k}& \bigotimes_{i=1}^{|k|-1} \d s_i \,\1\big\{\#\{i \colon t_i \in I\}\geq j\big\} \\
& =\int_{[0,T]^{|k|-1}} \bigotimes_{i=1}^{|k|-1} \d t_i \,\1\big\{\#\{i \colon t_i \in I\}\geq j\big\} \1\{t_1 < t_2 < \ldots < t_{|k|-1}\} \\
&= \frac{1}{(|k|-1)!}\int_{[0,T]^{|k|-1}} \bigotimes_{i=1}^{|k|-1} \d t_i  \,\1\big\{\#\{i \colon t_i \in I\}\geq j\big\} = \sum_{\ell\geq j}\frac{|I|^\ell (T-|I|)^{|k|-1-\ell}}{\ell! (|k|-1-\ell)!},
\end{align*}
where the last term is smaller than $|I|^j (2T)^{|k|-1}/(|k|-1)!$, if $|I|$ is small, which we can assume without loss of generality.
The terms of $\Q^\ssup{T}(\d \xi)$ that depend on $(y_i, y^\prime_i,z_i)_{i=1, \ldots, |k|-1}\in \mathfrak{X}_k$ can be estimated as in Lemma~\ref{lem-taubound}. Altogether we get that
\begin{equation}
M_{b\mu}^\ssup{T}(\{\xi\colon J_{[s-\delta,t+\delta]\cap[0,T]}(\xi)\geq j,\, |\xi_0|\leq L\}) \leq  (\eps+2\delta)^j \sum_{n=1}^{L} b^n (2HT)^{n-1}\frac{n^{2(n-1)}}{n!(n-1)!}.
\end{equation}

Via  \eqref{nuboundM} and the assumption that $H(\nu| M_{b\mu}^\ssup{T})<\infty$, this implies that also the right-hand side of \eqref{ProofLemma(i)} is not larger than $C_L (\eps+2\delta)$ for some $C_L>0$ that does not depend on $s$ nor on $t$.

Summarizing, we have shown that, for any $L\in\N$ and $\eps,\delta>0$,
$$
\limsup_{n\to\infty}\sup_{s,t\in[0,T]\colon |s-t|\leq \eps}\mathfrak{d}(\rho_t(\nu_n),\rho_s(\nu_n))\leq 2C_L(\eps+\delta)+2\beta\eps_L.
$$
We first pick $L$ large enough such that $\eps_L$ is small enough. Since the left-hand side does not depend on $\delta$, we may make $\delta\downarrow0$ on the right-hand side. In an analogous way, we derive the same bound for $\sup_{t,s\in[0,T]\colon |s-t|\leq \eps}  \mathfrak d(\rho_t(\nu),\rho_s(\nu))$. This implies via \eqref{triangle} our assertion.
\end{proof}

As a byproduct of the above proof, we have the following result. 
\begin{cor}[$\rho(\nu)$ is a (uniformly) continuous path] Fix any  $\Aold>0$ and some function $f \colon \N\to [0,\infty)$ that grows at infinity faster than linear, i.e., $f(r)/r\to\infty$ as $r\to\infty$. Then, if $\nu\in\Mcal(\Gamma_T^{\ssup 1})$ is such that $H(\nu| M_{b\mu}^\ssup{T})<\infty$, then $[0,T]\ni t\mapsto \rho_t(\nu)\in \Mcal(\Scal\times \N) $ is  uniformly  continuous.
\end{cor}

\section{Proof of Theorem~\ref{thm-LDP}: the LDP  for $\Vcal_{N}^{\ssup{T}}$}\label{sec-LDPproof}

\noindent In this section, we prove our second main result, the LDP for $\Vcal_{N}^{\ssup{T}}$ of Theorem~\ref{thm-LDP}. Since we will rely on the representation of Theorem~\ref{thm-distMi} in terms of the PPP $Y_N$ with intensity measure $NM_{b\mu}^{\ssup T}$, we first need an LDP for $Y_N$, which we provide in Section~\ref{sec-Yproof}. The proof of our LDP is carried through in Section~\ref{sec-LDPproof}.

\subsection{LDP for the PPP $Y_N$}\label{sec-Yproof}
We state and prove an LDP for the random variable $Y_N$ introduced in Theorem~\ref{thm-distMi}. We feel that this result is not new, but we did not find a reference, hence, we give an outline of  a proof. Recall the relative entropy from  \eqref{entropy}.

\begin{lemma}[LDP for $Y_N$ under PPP($N \mathfrak m$)]\label{lem-YLDP}
Assume that $\Xcal$ is a Polish space, and pick a finite and positive measure $\mathfrak m$ on $\Xcal$ and  assume that $Y_N$ is a Poisson point process with intensity measure $N\mathfrak m$. Then $(\frac 1N Y_N)_{N\in\N}$ satisfies an LDP on $\Mcal(\Xcal)$ with rate function $\nu\mapsto H(\nu\mid \mathfrak m)$. All the level sets $\{ \nu\colon H(\nu\mid\mathfrak m)\leq C\}$ with $C\in\R$ are compact.
\end{lemma}

\begin{proof} The abstract version of Cram\'er's theorem gives immediately an LDP for $\frac 1N Y_N$. Indeed, note that we have in distribution that $Y_N=Z_1+\dots+Z_N$, where $Z_1,\dots,Z_N$ are independent PPPs with intensity measure $\mathfrak m$. Then $\frac 1N Y_N$, as the average of i.i.d.~random objects, satisfies the LDP with  rate function equal to the Legendre transform of the logarithm of the moment generating function of $Z_1$, which reads
$$
\begin{aligned}
\Mcal^*(\Xcal)\ni \nu&\mapsto \sup_{f\in\Ccal_0(\Xcal)}\Big(\langle\nu, f\rangle-\log {\tt E}\big[\e^{\langle f,Z_1\rangle}\big]\Big)\\
&=\sup_{f\in\Ccal_0(\Xcal)}\Big(\langle\nu, f\rangle-\langle \e^{f}-1,\mathfrak m\rangle\Big),
\end{aligned}
$$
where we used a well-known formula for exponential Poisson moments, and  $\Ccal_0(\Xcal)$ is the closure in the uniform norm of the set of all continuous, compactly supported functions $f\colon \Xcal\to\R$, the dual of which is set $\Mcal_\pm(\Xcal)$ of all signed measures on $\Xcal$. 
Now that we have formulated our problem in a Banach space setting it is special case of the G\"artner--Ellis theorem. In order to finish the proof of the LDP that we need (i.e., to restrict from $\Mcal_\pm(\Xcal)$ to $\Mcal(\Xcal)$), we need to check that the restriction of the above rate function to $\Mcal_\pm(\Xcal)\setminus \Mcal(\Xcal)$ is constantly equal to $+\infty$, which we leave to the reader. 

Now we identify the rate function as $H(\nu|\mathfrak m)$ by standard means. Indeed, if $\nu\ll \mathfrak m$, then we may insert (continuous bounded approximations of) $f=\log\frac{\d\nu}{\d\mathfrak m}$ and obtain that the rate function is $\geq H(\nu|\mathfrak m)$, and the opposite inequality is seen by 
$$
H(\nu|\mathfrak m)=H(\nu|\e^f\,\d\mathfrak m)+\langle \nu, f\rangle-\langle \e^f,\mathfrak m\rangle +\mathfrak m(\mathcal X)\geq \langle \nu, f\rangle-\langle \e^f-1,\mathfrak m\rangle,
$$
since the entropy is nonnegative. In the case that $\nu$ is not absolutely continuous with respect to $\mathfrak m$, we take $f$ as a continuous and bounded approximation of $M\1_A$ with a large $M$ and a measurable set $A$ that satisfies $\mathfrak m(A)=0<\nu(A)$. See {\cite[Lemma 3.2.13]{DeSt89}} for details.

The sets $\{\frac{\d \nu}{\d \mathfrak{m}}\colon H(\nu\mid \mathfrak{m})\leq C \}$ are weakly compact in $L^1(\mathfrak{m})$ by uniform integrability since we can write $H(\nu\mid \mathfrak{m}) = \int \left(\frac{\d \nu}{\d \mathfrak{m}} \log \frac{\d \nu}{\d \mathfrak{m}} - \frac{\d \nu}{\d \mathfrak{m}} + 1  \right)\d \mathfrak{m}$. 
From this we see that the sets $\{\nu\colon H(\nu\mid \mathfrak{m})\leq C \}$ are compact in $\Mcal(\Xcal)$ with respect to its (functional analytic) weak topology and therefore also in the topology generated by testing against bounded measurable functions {\cite[Thrm~4.7.25]{Bog07}} , which certainly include all continuous bounded functions.

\end{proof}

\subsection{Proof of Theorem~\ref{thm-LDP}}\label{sec-finish}

Now let us derive the large-deviations principle for the distribution of $\Vcal_{N}^{\ssup{T}}$ with poissonised initial distribution under conditioning on $\Vcal_{N}^{\ssup{T}}\in\Acal_{f,\Aold}$ for any $\Aold\in(0,\infty)$, where we recall that we fixed a majorizing function $f$ such that $f(r)\geq r$ for any $r$ and $f(r)/r\to\infty$ as $r\to\infty$. Recall that we add an additional superindex $^{\ssup N}$ to indicate that $K$ is replaced by $\frac 1N K$. We recall that we are under the assumption on $K$ in \eqref{AssK1} and consider the distribution of $\Vcal_{N}^{\ssup{T}}$ under the poissonised initial measure, $\P^{\ssup N}_{\Poi_{N\mu}}$ with some $\mu\in\Mcal_1(\Scal)$.  We write $\Vcal_N$ instead of $\Vcal_N^{\ssup{T}}$. 

We pick $b\in(0,\infty)$ so small that $|M_{b\mu}^{\ssup T}|<\infty$, according to Lemma~\ref{lem-refmeasure-finite}. 

First we point out that the rate function in \eqref{ratefunctionA} has compact level sets. Indeed, as is seen from \eqref{Imuchar}, $I_\mu^{\ssup T}(\nu)$ is a sum of terms each of which is a lower semi-continuous function of $ \nu$ on $\Acal_{f,\Aold}$. Indeed, the quadratic term is lower semi-continuous by Lemma~\ref{Lem-RfrakProp}(2), the maps $\nu\mapsto \int \nu_0(\d k)\, |k|=|c_{\nu_0}|$ and $\nu\mapsto |\nu |$ are continuous on $\Acal_{f, \Aold}$ as is seen using the arguments from the proof of Lemma~\ref{Lem-RfrakProp}(3). Finally, by Lemma~\ref{lem-YLDP}, the entropy $H(\cdot|M_{b\mu}^{\ssup T})$ has even compact sublevel sets (in $\Mcal(\Gamma^\ssup{1}_T)$ and hence in $\Acal_{f, \Aold}$). 

Now we turn to the proofs of the upper and lower bounds. We start from Theorem~\ref{thm-distMi} with $K$ replaced by $\frac 1N K$, more precisely, from Corollary~\ref{DistVcal-rescaledkernel}. This gives, for any measurable set $E\subset \Gamma_T^{\ssup 1}$,
\begin{equation}\label{Prob(E)}
\P^{\ssup N}_{\Poi_{N\mu}}(\Vcal_N\in E)={\tt E}_{N M_{b\mu}^{\ssup T}}\Big[\e^{N\phi_{b}(\frac 1N Y_N)}\1\{\smfrac 1N Y_N\in E \}\e^{\frac 1{2N} \sum_i R^{\ssup T}(\Xi_i,\Xi_i)}\e^{- \int \varphi_{\xi_0}(\xi)\,\frac 1N Y_N(\d\xi)}\Big],
\end{equation}
where
\begin{equation}\label{phiNdef}
\phi_{b}(\nu) = - \frac{1}{2}\langle \nu,\mathfrak{R}^\ssup{T}(\nu) \rangle + |\nu|(b-1) -|c_{\nu_0}| \log b+|M_{b\mu}^{\ssup T}|-1.
\end{equation}
Observe that the last two terms in the expectation on the right-hand side are $\e^{o(N)}$ as $N\to\infty$, uniformly on $\{\frac 1N Y_N\in\Acal_ {f,\Aold}\}$ by  Lemmas~\ref{lem:diag-term} and \ref{lem-limitconnection}. 
Furthermore, $\phi_b$ is bounded and continuous on $\Acal_{f,\Aold}$ in the weak topology, according to the above remarks on the compactness of level sets of $I_\mu^{\ssup T}$. Futhermore, note that $\Acal_{f,\Aold}$ is closed, due to the continuity of $\xi\mapsto |\xi_0|$ and nonnegativity of $f$, using Fatou's lemma.

Now the LDP for $\frac 1N Y_N$ under ${\tt P}_{N M_{b\mu}^{\ssup T}}$ from Lemma~\ref{lem-YLDP}, together with Varadhan's lemma (Lemma 4.3.6 in \cite{DZ10}) implies, for any closed set $F\subset \Mcal(\Gamma_T^{\ssup 1})$ (implying that also $F\cap\Acal_{f,\Aold}$ is closed),
that
\begin{equation}\label{Upperbound}
\limsup_{N\to\infty}\frac 1N\log \P^{\ssup N}_{\Poi_{N\mu}}(\Vcal_N\in F\cap\Acal_{f,\Aold} )\leq -\inf\{H(\nu|M_{b\mu}^{\ssup T})-\phi_b(\nu)\colon \nu\in F\cap\Acal_{f,\Aold}\}.
\end{equation}
Observe from \eqref{Imuchar} that $H(\nu|M_{b\mu}^{\ssup T})-\phi_b(\nu)=I_\mu^{\ssup T}(\nu)$ for any $\nu\in\Mcal(\Gamma_T^{\ssup 1})$. In particular, recalling that $\chi_\Aold=\inf_{\nu\in\Acal_{f,\Aold}} I_\mu(\nu)$,
\begin{equation}\label{Upperboundchi}
\limsup_{N\to\infty}\frac 1N\log \P^{\ssup N}_{\Poi_{N\mu}}(\Vcal_N\in \Acal_{f,\Aold} )\leq- \chi_\Aold.
\end{equation}
In order to finish the proof of Theorem~\ref{thm-LDP}, we need to argue that also the complementary lower bound holds for  \eqref{Upperbound} for $F$ replaced by some open set $G\subset \Mcal(\Gamma_T^{\ssup 1})$. This will then imply the corresponding lower bound in \eqref{Upperboundchi}, which finishes the proof.

In this point, there is a technical problem, since the set $\Acal_{f,\Aold}$ is not open in $\Mcal(\Gamma_T^{\ssup 1})$. An obvious idea is to go to the set 
\begin{equation}\label{Acaldef-alt}
\Acal_{f,<\Aold}=\Big\{\nu\in \Mcal(\Gamma^\ssup{1}_T)\colon \int \nu(\d\xi)\,f(|\xi_0|)<\Aold\Big\}.
\end{equation}
However, $\Acal_{f,<\Aold}$ is still not an open set, since the map $\nu\mapsto  \int \nu(\d\xi)\,f(|\xi_0|)$ is not continuous (the map $\xi\mapsto f(|\xi_0|)$ not bounded). We solve this by applying some restriction argument. Indeed, for a large cutting  parameter $L\in\N$, we insert an indicator on the event that the PPP $Y_N$ has no particles that are larger than $L$, i.e., that it is concentrated on $\Gamma_{T,\leq L}^{\ssup 1}=\{\xi\in\Gamma_T^{\ssup 1}\colon \xi_t\in\Mcal_{\N_0}(\Scal\times [L])\,\forall t\in[0,T]\}$. Then we will condition on this event and make a change of measure from the PPP $Y_N$ with intensity measure $N M_{b\mu}^{\ssup T} $ to the intensity measure 
$NM_{b\mu}^{\ssup {T,\leq L}}$, the restriction of $N M_{b\mu}^{\ssup T} $ to $\Gamma_{T,\leq L}^{\ssup 1}$. 

In the following we often identify measures $\nu$ on $\Gamma_T^{\ssup 1}$ that satisfy $\nu(\Gamma^\ssup{1}_T \setminus \Gamma^\ssup{1}_{T, \leq L}) = 0$ with measures on   $\Gamma_{T,\leq L}^{\ssup 1}$.
We introduce the restriction operator $\Pi_L\colon\Mcal(\Gamma_T^{ \ssup 1})\to \Mcal(\Gamma_{T}^{\ssup 1})$, which maps $\nu$ to the measure $\Pi_L(\nu)= \nu^\ssup{\leq L}$, defined by $\nu^\ssup{\leq L}(\cdot) = \nu(\cdot\cap \Gamma_{T,\leq L}^{ \ssup 1})$. Note that the mapping $\Pi_L$ is continuous with respect to weak convergence. 
One easily sees, by definition of PPP, that the distribution of $Y_N$ under ${\tt E}_{NM_{b\mu}^{\ssup {T}}} $, conditioned on $\{Y_{N}(\Gamma_T ^{\ssup 1} \setminus \Gamma_{T,\leq L}^{\ssup 1})=0\}$, is equal to its distribution under ${\tt E}_{NM_{b\mu}^{\ssup {T,\leq L}}}$. This implies that
\begin{equation}\label{lowboundPoi2}
\begin{aligned}
\P_{\Poi_{N\mu}}^{\ssup{N}}&(\Vcal_{N}\in G\cap\Acal_{f,\Aold}) \geq \P_{\Poi_{N\mu}}^{\ssup{N}}(\Vcal_{N}\in G\cap\Acal_{f,<\Aold})\\
&\geq \e^{o(N)}\,{\tt E}_{NM_{b\mu}^{\ssup T}}\Big[\e^{N\phi_b(\smfrac{1}{N}Y_N)}\1\{\smfrac 1N Y_N\in G\cap\Acal_{f,<\Aold}\}\Big| Y_{N}(\Gamma_T ^{\ssup 1} \setminus \Gamma_{T,\leq L}^{\ssup 1})=0 \Big]\\
&\qquad\qquad\times  \Poi_{NM_{\mu}^{\ssup {T}}}\big(Y_{N}(\Gamma_T ^{\ssup 1} \setminus \Gamma_{T,\leq L}^{\ssup 1})=0 \big)\\
&= \e^{o(N)}\,{\tt E}_{NM_{b\mu}^{\ssup {T,\leq L}}}\Big[\e^{N\phi_b(\smfrac{1}{N}Y_N)}\1\{\smfrac 1N Y_N\in G\cap\Acal_{f,<\Aold}\}\Big]\, \\
&\qquad\qquad\times\Poi_{NM_{\mu}^{\ssup {T}}}\big(Y_{N}(\Gamma_T ^{\ssup 1} \setminus \Gamma_{T,\leq L}^{\ssup 1})=0 \big).
\end{aligned}
\end{equation}
Note that
\begin{equation}
 -\delta_L=\liminf_ {N\to\infty}\frac 1N\log  \Poi_{NM_{\mu}^{\ssup {T}}}\big(Y_{N}(\Gamma_T ^{\ssup 1} \setminus \Gamma_{T,\leq L}^{\ssup 1})=0 \big)
\end{equation}
increases to zero as $L\to\infty$. Indeed, this void probability is equal to $\e^{-N M_{b\mu}^{\ssup T}( \Gamma_T ^{\ssup 1} \setminus \Gamma_{T,\leq L}^{\ssup 1})}$, and the rate in the  vanishes as $L\to\infty$.

Recall that $Y_N$ in the expectation on the right-hand side of \eqref{lowboundPoi2} is a PPP with intensity measure $M_{b\mu}^{\ssup {T,\leq L}}$ and hence $Y_N = Y_N^\ssup{\leq L}$ almost surely. Hence, we can rewrite the condition in the indicator on the right-hand side of \eqref{lowboundPoi2} as $\{\smfrac 1N Y_N\in \Pi_L^{-1}(G\cap \Acal_{f,<\Aold})\}$.
Now, $\Pi_L^{-1}(G)$ is open, since $\Pi_L$ is continuous and $G$ is open. Further,  $\Pi_L^{-1}(\Acal_{f,<\Aold})$ is equal to the set of all $\nu\in \Mcal(\Gamma^\ssup{1}_T)$ that satisfy $\int \nu(\d \xi) f(|\xi_0|)\1\{|\xi_0|\leq L\} < \Aold$, which is an open set in $\Mcal(\Gamma^\ssup{1}_T)$, since the map $\xi \mapsto f(|\xi_0|)\1\{|\xi_0|\leq L\}$ is dominated by a continuous and bounded function.

Hence, we can apply now the lower-bound part of  Varadhan's lemma (Lemma 4.3.4 in \cite{DZ10}) and the LDP for $\frac 1N Y_N$ from Lemma~\ref{lem-YLDP} with $\mathfrak m=M_{b\mu}^{\ssup {T,\leq L}}$, to obtain that
\begin{equation}\label{lowboundPoi3}
\begin{aligned}
\liminf_{N\to \infty} &\frac 1N \log \P_{\Poi_{N\mu}}^{\ssup{N}}(\Vcal_{N}\in G\cap\Acal_{f,\Aold}) \\
& \geq -\inf\Big\{H(\nu\mid M_{b\mu}^{\ssup {T,\leq L}})-\phi_b(\nu)\colon \nu\in \Pi_L^{-1}(G\cap\Acal_{f,<\Aold}) \Big\} -\delta_L \\
& \geq -\inf\Big\{H(\nu^\ssup{\leq L}\mid M_{b\mu}^{\ssup {T,\leq L}})-\phi_b(\nu^\ssup{\leq L})\colon \nu\in G\cap\Acal_{f,<\Aold} \Big\} -\delta_L,
\end{aligned}
\end{equation}
where for the last equality we used that $H(\nu\mid M_{b\mu}^{\ssup {T,\leq L}}) = \infty$ if $\nu(\Gamma^\ssup{1}_{T}\setminus \Gamma^\ssup{1}_{T,\leq L}) >0$, and hence it suffices to take the infimum over all $\nu$ that satisfy $\nu = \nu^\ssup{\leq L}$. 
It is easy to see that 
\begin{multline}\label{LtoinftyinVP}
\liminf_{L\to\infty}\inf\Big\{H(\nu^\ssup{\leq L}\mid M_{b\mu}^{\ssup {T,\leq L}})-\phi_b(\nu^\ssup{\leq L})\colon \nu\in G\cap\Acal_{f,<\Aold} \Big\}\\
\leq 
\inf\Big\{H(\nu\mid M_{b\mu}^{\ssup {T}})-\phi_b(\nu)\colon \nu\in G\cap\Acal_{f,<\Aold} \Big\}
= \inf\Big\{H(\nu\mid M_{b\mu}^{\ssup {T}})-\phi_b(\nu)\colon \nu\in G\cap\Acal_{f,\Aold} \Big\}
\end{multline}
where the last step follows from a simple approximation step (approach $\nu$ satisfying $\int\nu(\d\xi)\, f(|\xi_0|)=\Aold$ by $(1-\eps) \nu$ with $\eps\downarrow 0$).

This shows that the complementary lower bound in \eqref{Upperbound} holds for $F$ replaced by an open set $G$, and finishes the proof of Theorem~\ref{thm-LDP}.

The preceding proof of Theorem~\ref{thm-LDP}(1)  can easily be adapted to the proof of the (2)-part. Actually, the proof of the lower bound in case (1) gave already many details of the proof of the lower bound in Case (2). The proof of the upper bound in Case (1) needs more than the proof of the upper bound in Case (2). Finally, the compactness of the sublevel sets of the rate function in \eqref{ILdef} is clear, since $|M_{\mu}^{\ssup {T,\leq L }}|$ is finite (hence, one can take $b=1$; compare to \eqref{IidentEL}), and $\nu\mapsto  \langle \nu,\mathfrak R^{\ssup T}(\nu)\rangle$ is lower semicontinuous. This ends the proof of Theorem~\ref{thm-LDP}.

\section{Analysis of $I_\mu^{\ssup T}$,  (non-)gelation, and the Smoluchowski equation}\label{sec-minimizerproof}

\noindent In this section, we analyse the minimiser(s) of the rate function $I_\mu^{\ssup T}$ appearing in Theorem~\ref{thm-LDP} (defined in \eqref{Idef}) and prove Theorem~\ref{thm-Lossofmass}. In particular, in Section~\ref{sec-bounds} we derive bounds on moments of the reference measure $M_\mu^{\ssup T}$ and lower bounds on $I_\mu^{\ssup T}$ and give criteria for the existence of minimisers of $I_\mu^{\ssup T}$. In Section~\ref{sec-ELeq} we derive the Euler--Lagrange equations for these minimisers and use them to prove some estimates for its moments. Then Section~\ref{sec-subcrphase} is devoted to the proof of non-gelation at small times (finishing the proof of Theorem~\ref{thm-Lossofmass}(1)), and Section~\ref{sec-PhaseTransition} to the proof of loss of mass (i.e., existence of gelation) at a late time (finishing the proof of Theorem~\ref{thm-Lossofmass}(2)). On the way, we also prove Proposition~\ref{prop-criteria}(2) in Section~\ref{sec-ELeq} and Proposition~\ref{prop-criteria}(1) at the end of Section~\ref{sec-subcrphase}, furthermore we derive the Smoluchowski equation in Section~\ref{sec-Smolproof}.

For the remainder of this section, we keep $\mu\in\Mcal_1(\Scal)$ and $T\in (0,\infty)$ fixed and assume only that the kernel $K$ is nonnegative and measurable in its four arguments. Recall the reference measure $M_\mu^{\ssup T }$ from \eqref{Mmudef}. We recall from \eqref{Idef} and \eqref{Imuchar} that, for $\nu\in\Mcal(\Gamma_T^{\ssup 1})$ that is absolutely continuous with respect to $M_\mu^{\ssup T}$ (otherwise, $I^{\ssup T}_\mu( \nu)=\infty$),
\begin{equation}\label{IidentEL}
\begin{aligned}
I^{\ssup T}_\mu( \nu)
&=\Big\langle\nu,\log \frac{\d \nu}{\d M_\mu^{\ssup T}}\Big\rangle+ \frac 12\langle \nu,\mathfrak R^{\ssup T}(\nu)\rangle+ 1- |\nu|\\
&=H(\nu|M^{\ssup T}_{b\mu}) +1- |M_{b\mu}^{\ssup T}|+(1-b) |\nu|+\int\nu_0(\d k)\, |k|\,\log b+ \frac 12 \langle \nu,\mathfrak R^{\ssup T}(\nu)\rangle,
\end{aligned}
\end{equation}
where the second line assumes that $b\in(0,\infty)$ is so small that $|M_{b\mu}^{\ssup T}|<\infty$. (A sufficient criterion for this is given in Lemma~\ref{lem-refmeasure-finite} below.)

We mentioned already below \eqref{entropy} that the map $\nu\mapsto H(\nu|M_{b\mu}^{\ssup T})$ is lower semicontinuous and convex and has compact sublevel sets for any $b\in(0,\infty)$ small enough such that $|M_{b\mu}^{\ssup T}|<\infty$. We proved this in Lemma~\ref{lem-YLDP}. Furthermore, according to Lemma~\ref{Lem-RfrakProp}, also the map $\nu\mapsto\langle \nu,\mathfrak R^{\ssup{T}}(\nu)\rangle$ is lower semicontinuous. However, since $b$ might be less than one, the lower semicontinuity of $I_\mu^{\ssup T}$ is not {\em a priori} clear, since this would require the continuity of $\nu\mapsto \int\nu_0(\d k)\, |k|$. We are able to show this only on suitable subsets (e.g., on $\Acal_{f,\Aold}$) or after restriction to a cut-off version that makes $|k|$ bounded (in general it is not true).

If $K$ is positive definite, then Lemma~\ref{Lem-RfrakProp} implies the strict convexity of $I^{\ssup T}_\mu$ and hence the uniqueness of the minimiser, since the domain of $I^{\ssup T}_\mu$ is convex. However, we are not going to use this assumption in this section, therefore we might have several minimisers.

\subsection{Bounds on $M_\mu^{\ssup T}$ and on $I_\mu^{\ssup T}$}\label{sec-bounds}

Now we can state conditions under which $I_\mu^{\ssup T}$ has a minimiser. 
Recall $q_\mu^{\ssup T}$ from \eqref{qTdef} and that the sublevel sets of $H(\cdot|M)$ are compact for any finite measure $M$ (see the proof in Section~\ref{sec-Yproof}). The following will be used (under the assumption in \eqref{AssK1}) for all sufficiently small $T$.

\begin{lemma}[Moments of $M_\mu^\ssup{T}$ under $q_\mu^\ssup{T}<1$]\label{lem-estiM}
Fix $T\in(0,\infty)$. If $q_\mu^\ssup{T}<1$, then $\int M_\mu^{\ssup T }(\d \xi)\, |\xi_0|^\alpha<\infty$ for any $\alpha\in [0,\infty)$. Furthermore, the sublevel sets of $I_\mu^{\ssup T}$ are then compact, i.e., $I_\mu^{\ssup T}$ has at least one minimiser on $\Mcal(\Gamma_T^{\ssup 1})$.
\end{lemma}
\begin{proof}Just note that
$$
\int M_\mu^{\ssup T }(\d \xi)\, |\xi_0|^\alpha
=\sum_ {n\in\N} M_\mu^{\ssup T }\big(\{\xi\colon |\xi_0|=n\}\big)\,n^\alpha=\sum_ {n\in\N} (q_\mu^\ssup{T})^{n+o(n)}\,n^\alpha,
$$
and a comparison to the geometric series gives the result.

Now we may use the second line of \eqref{IidentEL} for $b=1$ and see that $I_\mu^{\ssup T}$ is equal to the sum of $H(\cdot|M_\mu^{\ssup T})$ (which has compact sublevel sets) and  $1-|M^{\ssup T}_\mu|+\langle \nu,\mathfrak R^{\ssup T}(\nu)\rangle$,  which is lower semi-continuous in $\nu$ by Lemma~\ref{Lem-RfrakProp}. Hence $I_\mu^{\ssup T}$ has  compact sublevel sets as well and possesses therefore a minimiser.
\end{proof}

It follows an elementary and useful bound on the moments of $M_{b\mu}^{\ssup T}$.

\begin{lemma}[Moments of $M_{b\mu}^\ssup{T}$ under  \eqref{AssK1}]\label{lem-refmeasure-finite}
Assume that  \eqref{AssK1} holds and fix any $T,b\in(0,\infty)$. Then for any $n\in \N_0$
\begin{equation}\label{eq:Mesti-n}
M_{b\mu}^\ssup{T}(\{ \xi \in \Gamma^\ssup{1}_T\colon |\xi_0|=n\}) \leq \frac{\e^{1-b}}{2\pi TH} (bTH\e^2)^n n^{-2}
\end{equation}
Consequently, if $b< 1/H\e^2T$, then $\int M_{b\mu}^{\ssup T}(\d \xi)\, |\xi_0|^\alpha<\infty$ for any $\alpha\in [0,\infty)$.
Further, if $TH\e^2<1$, then $q_\mu^\ssup{T} <1$.

\end{lemma}

\begin{proof}
With the help of the estimate  \eqref{tauestiprecise} for $\Q^{\ssup{T}}_k(\Gamma_{T,k}^{\ssup 1})$ derived in Lemma \ref{lem-taubound} under assumption \eqref{AssK1} and also using that $\Q^{\ssup{T}}_0(\Gamma_{T,0}^{\ssup 1})=0$, we obtain that for any $n\in \N_0$
\begin{equation}\label{eq:massrefmeas-esti}
\begin{aligned}
M_{b\mu}^\ssup{T}(\{ \xi \in \Gamma^\ssup{1}_T\colon |\xi_0|=n\})  &=\e \int \Poi_\mu(\d k) \e^{1-b} b^{|k|}\Q^{\ssup{T}}_k(\Gamma_{T,k}^{\ssup 1}) \1\{|k|=n\}\\
&\leq \e^{1-b}\, \frac {b^{n}}{n!}\frac{(TH)^{n-1}}{(n-1)!}\,n^{2(n-1)} \\
&= \frac{\e^{1-b}}{TH} \,\frac {(bTH)^{n}}{(n!)^2}\,n^{2n-1} 
\leq \frac{\e^{1-b}}{2\pi\,TH}\, (bTH\e^2)^{n}n^{-2} .
\end{aligned}.
\end{equation}
where we applied the Stirling bound $n!\geq n^n \e^{-n} \sqrt{2\pi n}$. Hence, \eqref{eq:Mesti-n} holds.
For any $\alpha \in [0,\infty)$ we get that 
\[
\int M_{b\mu}^{\ssup T}(\d \xi)\, |\xi_0|^\alpha\leq \frac{\e^{1-b}}{TH}\sum_{n=0}^\infty \frac {(bTH)^{n}}{(n!)^2}n^{2n-1} 
\leq \frac{\e^{1-b}}{2\pi\,TH} \sum_{n=0}^\infty(bTH\e^2)^{n}n^{\alpha-2}<\infty
\] 
if $bTH\e^2<1$, since the geometric series with that parameter converges. Choosing $b=1$ we get that $q_\mu^\ssup{T} \leq TH\e^2$ which implies the last claim.
\end{proof}


The following lower bound will be used in Section~\ref{sec-PhaseTransition} for proving gelation for large $T$; more precisely for all $T$ such that $I_\mu^{\ssup T}$ is bounded away from zero on $\Mcal(\Gamma_T^{\ssup 1})$.

\begin{lemma}[Lower bound on $I_\mu^{\ssup T}$]\label{lem_lowboundImu}
Under the assumptions in \eqref{AssK1} and  \eqref{AssK2}, \begin{equation}\label{Ilowbound}
\inf_{\nu\in\Mcal(\Gamma_T^{ \ssup 1})} I^{\ssup T}_\mu(\nu)\geq 1- \frac 1{2T}\Big( \frac{\e}{\pi H} + \frac{(\log(2TH\e^2))^2}{h} \Big).
\end{equation}
In particular, the infimum is positive for all sufficiently large $T$ and tends to one as $T\to\infty$.
\end{lemma}

\begin{proof}
For $T\in(0,\infty)$, pick some $b\in(0,1)$ such that $bTH\e^2 <1$. Then from Lemma \ref{lem-refmeasure-finite} with $\alpha=0$ we have (dropping the factor $n^{-2}$ in the sum) 
\begin{equation}\label{eq:Mtotalmass}
|M_{b\mu}^{\ssup T}| \leq \frac{\e^{3}b\e^{-b}}{2\pi\,(1-bTH\e^2)}<\infty.
\end{equation}
We derive from the second line of  \eqref{IidentEL} and the non-negativity of the relative entropy, that for any $\nu\in\Mcal(\Gamma_T^{ \ssup 1})$, the following holds
\begin{equation}\label{estiImu}
\begin{aligned}
I^{\ssup T}_\mu(\nu)&\geq H(\nu|M_{b\mu}^{\ssup T})+  1- |M_{b\mu}^{\ssup T}|+D\log b+ \frac 12 \langle \nu,\mathfrak R^{\ssup T}(\nu)\rangle\\
&\geq 1- |M_{b\mu}^{\ssup T}|+D\log b+ \frac 12 \langle \nu,\mathfrak R^{\ssup T}(\nu)\rangle,
\end{aligned}
\end{equation}
where we abbreviated $D=\int \nu_0(\d k)\,|k| =|c_{\nu_0}|$. With the help of  \eqref{AssK2}, we obtain 
\begin{equation}\label{eq:Rloweresti}
\begin{aligned}
\langle \nu,\mathfrak R^{\ssup T}(\nu)\rangle
&=\int_0^T \d t\int \nu(\d \xi)\int \nu(\d \xi')\,\langle \xi_t, K\xi_t'\rangle
\geq \int_0^T \d t\int \nu(\d \xi)\int \nu(\d \xi')\,h \|\xi_t\|_1\,\|\xi_t'\|_1\\
&=h\int_0^T \d t\int \nu(\d \xi)\int \nu(\d \xi')\,|\xi_0|\,|\xi_0'| 
= hT D^2.
\end{aligned}
\end{equation}
The polynomial $D\mapsto D\log b+ \frac 12 hT D^2$ assumes its minimal value at $D=\frac{-\log b}{hT}$ with value $-\frac{(\log b)^2}{2hT}$. Hence, 
$$
\inf_{\nu\in\Mcal(\Gamma_T^{ \ssup 1})} I^{\ssup T}_\mu(\nu)\geq 1-\frac{\e^{3}b\e^{-b}}{2\pi\,(1-bTH\e^2)}-\frac{(\log b)^2}{2hT},\qquad b\in(0,\smfrac 1{TH \e^2}).               
$$
Picking $b = 1/(2TH\e^2)$, we get
\[
\begin{aligned}
\inf_{\nu\in\Mcal(\Gamma_T^{ \ssup 1})} I^{\ssup T}_\mu(\nu)&\geq 1- \frac 1{2T}\Big( \frac{\e^{1-(2TH\e^2)^{-1}}}{\pi H} + \frac{(\log(2TH\e^2))^2}{h} \Big)\\
&\geq 1- \frac 1{2T}\Big( \frac{\e}{\pi H} + \frac{(\log(2TH\e^2))^2}{h} \Big).
\end{aligned}
\]
\end{proof}

\subsection{Euler--Lagrange equations}\label{sec-ELeq}

In this section we characterise minimisers of  $I^{\ssup T}_\mu$ via the variational equalities, which we also call {\em Euler--Lagrange equations}. This  will also lead to a proof of Proposition~\ref{prop-criteria}(2).

\begin{lemma}\label{lem_ELequations}
For any $T\in(0,\infty)$, any minimiser $\nu^{\ssup T}$ of $I^{\ssup T}_\mu$ on the set $\Mcal(\Gamma_T^{\ssup 1})$ satisfies the Euler--Lagrange equation
\begin{equation}\label{ELPoissonized}
\nu^{\ssup T}(\d \xi)= M_\mu^{\ssup T}(\d \xi)\e^{-\mathfrak R^{\ssup T}(\nu^{\ssup T})(\xi)},\qquad\xi\in\Gamma_T^{\ssup 1}.
\end{equation}
Furthermore, for any $L\in\N$, any minimiser $\nu^{\ssup {T,\leq L}}$ of $I^{\ssup {T,\leq L}}_\mu$ defined in \eqref{ILdef} on the set $\Mcal(\Gamma_{T,\leq L}^{\ssup 1})$ satisfies the Euler--Lagrange equation
\begin{equation}\label{ELwithL}
\nu^{\ssup T}(\d \xi)= M_\mu^{\ssup {T,\leq L}}(\d \xi)\e^{-\mathfrak R^{\ssup T}(\nu^{\ssup T})(\xi)},\qquad\xi\in\Gamma_{T,\leq L}^{\ssup 1}.
\end{equation}
\end{lemma}

\begin{proof} We drop the superscript $^\ssup{T}$ and the index $\mu$ from the notation.

Let $\nu$ be a minimiser of $I$. Since $I(\nu)$ is finite, $\nu$ has a nonnegative density $\varphi$ with respect to $M$. We show that $\varphi$ is positive $M$-almost surely. Indeed, if there is a measurable set $B\subset\Mcal(\Gamma_T^{\ssup 1})$ with positive  $M$-measure such that $\varphi=0$ on $B$, then we are going to see that $\nu_\eps(\d \xi)=M(\d \xi)(\varphi(\xi)+\eps \1_B(\xi))$ has a strictly smaller $I$-value, in contradiction to the minimality of $\nu$. Indeed, observe that
$$
\begin{aligned}
I(\nu_\eps)-I(\nu)&=\frac 12\int\int M(\d\xi) M(\d\widetilde \xi)R(\xi,\widetilde \xi)\Big[2\eps\1_B(\xi)\varphi(\widetilde \xi)+\eps^2\1_B(\xi)\1_B(\widetilde \xi)\Big] \\
&\qquad +M(B)\eps\log \eps- \eps M(B).
\end{aligned}
$$
Since this is $\leq M(B)\eps\log \eps +O(\eps)$ for $\eps \downarrow 0$, it is negative for sufficiently small $\eps>0$. Hence, $\varphi$ is positive $M$-almost surely.

Now we calculate the directional derivative of $I$ in $\nu(\d \xi)=M(\d\xi)\,\varphi(\xi)$ in direction of $\nu_\eps(\d\xi)=M(\d\xi)(\varphi(\xi)+\eps\gamma(\xi))$ (with $\eps\in\R$) for a large class of measurable and bounded functions $\gamma\colon \Gamma_T^{\ssup 1}\to\R$. We fix $\delta>0$ and $L\in\N$ and assume that $\gamma=0$ on $\{\xi\colon |\xi_0|> L, \varphi(\xi)\leq \delta\}$. Then $\varphi+\eps\gamma>0$ for all $\eps\in \R$ with sufficiently small $|\eps|$. By minimality in $\eps=0$, we have
\begin{equation}\label{Imuderivative}
\begin{aligned}
0&=\frac{\d}{\d \eps} I(\nu_\eps) \Big|_{\eps = 0} 
=\frac{\d}{\d \eps}\Big(\big\langle M(\varphi+\eps\gamma), \log( \varphi+\eps\gamma)\big\rangle +\frac 12 \big\langle M(\varphi+\eps\gamma),\mathfrak R(M(\varphi+\eps\gamma))\big\rangle\\
&\qquad - \langle M(\varphi+\eps\gamma),\1\big\rangle\Big)\Big|_{\eps = 0}\\
&= \big\langle M \gamma, \log \varphi\big\rangle +\langle M \gamma, \1\rangle+ \langle M\gamma, \mathfrak R(\nu)\rangle -\langle M \gamma, \1\rangle\\
&= \Big\langle M \gamma, \log \varphi+ \mathfrak R(\nu) \Big\rangle.
\end{aligned}
\end{equation}
Since this holds for any bounded measurable function $\gamma$ with $\supp(\gamma) \subset \{\xi\colon |\xi_0|\leq L\mbox{ or }\varphi(\xi)>\delta\}$, we obtain that
$$
0= \log \varphi(\xi)+ \mathfrak R(\nu)(\xi), \qquad M\mbox{-almost surely},
$$
first only on the set $\{\xi\colon |\xi_0|\leq L\mbox{ or }\varphi(\xi)>\delta\}$, and hence on 
\[\bigcup_{L\in\N}\{\xi\colon |\xi_0|\leq L\} \cup \bigcup_{\delta>0}\{\xi\colon \varphi(\xi)>\delta\}= \Gamma^\ssup{1}_T \cup \{\varphi>0\},
\]
which is equal to $\Gamma^\ssup{1}_T$ $M$-almost surely. 

This implies the claim in \eqref{ELPoissonized}. The proof of \eqref{ELwithL} is analogous.
\end{proof}

\begin{lemma}[Bounds on $\nu^{\ssup T}$]\label{lem-estinuT}
Assume that  $\nu^{\ssup T}$ is a minimiser of $I^{\ssup T}_\mu$ on $\Mcal(\Gamma_T^{\ssup 1})$ or a minimiser of $I_\mu^{\ssup {T,\leq L}}$ on $\Mcal(\Gamma_{T,\leq L}^{\ssup 1})$ for some $L\in\N$. 
\begin{enumerate}
\item Under the assumption in \eqref{AssK1}, and if $T<1/\e^2 H$,
\begin{equation}\label{2ndmomentnufinite}
\int \nu^{\ssup T}(\d \xi)\, |\xi_0|^2\leq \frac {\e^2}{2 \pi(1-\e^2 TH)}.
\end{equation}
\item Under the assumptions in \eqref{AssK1} and \eqref{AssK2}, for any $T>0$, 
\begin{equation}\label{uppboundX}
\int \nu^{\ssup T}(\d \xi)\, |\xi_0|\leq \max\left\{\frac{1}{hT}\log(2 TH\e^2), \frac 1{2\pi HT} \right\}.
\end{equation}
\end{enumerate}
\end{lemma}

\begin{proof} Again, we write $\nu$ instead of $\nu^{\ssup T}$ and $M$ instead of $M_\mu^{\ssup T}$ and $\mathfrak R$ instead of $\mathfrak R^{\ssup T}$. 

We start with (1). From the EL-equations in \eqref{ELPoissonized}, we see that 
$\int \nu(\d\xi)\,|\xi_0|^2\leq \int M(\d \xi)\,|\xi_0|^2$. Since $M^{\ssup {\leq L}}\leq M$, we only need to prove the statement for $I^{\ssup T}_\mu$.
We assume only \eqref{AssK1}. Under the assumption \eqref{AssK1}, applying Lemma~\ref{lem-refmeasure-finite} for $\alpha=2$ and $b=1$ finishes the proof of \eqref{2ndmomentnufinite}, recalling \eqref{eq:Mtotalmass}.

We continue with (2).  We handle simultaneously the minimisers of $I^\ssup{T,\leq L}_\mu$, for $L\in \N$, and $I^\ssup{T}_\mu$ and denote them by $\nu^\ssup{\leq L}$ for both $L\in \N$ and $L=\infty$. We abbreviate $D^\ssup{\leq L}= \int \nu^\ssup{\leq L}(\d \xi) |\xi_0| = \int M^\ssup{\leq L} (\d \xi) \e^{-\mathfrak{R}(\nu^\ssup{\leq L})(\xi)} |\xi_0|$ (using \eqref{ELPoissonized}, respectively  \eqref{ELwithL}.)  If $D^\ssup{\leq L}\leq \frac{1}{hT}\log(2TH\e^2)$, then we are done. If the converse is true, then $TH\e^2 \e^{-ThD^\ssup{\leq L}}<\frac 12$. We begin by noting that  
$$
\mathfrak R(\nu^{\ssup{\leq L}})(\xi)\geq  Th |\xi_0| D^{\ssup{\leq L}}
$$
holds under the assumption in \eqref{AssK2} and is derived using the same steps as in \eqref{eq:Rloweresti} (without the additional integration over $\nu(\d \xi)$). This already implies that 
\[D^{\ssup{\leq L}}\leq \int M(\d \xi)\,\e^{- Th |\xi_0| D^{\ssup{\leq L}}}|\xi_0|.\]
This upper bound can be further estimated from above. Indeed, using the same arguments as in \eqref{eq:massrefmeas-esti} we obtain that 
\begin{equation}\label{xex-upperbd}
\begin{aligned}
D^{\ssup{\leq L}} &\leq 
\frac 1{2\pi HT}\sum_{n=1}^\infty \frac 1n \big(TH\e^2 \e^{-ThD^{\ssup{\leq L}}}\big)^n
=-\frac 1{2\pi HT}\log\big(1- TH\e^2 \e^{-ThD^{\ssup{\leq L}}}\big)\\
&\leq \frac 1{2\pi HT} \frac{TH\e^2 \e^{-ThD^{\ssup{\leq L}}}}{1-TH\e^2 \e^{-ThD^{\ssup{\leq L}}}} \leq  \frac 1{2\pi HT} ,
\end{aligned}
\end{equation}
where we used that, due to $TH\e^2 \e^{-ThD^{\ssup{\leq L}}}<\frac 12<1$, we can apply the formula $\sum_{n=1}^\infty \frac{1}{n}q^n = - \log(1-q)$ that holds for $q\in [0,1)$. Then, we used the estimate $-\log(1-x) \leq \frac{x}{1-x}$, for $x<1$ and after that we used the monotonicity of $x\mapsto \frac{x}{1-x}$. Hence, we proved \eqref{uppboundX}. 
\end{proof}

\begin{proof}[Proof of Proposition~\ref{prop-criteria}(2)]
Assertion (a) is shown in Lemma \ref{lem_lowboundImu} and assertion (b) is from Lemma \ref{lem-estinuT}.
\end{proof}

We see also from \eqref{2ndmomentnufinite} that, for any minimiser $\nu^{\ssup T}$ of $I_\mu^{\ssup T}$, we have that $\nu^{\ssup T}\in\Acal_{f,\Aold}$ for $f(r)=r^2$, all $T\in(0,1/\e^2H)$ and any sufficiently large $\Aold$. Here is another benefit from \eqref{2ndmomentnufinite}:

\begin{lemma}[Uniqueness of solutions to EL-equations]\label{lem-ELunique}
Assume that $K$ satisfies \eqref{AssK1}. For any $T\in(0,\frac 1{H \e^2\,}\frac\pi{1+\pi})$, there is at most one solution $\nu$ to \eqref{ELPoissonized}.
\end{lemma}

\begin{proof} Assume that $ \nu$ and $\widetilde \nu$ are two solutions to the EL equation in  \eqref{ELPoissonized}. Using the estimate $|\e^{-x}-\e^{-y}|\leq |x-y|\min\{\e^{-x},\e^{-y}\}\leq |x-y|(\e^{-x}+\e^{-y})$ for $x,y\in\R$, we obtain that
$$\begin{aligned}
\int |\nu-\widetilde \nu|(\d\xi)\,|\xi_0|&=\Big|\int M_\mu^{\ssup T}(\d \xi)\,|\xi_0|\,\Big(\e^{-\mathfrak R^{\ssup T}(\nu)(\xi)}-\e^{-\mathfrak R^{\ssup T}(\widetilde \nu)(\xi)}\Big)\Big|\\
&\leq\int M_\mu^{\ssup T}(\d \xi)\,|\xi_0|\,|\mathfrak R^{\ssup T}(\nu-\widetilde \nu)(\xi)|\Big(\e^{-\mathfrak R^{\ssup T}(\nu)(\xi)}+\e^{-\mathfrak R^{\ssup T}(\widetilde \nu)(\xi)}\Big).
\end{aligned}
$$
Now we use \eqref{AssK1} to get 
$$
|\mathfrak R^{\ssup T}(\nu-\widetilde \nu)(\xi)|\leq \mathfrak R^{\ssup T}(|\nu-\widetilde \nu|)(\xi)\leq HT \int |\nu-\widetilde \nu|(\d\widetilde \xi)\,|\widetilde \xi_0|\,|\xi_0|.
$$
From now on we assume that $\e^2TH<1$. Then, we can combine that last two estimates, use the EL equations again, and afterwards \eqref{2ndmomentnufinite} and obtain
\[
\begin{aligned}
\int |\nu-\widetilde \nu|(\d\xi)\,|\xi_0|
&\leq HT \int |\nu-\widetilde \nu|(\d\widetilde \xi)\,|\widetilde \xi_0|\,\Big(\int\nu(\d\xi)\,|\xi_0|^2+ \int\widetilde \nu(\d\xi)\,|\xi_0|^2\Big)\\
&\leq 2 HT \int |\nu-\widetilde \nu|(\d\xi)\,| \xi_0|\,\frac {\e^2}{2\pi (1-\e^2 TH)}.
\end {aligned}
\]
If $T$ is so small that $\e^2 TH<1$ and $\frac {\e^2HT}{\pi(1-\e^2 TH)}<1$, this  implies that $\int |\nu-\widetilde \nu|(\d\xi)\,|\xi_0|=0$, which implies that $\nu=\widetilde \nu$. The condition on $x=HT$ reads $0<x<1/\e^2$ and $x<\frac \pi{\e^2(1+\pi)}$. Hence, the latter inequality holds for $TH\in(0,\frac 1{\e^2}\,\frac \pi{1+\pi })$. This implies the assertion.
\end{proof}

\subsection{Subcritical phase: convergence and non-gelation}\label{sec-subcrphase}

In this section, we provide the proofs of Theorem~\ref{thm-Lossofmass}(1) and Proposition~\ref{prop-criteria}(1).
Throughout the section, we fix $T>0$ and $\mu\in \Mcal_1(\Scal)$. Note that Lemma \ref{lem-estiM} already covers Theorem~\ref{thm-Lossofmass}(1)(a) about the compact sublevel sets of $I_{\mu}^\ssup{T}$ and the existence of minimisers, and Lemma~\ref{lem_ELequations} implies Assertion (c) about the validity of the Euler--Lagrange equations for minimisers. 

The outline of this section is as follows. The tightness assertion about $\Vcal_N^\ssup{T}$ under $\P_{\Poi_{N\mu}}^\ssup{N}$ in Theorem~\ref{thm-Lossofmass}(1)(d) is proved in Lemma \ref{subcrLLN}, and the tightness of $c_{\Vcal_{N,0}^\ssup{T}}$, as well as  Theorem~\ref{thm-Lossofmass}(1)(e) is proved in Corollary \ref{subcrLLNmass}. Finally, the asssertion about non-gelation in Theorem~\ref{thm-Lossofmass}(1)(b) is proved in Corollary \ref{subcrfullmass}. The proof of Proposition~\ref{prop-criteria}(1) is finished at the end of this section.

Let us start by explaining the strategy for proving the tightness result. Usually tightness is directly implied by the LDP, as we stated in Corollary \ref{cor-accu}. The problem is that in our LDP from Theorem~\ref{thm-LDP} we conditioned on $\{\Vcal_N^\ssup{T}\in \Acal_{f,\Aold}\}$. However, we want to prove tightness for the unconditioned distribution of $\Vcal_N^\ssup{T}$. Note that we are free to choose $f(r)=r^2$. For that particular choice can argue that the probability of the event $\{\Vcal_N^\ssup{T}\notin \Acal_{f,\Aold}\}$ vanishes (see Lemma \ref{lem-secmom}), for large $\Aold$. Further, we can show that the minimisers of $I_\mu^\ssup{T}$ are also minimisers of the $\Aold$-dependent rate function from \eqref{ratefunctionA}, if $\Aold$ is large enough. Finally, we will use this to argue that the unconditioned distribution of $\Vcal_N^\ssup{T}$ converges to a distribution that is concentrated on minimisers of $I_\mu^\ssup{T}$.

\begin{lemma}\label{lem-secmom}
Let $T>0$ and $\mu\in \Mcal_1(\Scal)$ be such that $q_\mu^\ssup{T} <1$. Then 
\begin{equation}
\sup_{N\in \N} \P_{\Poi_{N\mu}}^\ssup{N}\big(\Vcal_N^\ssup{T}\notin \Acal_{f,\Aold}\big) \leq \frac{C}{\Aold},\qquad \Aold\in (0,\infty),
\end{equation}
where $C = \int M_\mu^{\ssup T}(\d \xi)\, |\xi_0|^2 <\infty$.
\end{lemma}

\begin{proof} By  Markov inequality, it is enough to show that the expectation of $\int\Vcal_{N,0}^\ssup{T}(\d k)\, |k|^2$ under $\E_{\Poi_{N\mu}}^\ssup{N}$ is bounded in $N$.  Abbreviate ${\tt P}_N = {\tt P}_{N M_{\mu}^\ssup{T}}$. Applying Corollary \ref{DistVcal-rescaledkernel} with $b=1$ and choosing $f$ as the constant function $\nu\mapsto 1$, gives us 
\begin{align*}
1 = {\tt E}_N\big[\e^{-\frac {1}{2N} \sum_{i,j \colon i\not=j}R(\Xi_i,\Xi_j)}\e^{-\frac 1N \int\varphi_{\xi_0}(\xi)\,Y_N(\d\xi)}\big] \e^{N(|M^\ssup{T}_{\mu}|-1)}.
\end{align*}
Notice that, since $q_\mu^\ssup{T} <1$, we have that $|M^\ssup{T}_{\mu}|<\infty$ and the expression above is well-defined for $b=1$. Then, we apply Corollary \ref{DistVcal-rescaledkernel} a second time with $b=1$ to the function $\nu\mapsto \int \nu_0(\d k) |k|^2$ and combine the formulas to get
\begin{equation}\label{secmom-eq1}
\begin{aligned}
\E_{\Poi_{N\mu}}^\ssup{N}\Big[\int\Vcal_{N,0}^\ssup{T}(\d k)\, |k|^2\Big]
&=\frac{{\tt E}_N\Big[\int\smfrac{1}{N}Y_{N,0}(\d k)\, |k|^2\,\e^{-\frac{1}{2N} \sum_{i,j \colon i\not=j}R(\Xi_i,\Xi_j)}\e^{-\frac 1N \int\varphi_{\xi_0}(\xi)\,Y_N(\d\xi)}\Big]}{{\tt E}_N\big[\e^{-\frac {1}{2N} \sum_{i, j \colon i\not=j}R(\Xi_i,\Xi_j)}\e^{-\frac 1N \int\varphi_{\xi_0}(\xi)\,Y_N(\d\xi)}\big]}\\
& = \frac{{\tt E}_N \big[f_N(Y_N) g_N(Y_N)\big]}{{\tt E}_N\big[g_N(Y_N)\big]}.
\end{aligned}
\end{equation}
where, for $\nu= \sum_i \delta_{\xi_i}\in \Mcal_{\N_0}(\Gamma^\ssup{1}_T)$, we defined 
\begin{align*}
f_N(\nu) = \int \smfrac{1}{N}\nu_0(\d k) \,|k|^2, \quad g_N(\nu)= \e^{-\frac{1}{2N}\sum_{i,j \colon i\neq j} R(\xi_i,\xi_j)}\e^{-\frac 1N \int\varphi_{\xi_0}(\xi)\,\nu(\d\xi)}.
\end{align*}
Observe that $f_N$ and $-g_N$ are increasing on $\Mcal_{\N_0}(\Gamma^\ssup{1}_T)$ under the addition of points. Thus, we can apply the Harris-FKG inequality (see Theorem 20.4 in \cite{LaPe17}) to bound the right-hand side of \eqref{secmom-eq1} by
\begin{equation}\label{esti-2ndmom}
\begin{aligned}
{\tt E}_N \big[f_N(Y_N)\big] = {\tt E}_N\Big[\int\smfrac{1}{N}Y_{N,0}(\d k)\, |k|^2\Big] = \int_{\Gamma^\ssup{1}_T}M_\mu^{\ssup T}(\d \xi)\, |\xi_0|^2,
\end{aligned}
\end{equation}
where we used Campbell's formula (see Proposition 2.7 in \cite{LaPe17}). Note that the right-hand side is finite under the assumption $q^\ssup{T}_\mu < 1$ due to Lemma \ref{lem-estiM}. Hence, we have a established a bound that is uniform in $N\in \N$.
This finishes the proof.
\end{proof}

It is standard (see Corollary~\ref{cor-accu}) that, given the LDP of Theorem~\ref{thm-LDP}, accumulation points of $(\Vcal^\ssup{T}_N)_{N\in\N}$ exist and are concentrated on the set of minimisers of the rate function. However, this holds {\em a  priori} only under conditioning on $\Vcal^\ssup{T}_{N} \in \Acal_{f,\Aold}$. However, we now derive  that under $q_\mu^\ssup{T}<1$ this holds under the unconditioned measure $\P_{\Poi_{N\mu}}^\ssup{N}$ as well.

\begin{lemma}[Law of large numbers]\label{subcrLLN}
Fix $T>0$ and assume that $q_\mu^\ssup{T}<1$. Let $\P_{\Vcal_N}$  denote the distribution of $\Vcal_N^\ssup{T}$ under $\P^\ssup{N}_{\Poi_{N\mu}}$. Then the sequence of measures $(\P_{\Vcal_N})_{N\in \N}$ is tight (and thus relatively compact) and each limit point $\P$ is concentrated on the set of minimisers of $I_\mu^{\ssup T}$, i.e.,
\begin{equation}\label{support-limdistri}
\supp( \P )\subset D_0 := \{\nu \in \Mcal(\Gamma^\ssup{1}_T)\colon I_\mu^\ssup{T}(\nu) = \inf I_\mu^{\ssup T}\} .
\end{equation}
\end{lemma}

\begin{proof} Recall that we are working with $f(r)=r^2$. Abbreviate the distribution of $\Vcal_N^\ssup{T}$ under $\P^\ssup{N}_{\Poi_{N\mu}}(\cdot \mid \Vcal_N^\ssup{T} \in \Acal_{f,\Aold})$ by $\P_{\Vcal_N, \Aold}$. According to Theorem~\ref{thm-LDP}, $(\P_{\Vcal_N, \Aold})_{N\in\N}$ satisfies an LDP on $\Acal_{f,\Aold}$ with good rate function $I_{\Aold,\mu}^{\ssup T}$ given by $I_{\Aold,\mu}^{\ssup T}(\nu)=I_{\mu}^{\ssup T}(\nu)-\chi_\Aold$ for $\nu\in \Acal_{f,\Aold}$. According to Lemma~\ref{lem-estiM}, because of $q_\mu^\ssup{T}<1$, $I_\mu^{\ssup T}$ possesses at least one minimiser on $\Mcal(\Gamma_T^{\ssup 1})$. For sufficiently large $\Aold$, every minimiser $\nu$ of $I_\mu^{\ssup T}$ lies in $\Acal_{f,\Aold}$, since it satisfies the EL-equation in \eqref{ELPoissonized} and satisfies therefore the estimate in \eqref{2ndmomentnufinite} for the second moment, which does not depend on $\nu$. Hence, $\{I_{\Aold,\mu}^{\ssup T}=0\}=D_0$ for all sufficiently large $\Aold$.

As a consequence of Corollary~\ref{cor-accu},  $(\P_{\Vcal_N, \Aold})_{N\in\N}$  is tight, and any accumulation point $\P_\Aold$ is concentrated on $\{I_{\Aold,\mu}^{\ssup T}=0\}$, that is, on $D_0$. We now show that $(\P_{\Vcal_N, \Aold})_{N\in\N}$ and $(\P_{\Vcal_N})_{N\in\N}$ have the same limiting behaviour, that is, they are tight and every accumulation point is concentrated on $ D_0$. Indeed, for any open neighbourhood $U$ of $D_0$, we have
$$
\P_{\Vcal_N}(U^{\rm c})\leq \P_{\Vcal_N}(U^{\rm c}\cap\Acal_{f,\Aold})+\P_{\Vcal_N}(\Acal_{f,\Aold}^{\rm c})
\leq \P_{\Vcal_N,\Aold}(U^{\rm c})+\frac C\Aold,
$$
according to Lemma~\ref{lem-secmom}, where $C/\Aold$ can be taken arbitrarily small, uniformly in $N$. Using the above, we see that $\P_{\Vcal_N}(U^{\rm c})$ vanishes as $N\to\infty$. Hence, any accumulation point $\P$ of $(\P_{\Vcal_N})_{N\in\N}$ is concentrated on $U$ and hence on $ D_0$. 
\end{proof}

%

\begin{cor}[Law of large numbers of the total mass]\label{subcrLLNmass}
Fix $T>0$ and assume that $q_\mu^\ssup{T}<1$. Let $\P_{\Vcal_N}$  denote the distribution of $\Vcal_N^\ssup{T}$ under $\P^\ssup{N}_{\Poi_{N\mu}}$ and recall that $(\P_{\Vcal_N})_{N\in \N}$ is tight by Lemma \ref{subcrLLN}. Take a subsequence, also denoted  $(\P_{\Vcal_N})_{N\in \N}$, with limit point $\P$ and let $\Vcal$ be a random variable with distribution $\P$. Recalling the definition of $c_{\nu_0}$ from \eqref{clambdadef} and that $|c_{\nu_0}| = c_{\nu_0}(\Scal)$, for $\nu\in \Mcal(\Gamma^\ssup{1}_T)$, we have that
\begin{equation}\label{eq:cdistconv}
\big|c_{\Vcal^{\ssup {T}}_{N,0}}\big| \to |c_{\Vcal_0}|\text{ in distribution, as } N\to \infty
\end{equation}
Further, we have that $|c_{\Vcal^{\ssup {T}}_{N,0}}|\to 1$ in probability, as $N\to \infty$. Hence $|c_{\Vcal_0}|=1$ $\P$-almost surely. 
\end{cor}

\begin{proof} 
We start with proving \eqref{eq:cdistconv}. In the proof of Lemma \ref{subcrLLN} we have seen that $\P$ is concentrated on $D_0\subset  \Acal_{f,\Aold}$, for $\Aold$ large enough.
Since the map $\nu\mapsto |c_{\nu_0}|$ is continuous on $\Acal_{f,\Aold}$, the claim follows by Lemma \ref{subcrLLN} and  the continuous mapping theorem.


It remains to argue the last statement. Recall that $N c_{\Vcal^\ssup{T}_{N,0}}$ is equal to the number of atoms $n(0)$ of the coagulation process, which is $\Poi_N$-distributed under $\P_{\Poi_{N \mu}}^{\ssup N}$. By the law of large numbers $c_{\Vcal^\ssup{T}_{N,0}} \to 1$ in probability. We combine this with the result \eqref{eq:cdistconv}, which implies that $\P(||c_{\Vcal_0}|-1|>\eps)\leq \liminf_{N\to \infty}\P_{\Poi_{N\mu}}^\ssup{N}(||c_{(\Vcal^\ssup{T}_N)_0}|-1|>\eps)=0$ for any $\eps>0$. Hence, $|c_{\Vcal_0}| =1$ $\P$-almost surely.
\end{proof}

Recall the definition of $\G_T^{\ssup \mu}$ from \eqref{nongel} and observe that via \eqref{totalmass-identity} we have that 
\[
\G_T^\ssup{\mu} = \lim_{L\to\infty}\limsup_{N\to\infty}\E^{\ssup N}_{\Poi_{N\mu}}\Big[ \Big|c_{(\Vcal_N^\ssup{T})_0}^\ssup{\leq L} \Big|\Big]
\]

\begin{cor}[No gelation if $q^\ssup{T}_\mu <1$]\label{subcrfullmass}
Fix $T\in(0,\infty)$ and $\mu  \in\Mcal_1(\Scal)$ and  assume that $q^\ssup{T}_\mu <1$. Then $\G_T^{\ssup \mu}=1$. 
\end{cor}

\begin{proof}
Since $|c_{\Vcal^\ssup{T}_{N,0}}^\ssup{\leq L}| \leq |c_{\Vcal^\ssup{T}_{N,0}}|$ holds for any $L,N\in\N$, it is already clear that $\G_T^{\ssup \mu}\leq 1$.

Now, we argue that $\G_T^{\ssup \mu}\geq 1$.
Fix any limit point $\P$ of $(\P_{\Vcal_N})_{N\in \N}$. Then 
\begin{align*}
\lim_{L\to \infty} \liminf_{N\to \infty} \E_{\Vcal_N}[|c_{\Vcal_{N,0}}^\ssup{\leq L}|] &\geq \lim_{L\to \infty} \E_\Vcal[|c_{\Vcal_{0}}^\ssup{\leq L}|] = \E_\Vcal[|c_{\Vcal_{0}}|] = 1,
\end{align*}
where we used monotone convergence in $L$ for the first equality, Corollary \ref{subcrLLNmass} for the second one. 
\end{proof}

This finishes the proof of Theorem~\ref{thm-Lossofmass}(1).

\begin{proof}[Proof of Proposition~\ref{prop-criteria}(1) ]
By Lemma \ref{lem-refmeasure-finite} we know that $TH<1/\e^2$ implies $q_\mu^\ssup{T}<1$, which gives us assertion (a).  
If $TH<\frac{1}{\e^2}\frac{\pi}{1+\pi}$, then the EL equation \eqref{ELPoissonized} has a unique solution $\nu^\ssup{T}$ according to Lemma \ref{lem-ELunique}, which gives the first assertion of (b). Then, by Lemma \ref{subcrLLN} every limit point $\P$ of $\P^\ssup{N}_{\Poi_{N\mu}}(\Vcal_N^\ssup{T}\in \, \cdot \,)$ is a probability measure concentrated on $\{\nu^\ssup{T}\}$, i.e. the only possible limit is $\P= \delta_{\nu^\ssup{T}}$. Consequently, $\P^\ssup{N}_{\Poi_{N\mu}}(\Vcal_N^\ssup{T}\in \, \cdot \,)\Longrightarrow\delta_{\nu^{\ssup T}}$ as $N\to \infty$. Now Lemma~\ref{lem-Contrho} and \eqref{Xi_as_function_of_nu} imply the last statement of (b).
\end{proof}

\subsection{Supercritical phase: loss of mass and gelation}\label{sec-PhaseTransition}

In this section we assume that the upper bound on the kernel $K$ in \eqref{AssK1} holds as well as the lower bound \eqref{AssK2},  and  show that gelation occurs under the assumption that $\kappa:=\inf_{\nu\in\Mcal(\Gamma_T^{\ssup 1})}I_\mu^{\ssup T}(\nu)>0$. That is, for the process $\Vcal_N^\ssup{T}$ we can observe mass loss in the sense that $\G_T^{\ssup {\mu }}<1$; see \eqref{nongel}, which is the assertion (a) in  Theorem~\ref{thm-Lossofmass}(2). Assertion (b) is implied by Lemma \ref{lem_ELequations}.

Fix $\mu\in\Mcal_1(\Scal)$. Pick any $b\in(0,1]$, such that the total mass of $M_{b\mu}^{\ssup T}$ is finite. Recall the equality about the $L\to \infty$ limit from \eqref{LtoinftyinVP} and note that this also holds if the infimum is taken over $\nu\in \Mcal(\Gamma^\ssup{1}_T)$. This implies that $L$ can be picked so large that $\kappa^{\ssup{\leq L}}=\inf_{\nu\in\Mcal(\Gamma_{T,\leq L}^{\ssup 1})}(H(\nu|M^{\ssup{T,\leq L}}_{b\mu})-\phi_{b}(\nu))\geq\frac \kappa 2$, where we recall the definition of $\phi_b$ from \eqref{phiNdef} and the fact that $H(\nu|M^{\ssup{T}}_{b\mu})-\phi_{b}(\nu)=I_\mu^{\ssup T}(\nu)$.

Note that $|c^{\ssup{\leq L}}_{\Vcal_{N,0}^{\ssup {T,N}}}| + |c^{\ssup{> L}}_{\Vcal_{N,0}^{\ssup {T}}}| $ (with obvious notation for the second quantity) is equal to $\frac 1N$ times a $\Poi_N$-distributed variable under $\E_{\Poi_{N\mu}}^\ssup{N}$. Hence, for any $\eps\in(0,1)$, on the event $\{|c^{\ssup{> L}}_{\Vcal_{N,0}^{\ssup {T}}}| >\eps\}$, the expectation of $|c^{\ssup{\leq L}}_{\Vcal_{N,0}^{\ssup {T}}}|$ is not larger than $1-\eps+o(1)$. Hence, it suffices to show that, for some $\eps\in(0,1)$,
\begin{equation}
\limsup_{N\to\infty}\E_{\Poi_{N\mu}}^\ssup{N}\Big[\big|c^{\ssup{\leq L}}_{\Vcal_N^\ssup{T}}\big|\1\big\{|c^{\ssup{> L}}_{\Vcal_{N,0}^{\ssup {T}}}| \leq \eps\}\Big]=0.
\end{equation}
Indeed, we will show that this expectation decays even exponentially fast. To show this, we apply Corollary~\ref{DistVcal-rescaledkernel} to obtain
$$
\begin{aligned}
\E_{\Poi_{N\mu}}^\ssup{N}\Big[&\big|c^{\ssup{\leq L}}_{\Vcal_N^\ssup{T}}\big|\1\big\{|c^{\ssup{> L}}_{\Vcal_{N,0}^{\ssup {T}}}| \leq \eps\}\Big]
={\tt E}_{N M_{b\mu}^{\ssup T}}\Big[\e^{-\frac 1{2N}\sum_{i\not= j}R^{\ssup T}(\Xi_i,\Xi_j)} b^{-N{D_N}} \e^{(b-1)|Y_{N,0}|}\\
&\qquad\times \e^{- \int\varphi_{\xi_0}(\xi)\,\frac 1N Y_N(\d\xi)} D_N^{\ssup{ \leq L}}\1\{D_N^{\ssup{>L}}\leq \eps\}\Big]\e^{N(|M_{b\mu}^{\ssup T}|-1)},
\end {aligned}
$$
where we introduced the abbreviation $D_N^{\ssup{\leq L}}=\frac 1N\int Y_{N,0}(\d k)\, |k|\1\{|k|\leq L\}$ and analogously defined $D_N^{\ssup{>L}}$ and $D_N=D_N^{\ssup{\leq L}}+D_N^{\ssup{>L}}$. 

We now decompose the  PPP into $Y_N=Y_N^{\ssup{\leq L}}+Y_N^{\ssup{>L}}$, where $Y_N^{\ssup{\leq L}}$ and $Y_N^{\ssup{>L}}$ are the restrictions to  $\Gamma_{T,\leq L}^{\ssup{1}}$ and to $\Gamma_T^{\ssup 1}\setminus \Gamma_{T,\leq L}^{\ssup{1}}$, respectively. Note that they are independent PPPs with intensity  measures $N M_{b\mu}^{\ssup {T,\leq L}}$ and $N M_{b\mu}^{\ssup {T,> L}}$, respectively.

We drop, in the first term in the exponent, the sum involving all $\Xi_i$'s with $|\Xi_{i,0}|>L$, and obtain that the right-hand side of the last display is not larger than
\begin{equation}\label{LOM2}
\begin{aligned}
  {\tt E}_{N M_{b\mu}^{\ssup {T,\leq L}}}&\Big[\e^{N \phi_b(\frac{1}{N}Y_N^\ssup{\leq L})} \e^{ -\int\varphi_{\xi_0}(\xi)\,\frac 1N Y_N^{\ssup{\leq L}} (\d\xi)} D_N^{\ssup{ \leq L}}\Big]\\
&\qquad\times {\tt E}_{N M_{b\mu}^{\ssup {T,> L}}}\Big[b^{-N{D_N^{\ssup{> L}}}} \e^{(b-1)|Y^{\ssup{> L}}_{N,0}|}\1\{D_N^{\ssup{>L}}\leq \eps\}\Big],
\end{aligned}
\end{equation}
where $\phi_b$ is defined in \eqref{phiNdef}.
It is clear that the last line is $\leq \e^{N O(\eps)}$, since we can assume $b<1$ without loss of generality. 
Define 
\begin{eqnarray}
\widehat \P_N^{\ssup{\leq L}}(\d \nu) &=& (\widehat Z_N^{\ssup{\leq L}})^{-1} \, \e^{N\phi_b(\nu)}\e^{ -\int\varphi_{\xi_0}(\xi)\, \nu^{\ssup{\leq L}} (\d\xi)}\, {\tt P}_{N M_{b\mu}^{\ssup {T,\leq L}}}(\smfrac{1}{N}Y_N^{\ssup{\leq L}} \in \d \nu), \\
\widehat Z_N^{\ssup{\leq L}} &=& {\tt E}_{N M_{b\mu}^{\ssup {T,\leq L}}}\Big[\e^{N \phi_b(\frac{1}{N}Y_N^\ssup{\leq L})} \e^{ -\int\varphi_{\xi_0}(\xi)\,\frac 1N Y_N^{\ssup{\leq L}} (\d\xi)} \Big].
\end{eqnarray}

Then \eqref{LOM2} is not larger than $\e^{N O(\eps)}\widehat Z_N^{\ssup{\leq L}} \widehat \E_N^{\ssup{\leq L}}(D_N^{\ssup {\leq L}})$. Now, in the same was as in the proof of the LDP from Theorem~\ref{thm-LDP}(2), we can derive the LDP for $\frac{1}{N}Y_N^{\ssup{\leq L}}$ under $\widehat \P_N^{\ssup{\leq L}}$, using Lemma~\ref{lem-YLDP} and Varadhan's lemm). From this, we first get that $\widehat Z_N^{\ssup{\leq L}}  \leq \e^{-N(\kappa^{\ssup{\leq L}}+o(1))}\leq \e^{-N \kappa/3}$ for all sufficiently large $N$. Furthermore, analogously to Corollary~\ref{cor-accu}, the sequence of distributions of $\frac{1}{N}Y_N^{\ssup{\leq L}}$  under $\widehat \P_N^{\ssup{\leq L}}$ has accumulation points, and each one lies in the set of minimisers of $I_\mu^{\ssup {T,\leq L}}$. Lemma~\ref{lem_ELequations} says that all these accumulation points $\nu^{\ssup{T,\leq L}}$ satisfy the Euler--Lagrange equation in \eqref{ELmitL}, and Lemma~\ref{lem-estinuT}(2) says that they all satisfy $\int \nu^{\ssup {T,\leq L}}(\d k)\,|k|\leq C_T$ for some constant $C_T$ that depends only on $T$ and the constants $H$ from \eqref{AssK1} and $h$ from \eqref{AssK2}. $D^{\ssup{\leq L}}$ being a continuous function of $\Vcal_N^\ssup{T}$ has accumulations points too under $\widehat \P_N^{\ssup{\leq L}}$ and their expectations are bounded by $C_T$ as well.

Hence, we have shown that, for all sufficiently large $N$,
$$
\begin{aligned}
\E_{\Poi_{N\mu}}^\ssup{N}\Big[\big|c^{\ssup{\leq L}}_{\Vcal_N^\ssup{T}}\big|\1\big\{|c^{\ssup{> L}}_{\Vcal_{N,0}^{\ssup {T}}}|\leq \eps\}\Big]&
\leq \eqref{LOM2}\\
&\leq \e^{NO(\eps)}\widehat Z_N^{\ssup{\leq L}} \widehat \E_N^{\ssup{\leq L}}(D_N^{\ssup {\leq L}})
\leq \e^{N[O(\eps)-\kappa/3]}(C_T+o(1)).
\end{aligned}
$$
Now we pick $\eps>0$ so small that the exponent on the right-hand side is strictly negative. This implies, for any sufficiently large $N$, that
\begin{equation}
\begin{aligned}
\E_{\Poi_{N\mu}}^\ssup{N}\Big[\big|c^{\ssup{\leq L}}_{\Vcal_N^\ssup{T}}\big|\Big] 
&\leq 1-\eps+o(1)+
\E_{\Poi_{N\mu}}^\ssup{N}\Big[\big|c^{\ssup{\leq L}}_{\Vcal_N^\ssup{T}}\big|\1\big\{|c^{\ssup{> L}}_{\Vcal_{N,0}^{\ssup {T}}}| \leq \eps\}\Big]\\
&=1-\eps+o(1).
\end {aligned}
\end{equation}
This implies that $\G_T^{\ssup {\mu }}<1$ and shows that gelation holds. This finishes the proof of Theorem~\ref{thm-Lossofmass}(2).

\subsection{The Smoluchowski equation}\label{sec-Smolproof}

Now we  prove Lemma~\ref{lem-Smol}. We pick a limit point $\nu^\ssup{T}$ of $(\Vcal_N^{\ssup {T}})_{N\in\N}$ under $\P_{\Poi_{N\mu }}^{\ssup N}$.
Note that, according to Proposition~\ref{prop-criteria}(1)(b) and Theorem~\ref{thm-Lossofmass}(1)(c), $\nu^\ssup{T}$ is uniquely determined as a solution to the Euler--Lagrange equation. We consider $\rho_t^{\ssup {T}}(\nu^{\ssup T})= \int\nu^{\ssup {T}}(\d\xi)\,\xi_t$ for $t\in[0,T]$, where we wrote $\rho^{\ssup T}$ for the map $\rho$ defined in \eqref{rhodef}. First let us argue that $\rho_t^{\ssup {T}}(\nu^{\ssup T})$ does not depend on $T$, as long as $t\leq T<\frac 1H\frac 1{\e^2}\frac\pi{\pi+1}$. Indeed,  by \eqref{Xi_as_function_of_nu}, we have,  for  any $N\in\N$, 
\begin{equation}\label{eq:Vcal_consist}
\rho_t^\ssup{T_1}(\Vcal_N^\ssup{T_1}) = \frac 1N \Xi_t = \rho_t^\ssup{T}(\Vcal_N^\ssup{T}),\qquad 0\leq t\leq T_1<T<\frac 1H\frac 1{\e^2}\frac\pi{\pi+1}.
\end{equation}
We want to pass to the limit $N\to \infty$ and use the continuities of  the maps $\nu\mapsto\rho^\ssup{T}(\nu)$ and  $\nu\mapsto\rho^\ssup{T_1}(\nu)$ on their respective domains (and the fact that also the marginal map $\nu\mapsto\rho_t^\ssup{T}(\nu)$ is continuous), we obtain that $\rho_t^{\ssup {T}}(\nu^{\ssup T}) =\rho_t^{\ssup {T_1}}(\nu^{\ssup {T_1 }})$. To justify the continuity, recall that the EL-equation \eqref{ELpoissonized} for $\nu^\ssup{T}$ imply that $\int \nu^\ssup{T}(\d \xi)\, |\xi_0|^2 \leq \int M_\mu^\ssup{T}(\d \xi)\, |\xi_0|^2$, which is finite by Lemma \ref{lem-refmeasure-finite} and our assumption $TH< \frac{1}{e^2}\frac{\pi}{1+\pi}$. Hence, $\nu^\ssup{T} \in \Acal_{f,\Aold}$ for $f(r)=r^2$ and some $\Aold\in(0,\infty)$. The EL-equations imply that $H(\nu^\ssup T|M_{\mu}^{\ssup T})$ is finite. Thus, Lemma \ref{lem-Contrho} applies and $\rho^{\ssup T}$ is continuous in $\nu^\ssup{T}$. Since $T_1\leq T$, the latter statements are also true for $\nu^\ssup{T_1}$, the (unique) limit point of $\Vcal_N^\ssup{T_1}$ under $\P^\ssup{N}_{\Poi_{N\mu}}$, and hence $\rho^{\ssup{T_1}}$ is also continuous in $\nu^\ssup{T_1}$. Now, equation \eqref{eq:Vcal_consist} implies that $\rho_t^{\ssup {T}}(\nu^{\ssup T}) =\rho_t^{\ssup {T_1}}(\nu^{\ssup {T_1 }})$ for any $t\leq T_1<T$.
 This shows that $\rho_t^{\ssup {T}}(\nu^{\ssup T})= \int\nu^{\ssup {T}}(\d\xi)\,\xi_t$ does not depend on $T$, as long as $t\leq T< \frac 1H\frac 1{\e^2}\frac\pi{\pi+1}$. Therefore, we write from now $\rho_t= \rho_t^{\ssup {T}}(\nu^{\ssup T})$ (in a small abuse of notation).

Our task is to show that, for any $m^*\in \N$ and any bounded continuous test function $g\colon \Scal \to \R$,
\begin{equation}\label{SmolEQtest}
\begin{aligned}
\frac{\d}{\d t} &\int_\Scal\rho_t(\d x^*,m^*)g(x^*)=-
\int_{\Scal}\rho_t(\d x^*,m^*) \,  K \rho_t(x^*,m^*)g(x^*)\\
&+\sum_{\overset{m,m'\in\N\colon}{m+m'=m^*}}\int_\Scal\int_\Scal\int_{\Scal}\rho_t(\d x,m)\rho_t(\d x',m') {\bf K}\big((x,m),(x',m'),\d x^*\big)g(x^*), \qquad  m^*\in\N.
\end {aligned}
\end{equation}

The base of this is the fact that $\nu^{\ssup T}$ satisfies the EL-equation in \eqref{ELpoissonized}. We need to rewrite that equation a bit. Recall that we introduced $\Q_k^{\ssup T}$ in \eqref{Qdef} as
\[
\Q_k^{\ssup T}(\d \xi) =  \P_k(\Xi\in \d\xi, |\Xi_T|=1)\, \e^{\varphi_k(\xi)},
\]
 where we rephrased the event that $\Xi$ lies in $\Gamma^\ssup{1}_{T}$ as the event $\{|\Xi_T|=1\}$ and recall that the density $\varphi_k$ was defined in \eqref{varphidef} as
\begin{align*}
\varphi_k(\xi) = \int_0^T \Phi(\xi_t)\,\d t,\qquad \mbox{where  }\Phi(\phi)=  &\frac{1}{2}\Big[\langle \phi, K\phi\rangle-\langle\phi, K^{\ssup{\rm diag}}\rangle\Big],\qquad \phi\in\Mcal_{\N_0}(\Scal\times \N).
\end{align*}
From \eqref{ROperator} and \eqref{Rdef} and Fubini's theorem, we see that
$$
\mathfrak R^{\ssup T}(\nu^\ssup{T})(\xi)=\int_0^T\d s\, \Big\langle \int_{\Gamma_T^{\ssup 1}}\nu^\ssup{T}(\d \xi')\,\xi_s',K \xi_s\Big\rangle
=\int_0^T\d s\, \langle\rho_s,K \xi_s\rangle.
$$
Hence, we derive from \eqref{ELpoissonized} that
\begin{equation}\label{ELeqrho2}
\rho_t=\e \E_{\Poi_\mu}\Big[\Xi_t \,\1 \{|\Xi_T|=1\}\,\e^{\int_0^T [\Phi(\Xi_s)-\langle \rho_s, K \Xi_s\rangle]\, \d s} \Big],\qquad t\in[0,T].
\end{equation}
 Recall that $\rho_t $ does not depend on $T$, as long as $t\leq T$. Hence the right-hand side of \eqref{ELeqrho2} does not depend on $T$, and we may put $T$ equal to $t$. 
We define the function 
\[
 \Mcal( \Scal\times \N)\ni \phi \mapsto f(\phi) = \int_\Scal \phi(\d x^*,m^*)\, \1\{|\phi|=1\} g(x^*),
\]
then
\begin{equation}\label{ELeqrho3}
\int_\Scal \rho_t(\d x^*,m^*)g(x^*)=\e \E_{\Poi_\mu}\Big[f(\Xi_t) \,\e^{\int_0^t [\Phi(\Xi_s)-\langle \rho_s, K \Xi_s\rangle]\, \d s} \Big],\qquad t\in[0,T].
\end{equation}

We are going to identify the $t$-derivative of both sides. The expectation on the right-hand side is with respect to a Markov chain in continuous time with only finitely many possible Markovian steps, together with an additional execution of a certain expectation (namely, the one with respect to $\Upsilon$) at every elapsure of one of the holding times; and this does not depend on the time. Hence, the right-hand side is differentiable with respect to $t$, as follows from general theory of Markov chains in continuous time on a discrete space, plus the said execution of another expectation that does not depend on time. Furthermore,  the derivate may be identified in terms of the generator $G$ of  the Marcus--Lushnikov process, using a kind of product differentiation rule, as 
\begin{equation}
\begin{aligned}
\frac{\d}{\d t}  \int_\Scal \rho_t(\d x^*,m^*)g(x^*)
=\e \E_{\Poi_\mu}\Big[\big((Gf)(\Xi_t) + f(\Xi_t)(\Phi(\Xi_t) -&\langle \rho_t,K\Xi_t\rangle)\big)\\
&\times \e^{\int_0^t [\Phi(\Xi_s)-\langle \rho_s, K \Xi_s\rangle]\, \d s}\Big].
\end{aligned}
\end{equation}
One easily checks that $f(\Xi_t)\Phi(\Xi_t) = 0$. 

To derive the Smoluchowski equation we will prove the following two equations, the first one dealing with the gain of particles of type $(\d x^*,m^*)$:
\begin{equation}\label{SmolGain}
\begin{aligned}
&\e \E_{\Poi_\mu}\Big[(Gf)(\Xi_t)\e^{\int_0^t [\Phi(\Xi_s)-\langle \rho_s, K \Xi_s\rangle]\, \d s}\Big] \\
&= \sum_{\substack{m, m' \colon \\ m+m' = m^*}} \int_\Scal \int_\Scal \int_\Scal \rho_t(\d x,m) \rho_t(\d x',m') {\bf K}\big((x,m),(x',m'), \d x^*\big) g(x^*),
\end{aligned}
\end{equation}
and the second one dealing with the loss of particles of type $(\d x^*,m^*)$
\begin{equation}\label{SmolLoss}
\e \E_{\Poi_\mu}\Big[f(\Xi_t)\langle \rho_t,K\Xi_t\rangle  \e^{\int_0^t [\Phi(\Xi_s)-\langle \rho_s, K \Xi_s\rangle]\, \d s}\Big] = \int_\Scal \rho_t(\d x^*,m^*)K\rho_t(x^*,m^* )g(x^*).
\end{equation}
The second one is an immediate consequence of the fact that, 
\[
f(\Xi_t) \langle \rho_t, K \Xi_t \rangle = \int_\Scal g(x^*) \Xi_t(\d x^*,m^*) \1\{|\Xi_t|=1\} \langle \Xi_t, K\rho_t\rangle 
\]
Then, equation \eqref{SmolLoss} follows from \eqref{ELeqrho3} by interchanging the integration of $x^*$ and the expectation $\E_{\Poi_\mu}$.

It remains to show equation \eqref{SmolGain}. From  Section~\ref{sec-model} we see that the generator of the Marcus--Lushnikov process may be written as
\begin{equation}\label{Generator}
\begin{aligned}
G(f)(\phi)&=\sum_{\{(x,m),(x',m')\}}\int_\Scal {\bf K}_\phi\big((x,m),(x',m'),\d z\big)\\
&\qquad\qquad\Big[f\big(\phi-\delta_{(x,m)}-\delta_{(x',m')}+\delta_{(z,m+m')}\big)-f(\phi)\Big],\qquad \phi\in\Mcal_{\N_0}(\Scal\times \N).
\end{aligned}
\end{equation}
where we sum over the possible (unordered) pairs  $(x,m),(x',m')\in \supp(\phi)$. Observe that ${\bf K}_\phi f(\phi) = 0$, since ${\bf K}_\phi =0$ ,if $|\phi|=1$, and $f(\phi)=0$, if $|\phi|\neq 1$, and hence
\begin{align*}
(G f) (\phi) = \sum_{\{(x,m),(x',m')\}}\int_\Scal {\bf K}_\phi\big((x,m),(x',m'),\d x^*\big)f\big(\phi-\delta_{(x,m)}-\delta_{(x',m')}+\delta_{(x^*,m+m')}\big)
\end{align*}
In fact $(G f) (\phi)$ is only non-trivial if $\phi=\delta_{(x,m)}+\delta_{(x',m')}$ for some $x,x'\in \Scal$ and $m,m'\in \N$ with $m+m' = m^*$ and in that case
\begin{equation}\label{eq:generator-id}
(G f) (\phi) = \int_\Scal {\bf K}_\phi\big((x,m),(x',m'),\d x^*\big) g(x^*) = \int_\Scal {\bf K}\big((x,m),(x',m'),\d x^*\big) g(x^*),
\end{equation}
where the second equality can be checked by distinguishing the two cases in the definition \eqref{Mkerneldef} of ${\bf K}_\phi$. 
For any $m\in \N$ we will use the short-hand notation $\phi =\delta_{(\cdot,m)}$ to denote that $\phi \in \{\delta_{(x,m)}\colon x\in \Scal\}$.
We have that
\begin{align*}
&\e \E_{\Poi_\mu}\Big[(Gf)(\Xi_t)\,\e^{\int_0^t [\Phi(\Xi_s)-\langle \rho_s, K \Xi_s\rangle]\, \d s}\Big] \\
 &= \sum_{\substack{(m,m') \colon \\ m+m' = m^*}}\e \E_{\Poi_\mu}\Big[\1\{\Xi_t = \delta_{(\cdot,m)}+\delta_{(\cdot,m')}\}(Gf)(\Xi_t)\,\e^{\int_0^t [\Phi(\Xi_s)-\langle \rho_s, K \Xi_s\rangle]\, \d s}\Big]\\
\end{align*}
With the help of Lemma \ref{lem:Z-decomp} we now show that 
\begin{equation}\label{Poimu-decomp}
\begin{aligned}
&\e \E_{\Poi_\mu}\Big[\1\{\Xi_t = \delta_{(\cdot,m)}+\delta_{(\cdot,m')} \}(Gf)(\Xi_t)\,\e^{\int_0^t [\Phi(\Xi_s)-\langle \rho_s, K \Xi_s\rangle]\, \d s}\Big]\\
& = \int_\Scal \int_\Scal \int_\Scal \mathbf{K}\big((x,m),(x',m'),\d x^*\big) \e \E_{\Poi_\mu}\Big[\Xi_t(\d x, m)\,\e^{\int_0^t [\Phi(\Xi_s)-\langle \rho_s, K \Xi_s\rangle]\, \d s}\Big] \\
&\hspace{3cm}\e \E_{\Poi_\mu}\Big[\Xi_t(\d x', m')\,\e^{\int_0^t [\Phi(\Xi_s)-\langle \rho_s, K \Xi_s\rangle]\, \d s}\Big]
\end{aligned}
\end{equation}
Note that the left-hand side is equal to 
\begin{equation}
\frac{1}{(m+m')!} \int \mu^{\otimes (m+m')}(\d \mathbf{x})\, \E_{\mathbf{x}}\Big[\1\{\Xi_t = \delta_{(\cdot,m)}+\delta_{(\cdot,m')}\}(Gf)(\Xi_t)\,\e^{\int_0^t [\Phi(\Xi_s)-\langle \rho_s, K \Xi_s\rangle]\, \d s}\Big]
\end{equation}
where we abbreviated $\Xi = \Xi(Z)$. The last expectation can be written as
\begin{equation}\label{expect-decomp}
\begin{aligned}
\sum_{\substack{A,B\colon A\dot{\cup} B=[m+m'], \\|A|=m,\, |B|=m'}} \E_{\mathbf{x}}\Big[\1\{A \nleftrightarrow B\}\1\{\Xi^\ssup{A}_t = \delta_{(\cdot,m)}\}&\1\{\Xi_t^\ssup{B} = \delta_{(\cdot,m')}\}\\
&\times (Gf)(\Xi_t)\,\e^{\int_0^t [\Phi(\Xi_s)-\langle \rho_s, K \Xi_s\rangle]\, \d s}\Big]
\end{aligned}
\end{equation}
where we used the short-hand notation $\Xi^\ssup{A}= \Xi^\ssup{T,A}= (\Xi_t(Z_A)_{t\in [0,T]}$ and recall that $Z_A$ denotes the subprocess of $Z$ that only deals with particles/sets $C$ with $C\subset A$ (and analogously for $B$). Recall that this decomposition is possible under the event $\{A \nleftrightarrow B\}$ and also implies that $\Xi =\Xi^\ssup{A} + \Xi^\ssup{B}$. Under the event $\{\Xi^\ssup{A}_t = \delta_{(\cdot,m)}, \Xi_t^\ssup{B} = \delta_{(\cdot,m')}\}$ we can use formula \eqref{eq:generator-id} and rewrite the right-hand side to get
\[
(Gf)(\Xi^\ssup{A}_t + \Xi^\ssup{B}_t) = \int_\Scal \int_\Scal \int_\Scal \Xi^\ssup{A}_t(\d x, m)\, \Xi^\ssup{B}_t(\d x',m')\,{\bf K}\big( (x,m),(x',m'),\d x^*\big) g(x^*)  
\]
Also, by basic calculations one gets that 
\begin{align*}
\int_0^t \Phi(\Xi^\ssup{A}_s + \Xi^\ssup{B}_s) \,\d s = \int_0^t \Phi(\Xi^\ssup{A}_s)\, \d s + \int_0^t \Phi(\Xi^\ssup{B}_s) \,\d s + R^\ssup{t}(\Xi^\ssup{A}, \Xi^\ssup{B}).
\end{align*}
For any fixed pair $A,B$, with $A\dot{\cup} B=[m+m']$, $|A|=m$ and $|B|=m'$, we can now apply Lemma \ref{lem:Z-decomp} and get that the expectation in \eqref{expect-decomp} is equal to 
\begin{align*}
&\E_{\mathbf{x}^\ssup{A}}\otimes \E_{\mathbf{x}^\ssup{B}} \Big[ \1\{\Xi^\ssup{A}_t = \delta_{(\cdot,m)}\}\1\{\Xi_t^\ssup{B} = \delta_{(\cdot,m')}\}(Gf)(\Xi^\ssup{A}_t + \Xi^\ssup{B}_t)\\
&\hspace{5cm}\e^{\int_0^t [\Phi(\Xi^\ssup{A}_s)-\langle \rho_s, K \Xi^\ssup{A}_s\rangle]\, \d s} \e^{\int_0^t [\Phi(\Xi^\ssup{B}_s)-\langle \rho_s, K \Xi^\ssup{B}_s\rangle]\, \d s}\Big]\\
& =\int_\Scal \int_\Scal \int_\Scal \mathbf{K}\big( (x,m),(x',m'),\d x^*\big) g(x^*) \\
&\hspace{3cm}\E_{\mathbf{x}^\ssup{A}}\Big[\Xi_t(\d x,m)\1\{\Xi^\ssup{A}_t = \delta_{(\cdot,m)}\}\,\e^{\int_0^t [\Phi(\Xi^\ssup{A}_s)-\langle \rho_s, K \Xi^\ssup{A}_s\rangle]\, \d s}\Big] \\
&\hspace{3cm}\E_{\mathbf{x}^\ssup{B}} \Big[\Xi^\ssup{B}_t(\d x',m')\1\{\Xi_t^\ssup{B} = \delta_{(\cdot,m')}\}\, \e^{\int_0^t [\Phi(\Xi^\ssup{B}_s)-\langle \rho_s, K \Xi^\ssup{B}_s\rangle]\, \d s}\Big],
\end{align*}
where we have written $\mathbf{x}^\ssup{A} = (x_i)_{i\in A}$, $\mathbf{x}^\ssup{B} = (x_i)_{i\in B}$ 
Note that 
\begin{align*}
\frac{1}{m!}&\int \mu^{\otimes m}(\d (\mathbf{x}^\ssup{A}))\E_{\mathbf{x}^\ssup{A}}\Big[\Xi_t(\d x,m)\1\{\Xi^\ssup{A}_t = \delta_{(\cdot,m)}\}\,\e^{\int_0^t [\Phi(\Xi^\ssup{A}_s)-\langle \rho_s, K \Xi^\ssup{A}_s\rangle]\, \d s}\Big] \\
&= \e \E_{\Poi_\mu}\Big[\Xi_t(\d x,m)\1\{\Xi_t = \delta_{(\cdot,m)}\} \e^{\int_0^t [\Phi(\Xi_s)-\langle \rho_s, K \Xi_s\rangle]\, \d s} \Big] = \rho_t(\d x, m )
\end{align*}
where the right-hand side does not depend on $A$ anymore. The same holds for the terms derived from the set $B$. Since the number of sets $A,B$ with $A\dot{\cup} B=[m+m']$, $|A|=m$ and $|B|=m'$ is equal to $\frac{(m+m')!}{m!m'!}$, all factorials cancel and we get equation \eqref{Poimu-decomp}. This finishes the proof.

\begin{funding}
This research has been funded by the
Deutsche Forschungsgemeinschaft (DFG) through grant 
CRC 1114 ``Scaling Cascades in Complex Systems'', 
Project C08, and by grant SPP 2265 ``Random Geometric Systems'', Project P01. LA and WK acknowledge partial support from GNAMPA, Indam, via the program ``Professori Visitatori 2022''.

\end{funding}

\bibliographystyle{imsart-number} 
\bibliography{article}       


\end{document}